\documentclass[11pt,leqno]{amsart}

\usepackage{latexsym,enumitem}
\usepackage{amssymb}
\usepackage[cp850]{inputenc}
\usepackage{epsfig}
\usepackage{psfrag}
\usepackage{amsthm}
\usepackage{amscd}
\usepackage{amsmath}
\usepackage{amsfonts}
\usepackage{graphics,caption}
\usepackage[all]{xy}
\usepackage{etoolbox}
\patchcmd{\quote}{\rightmargin}{\leftmargin 2em \rightmargin}{}{}
 \usepackage{mathtools}

\usepackage[english]{babel}
\usepackage{comment}
\usepackage{color}
\usepackage{hyperref}

\let\phi\varphi
\let\deg\grado

\DeclareMathOperator{\eqdef}{\coloneqq} 
\newcommand{\f}[2]{\frac{#1}{#2}} 

\let\epsilon\varepsilon
\let\subset\subseteq
\newcommand{\be}{\begin{equation*}}
 \newcommand{\ee}{\end{equation*}}
\newcommand{\bpf}{\begin{dimo}}
\newcommand{\epf}{\end{dimo}}
\newcommand{\bdefi}{\begin{defin}}
\newcommand{\edefi}{\end{defin}}
\newcommand{\bthm}{\begin{thm}}
\newcommand{\ethm}{\end{thm}}
\newcommand{\blem}{\begin{lem}}
\newcommand{\elem}{\end{lem}}
\newcommand{\bcor}{\begin{cor}}
\newcommand{\ecor}{\end{cor}}

\newcommand{\bprop}{\begin{prop}}
\newcommand{\eprop}{\end{prop}}
\newcommand{\bese}{\begin{ese}}
\newcommand{\eese}{\end{ese}}
\newcommand{\brem}{\begin{rem}}
\newcommand{\erem}{\end{rem}}
\newcommand{\bpfc}{\begin{dimoclaim}}
\newcommand{\epfc}{\end{dimoclaim}}

\newcommand{\rar}{\rightarrow} 
\newcommand{\Lar}{\,\Longleftarrow\,} 
\newcommand{\Rar}{\,\Longrightarrow\,} 
\newcommand{\imm}{\looparrowright}

\newcommand{\abs}[1]{\left\lvert#1\right\rvert}						
\newcommand{\set}[1]{\left\{#1\right\}}					
	
\newcommand{\quotient}[2]{\left.\raisebox{.1em}{$#1\!$}\middle/\raisebox{-.1em}{$#2$}\right.}

\DeclareMathOperator{\emp}{\varnothing} 
\DeclareMathOperator{\N}{\mathbb N}			
\DeclareMathOperator{\R}{\mathbb R}			
\DeclareMathOperator{\C}{\mathbb C}			
\DeclareMathOperator{\Z}{\mathbb Z}			

\DeclareMathOperator{\im}{Im} 

\DeclareMathOperator{\id}{id} 

\newcommand{\M}{\mathscr M}

\newcommand{\diffeo}{\xrightarrow{{}_{\simeq}}} 
\newcommand{\bH}{\mathbb H^3}
 \newcommand{\hyp}[1]{\quotient{\bH}{#1}}
  \newcommand{\iso}{\overset{\text{iso}}\simeq}
\newcommand{\subgroup}{\leqslant}

\newcommand{\cat}[1]{\textbf{\textsf{#1}}} 

\newenvironment{quot}
{
	\vspace{-0.2cm}
	\vspace{0.2cm}
}

\theoremstyle{definition}
\newtheorem{d1}{Definition}[section] 

\newenvironment{defin}
{
	\begin{quot}
		\begin{d1}
		}
		{\end{d1}
	\end{quot}

}

\theoremstyle{definition}
\newtheorem{r1}[d1]{Remark}

\newenvironment{rem}
{
	\begin{quot}
		\begin{r1}
		}
		{\end{r1}
	\end{quot}
}

\theoremstyle{definition}
\newtheorem{e1}[d1]{Exercise}

\theoremstyle{definition}
\newtheorem{ese1}[d1]{Example}

\newenvironment{ese}
{
	\begin{quot}
		\begin{ese1}
	}
	{	
		\end{ese1}
	\end{quot}
}

\theoremstyle{definition}
\newtheorem{f1}[d1]{Formulation}

\theoremstyle{definition}
\newtheorem{f2}[d1]{Fact}

\theoremstyle{definition}

\theoremstyle{definition}
\newtheorem*{con2}{Question}

\theoremstyle{definition}
\newtheorem{t1}[d1]{Theorem}

\newenvironment{thm}
{
	\begin{quot}
		\begin{t1}}
		{\end{t1}
	\end{quot}
}

\theoremstyle{definition}
\newtheorem*{T1*}{Theorem}

\newenvironment{teor*}
{
	\begin{quot}
		\begin{T1*}}
		{\end{T1*}
	\end{quot}
}

\newenvironment{dimo}
{\begin{proof}[Proof]
	}
	{\end{proof}}

\newenvironment{dimoclaim}{\emph{Proof of Claim:}\;}{\hfill$\square$}

	\theoremstyle{definition}
	\newtheorem{l1}[d1]{Lemma}
	
	\newenvironment{lem}
	{
		\begin{quot}
			\begin{l1}}
			{\end{l1}
		\end{quot}
	}
	\theoremstyle{definition}
	\newtheorem{p1}[d1]{Proposition}
	
	\newenvironment{prop}
	{
		\begin{quot}
			\begin{p1}}
			{\end{p1}
		\end{quot}
	}
	
	\theoremstyle{definition}
	\newtheorem{c1}[d1]{Corollary}
	
	\newenvironment{cor}
	{
		\begin{quot}
			\begin{c1}}
			{\end{c1}
		\end{quot}
	}

\newtheorem{innercustomthm}{Theorem}
\newenvironment{customthm}[1]
  {\renewcommand\theinnercustomthm{#1}\innercustomthm}
  {\endinnercustomthm}

 \newtheorem*{Theorem*}{Theorem}
 \newtheorem*{Proposition*}{Proposition}
 \newtheorem*{Lemma*}{Lemma}

\usepackage{lipsum} 
\usepackage{etoolbox}
\usepackage{mathrsfs}  
\usepackage{subfig}
\usepackage{floatrow}

\makeatletter
\def\subsection{\@startsection{subsection}{2}%
  \z@{.5\linespacing\@plus.7\linespacing}{.5\linespacing}
  {\normalfont\bfseries}}
\makeatother

\makeatletter
\def\section{\@startsection{section}{2}%
  \z@{.5\linespacing\@plus.7\linespacing}{1\linespacing}%
  {\normalfont\bfseries}}
\makeatother

\makeatletter
\def\paragraph{\@startsection{paragraph}{4}%
  \z@\z@{-\fontdimen2\font}%
  {\normalfont\bfseries}}
\makeatother

\renewcommand{\M}{\mathcal M}
\renewcommand{\hat}{\widehat}
\renewcommand{\P}{\mathcal P} 
\renewcommand{\tilde}{\widetilde}
\newcommand{\A}{\mathcal A}

\begin{document}

	 	\title{Hyperbolization of infinite-type 3-manifolds}
	 	\author{Tommaso Cremaschi}
	 	\date{\today}
		\thanks{The author gratefully acknowledges support from the U.S. National Science Foundation grants DMS 1107452, 1107263, 1107367 "RNMS: GEometric structures And Representation varieties" (the GEAR Network) and also from the grant DMS-1564410: Geometric Structures on Higher Teichm\"uller Spaces.}
	 	\maketitle
	 	
		\small 
	 		
		\paragraph*{Abstract:} We study the class $\M^B$ of 3-manifolds $M$ that have a compact exhaustion $M=\cup_{i\in\N} M_i$ satisfying: each $M_i$ is hyperbolizable with incompressible boundary and each component of $\partial M_i$ has genus at most $g= g(M)$. For manifolds in $\M^{B}$ we give necessary and sufficient topological conditions that guarantee the existence of a complete hyperbolic metric. 
\normalsize
		
\begin{center}
	 		\section*{Introduction}
	 	\end{center}
	 	Throughout this paper, $M$ is always an oriented, aspherical 3-manifold. A 3-manifold $M$ is said to be \textit{hyperbolizable} if it is homeomorphic to $\hyp \Gamma$ for $\Gamma\subgroup \text{Isom}(\bH)$ a discrete, torsion free subgroup, in which case $\Gamma$ is isomorphic to $\pi_1(M)$. We say that a 3-manifold $M$ is of \textit{finite type} if $M$ has finitely generated fundamental group. If otherwise, we say that $M$ is of \emph{infinite type}.
		
		A question of interest in low-dimensional topology is whether a manifold $M$ is hyperbolizable and what is the interplay between the geometry and the topology of $M$.\ Necessary and sufficient topological conditions for the existence of a complete hyperbolic metric in the interior of a compact 3-manifold have been known since Thurston's proof that the interior of every atoroidal Haken 3-manifold is hyperbolizable (1982, \cite{Kap2001}).\ The result was a step in Thurston's program on the study of geometric structures on 3-manifolds, known as the Geometrization conjecture, which was later completed by Perelman (2003, \cite{Per2003.3,Per2003.1,Per2003.2}). These results give a topological characterization of compact 3-manifolds admitting complete hyperbolic metrics in their interiors. On the other hand, by the Tameness Theorem (2004, \cite{AG2004,CG2006}) hyperbolic 3-manifolds with finitely generated fundamental group are \textit{tame}, that is they are homeomorphic to the interiors of compact 3-manifolds.\ By combining Geometrization and the Tameness Theorem we obtain a complete topological characterisation of hyperbolizable finite type 3-manifolds. We have that an irreducible finite type 3-manifold $M$ is hyperbolizable if and only if $M$ is the interior of a compact atoroidal 3-manifold $\overline M$ with infinite fundamental group. 
		
		In this work, we are concerned with the study of \textit{infinite-type} 3-manifolds. Some interesting examples of infinite-type 3-manifolds are Whitehead manifolds \cite{Wh1934,WH1935}, which were the first examples of non-tame open 3-manifolds, and Antoine's necklace \cite{An1921}, which is an non-tame complement of a Cantor set in $\mathbb S^3$.

		Geometric structures on infinite-type 3-manifolds are not widely studied. In particular, not much is known about the topology of hyperbolizable infinite-type 3-manifolds. Nevertheless, some interesting examples of such 3-manifolds are known (see \cite{SS2013,BMNS2016,Th1998}). In \cite{Th1998}, they arise as geometric limits of quasi-Fuchsian hyperbolic 3-manifolds. In \cite{BMNS2016}, the authors constructed infinite-type 3-manifolds by gluing together collections of hyperbolic 3-manifolds with bounded combinatorics via complicated pseudo-Anosov maps. An essential element of their proof is the model geometry developed to prove the Ending Lamination Conjecture \cite{Mi2010,BCM2012}. The boundedness comes from gluing together manifolds from a finite list of hyperbolizable 3-manifolds with incompressible boundary. Other examples arise in \cite{SS2013} as gluings of acylindrical hyperbolizable 3-manifolds with incompressible boundary and such that their boundary components have uniformly bounded genus.

		There are certain obvious obstructions to the existence of a complete hyperbolic metric. Indeed, let $M\cong\hyp\Gamma$ be an hyperbolizable 3-manifold, then by \cite{Fr2011} $\Gamma$ has no divisible subgroups (see Definition \ref{divisiblelementdefi}), hence neither does $\pi_1(M)$. Moreover, by definition covering spaces of hyperbolizable manifolds are hyperbolizable as well. We say that a manifold $M$ is \textit{locally hyperbolic} if every covering space $N\twoheadrightarrow M$ with $\pi_1(N)$ finitely generated is hyperbolizable.

Given the known obstructions and inspired by the above examples \cite{SS2013,BMNS2016}, we introduce the class $\M^B$, where $B$ stands for bounded, of 3-manifolds $M$ so that:

\vspace{0.3cm}

		\begin{enumerate}
		\item[(i)] $M$ admits a nested compact exhaustion $\set{M_n}_{n\in\N}$ by hyperbolizable 3-manifolds;
		\item[(ii)] for all $n\in\N$, the submanifold $M_n$ has incompressible boundary in $M$ so that $\pi_1(M_n)$ injects into $\pi_1(M_{n+1})$;
		\item[(iii)] each component $S$ of $\partial M_ n$ has uniformly bounded genus, that is $\text{genus}(S)\leq g= g(M)\in\N$.
		\end{enumerate}
			
We denote by $\M$ the class of 3-manifolds satisfying (i) and (ii). It is natural to address hyperbolization questions in this class since, by (i) and (ii), every $M\in\M$ is locally hyperbolic. Moreover, one can also show that for every manifold $M\in\M^B$ $\pi_1(M)$ does not contain any divisible subgroup (see Remark \ref{catzero}). Therefore, it is meaningful to look for a characterisation of hyperbolizable manifolds in $\M^B$. Since $\M^B$ already contains hyperbolizable 3-manifolds, namely the ones in \cite{BMNS2016,SS2013}, a first question is whether there exists non-hyperbolizable 3-manifolds in $\M^B$. In \cite{C20171} we built an example $M_\infty\in\M^B$ answering the following question of Agol \cite{DHM,Ma2007}:

		\medskip
			 	\begin{con2}[Agol] Is there a 3-dimensional manifold $M$ with no divisible subgroups in $\pi_1(M)$ that is locally hyperbolic but not hyperbolic?
	 	\end{con2}
		
	However, the 3-manifold $M_\infty$ is homotopy equivalent to a complete hyperbolic 3-manifold. In \cite{C2018c} we improved the above example by building a 3-manifold $N\in\M^B$ such that $N$ is not homotopy equivalent to any complete hyperbolic 3-manifold.

The main result of this paper is a complete topological characterisation of hyperbolizable manifolds in $\M^B$. Before stating the result we need to introduce some objects and notation.
				
		For all $M\in\M$, we construct a canonical \emph{maximal} bordified manifold $(\overline M,\partial\overline M)$, see Definition \ref{bordidefin}, where each component of $\partial\overline M$ is a surface, not necessarily of finite type nor closed. To construct $\overline M$, we compactify properly embedded $\pi_1$-injective submanifolds of the form $S\times[0,\infty)$ by adding $\text{int}(S)\times\set \infty $ to $M$. The bordification $\overline M$ only depends on the topology of $M$. Specifically, we have that $\text{int}(\overline M)$ is homeomorphic to $ M$, and that any two maximal bordifications for $M$ are homeomorphic. 
		We then say that an essential annulus $(A,\partial A)\rar (\overline M,\partial\overline M)$ is \emph{doubly peripheral} if both components of $\partial A$ are peripheral in $\partial\overline M$.
		
		Our main result is:

		\begin{customthm}{1} \label{maintheorem} Let $M\in\M^B$. Then, $M$ is homeomorphic to a complete hyperbolic 3-manifold if and only if  the associated maximal bordified manifold $\overline M$ does not admit any doubly peripheral annulus. 	\end{customthm}

	\medskip

		\paragraph{Detailed overview of the paper:} In Section \ref{background} we introduce some notation and recall some properties of 3-manifolds. In Section \ref{motivatingexample} we recall the example constructed in \cite{C20171} as a motivation for the topological constructions in Section \ref{section2}. The latter contains the bulk of the paper and is divided into two main subsections. 
		
		In Section \ref{secbordification} we construct the bordification $\overline M$ of  $M\in\M$. To construct $\overline M$ we compactify a maximal collection of pairwise disjoint properly embedded $\pi_1$-injective submanifolds of the form $S\times [0,\infty)$ by adding the `boundary at infinity' $\text{int}(S)\times\set\infty$. Here, $S$ is a surface with $\chi(S)\leq 0$. Thus, for each such product submanifold, this process adds a copy of $\text{int}(S)$ to the boundary of $M$. Adding the boundary at infinity does not change the topology of the interior, that is, we have a natural homeomorphism $\iota:M\rar\text{int}(\overline M)$. A bordification is a pair $(\overline M,\iota)$. We also require that bordifications $(\overline M,\iota)$ have no disk components in $\partial\overline M$, and that no two boundary components $S_1,S_2$ of $\partial \overline M$ contain cusp neighbourhoods $C_1\subset S_1$ and $C_2\subset S_2$ that, together with an annulus $C$ connecting $\partial C_1$ and $\partial C_2$, co-bound a submanifold of the form $(\mathbb S^1\times I)\times [0,\infty)$. These two conditions are to guarantee the existence of a maximal bordification. Bordifications of manifolds in $\M$ are defined up to the following equivalence relation: given $(\overline M,\iota)$ and $(\overline M',\iota')$, with interiors homeomorphic to $M$, then: $(\overline M,\iota)\sim(\overline M'\iota')$ if
there is a homeomorphism :

$$\psi:(\overline M,\partial\overline M)\diffeo(\overline M',\partial\overline M')$$ 

such that $\psi\circ \iota$ is isotopic to $\iota'$. We denote the collection of equivalence classes of bordification of $M$ by $\cat{Bor}(M)$.

We will show that to each $M\in\M$ we can assign a unique, up to homeomorphism, \emph{maximal bordification} $\overline M$. As a key property, all properly embedded product submanifolds $\P:S\times [0,\infty)\hookrightarrow \overline M$ are properly isotopic into collar neighbourhoods of $\partial\overline M$. Thus, the main result in Section \ref{secbordification} is:

 \begin{customthm}{2}
Let $M\in\M$. Then, there exists a unique bordification $\overline M\in\cat{Bor}(M)$ such that every properly embedded submanifold $S\times [0,\infty)$ in $\overline M$ is properly isotopic into a collar neighbourhood of a subsurface of $\partial \overline M$.
\end{customthm}


		 In Section \ref{sectionjsj}, we construct the \emph{characteristic submanifold} $(N,R)$ of the bordified manifold $(\overline M,\partial\overline M)$. This is a codimension-zero submanifold of $(\overline M,\partial\overline M)$ that, up to homotopy, contains all essential Seifert-fibered submanifolds of $\overline M$. The main result of the section is:
		 
		 \begin{customthm}{4} A maximal bordification $\overline M$ of $M\in \M$ admits a characteristic submanifold $(N,R)$ and any two characteristic submanifolds of $\overline M$ are properly isotopic.
\end{customthm}

		  The characteristic submanifold $(N,R)$ of $(\overline M,\partial\overline M)$ is obtained by studying how the characteristic submanifold $(N_n,R_n)$ of each compact component $( M_n,\partial  M_n)$ (see \cite{Jo1979,JS1978}) change as we go through the exhaustion.\ We construct the characteristic submanifold $(N,R)$ by taking maximal essential submanifolds $(Q_n,S_n)$ of $(N_n,R_n)$ with the property that in $M\setminus \text{int}(Q_n)$ we have a properly embedded submanifold homeomorphic to $S_n\times[0,\infty)$ in which $S_n\times\set 0$ corresponds to $S_n$.\ With this notion we can make sense of the condition in Theorem \ref{maintheorem} by looking at the characteristic submanifold of the maximal bordification $\overline M$: by a \textit{doubly peripheral annulus} $C$ we mean an essential annulus $C$ in $\overline M$ such that both boundary components of $C$ are peripheral in the components of $\partial\overline M$ containing them.

		  In Section \ref{necessarycond} by using the arguments of \cite{C20171} we show one direction of Theorem \ref{maintheorem}:

		 \begin{customthm}{5} If $M\in\mathcal M^B$ is hyperbolizable, then $\overline M$ cannot have a doubly peripheral annulus $C$.
 \end{customthm}
 
Finally, in section \ref{section3} we prove Theorem \ref{maintheorem}. The proof in inspired by ideas developed in \cite{SS2013}. For simplicity, we describe the case in which $\overline M$ is acylindrical. Given $ M_i$, we show that there exists $n_i>i$ such that $ M_i$ is contained in the acylindrical part of $ M_{n_i}$. Then, by choosing hyperbolic structures $\rho_i\colon \pi_1( M_i)\rar\text{ Isom}(\mathbb H^3)$ on all the $ M_i$'s and using the fact that for pared acylindrical finite-type hyperbolic 3-manifold $(X,P)$ the algebraic topology on $AH(X,P)$ is compact, see \cite[Thm 7.1]{Th1986}, we get that the sequences $\set{\rho_j\vert_{\pi_1(M_i)}}_{j\geq n_i}$ have converging subsequences. By a diagonal argument we obtain:

	\begin{customthm}{6} Given a manifold $M\in\M^B$, if the maximal bordification $\overline M$ is acylindrical then there exists a hyperbolic 3-manifold $N$ and a homotopy equivalence $f: M\rar N$.
\end{customthm}
	 
To conclude the proof of Theorem \ref{maintheorem} we show:
		 \begin{customthm}{7}Let $M\in\M^B$ and $\phi: \overline M\rar N$ be a homotopy equivalence with $N$ a complete hyperbolic manifold. If $\overline M$ is acylindrical, then we have a homeomorphism $\psi:M\rar N$ homotopic to $\phi$.
									
\end{customthm}

					\paragraph*{Acknowledgements:}  I would like to thank J.Souto for introducing me to the problem and for his advice without which none of the work would have been possible. I would also like to thank I.Biringer and M.Bridgeman for many helpful discussions, looking at some \textit{unreadable} early drafts and helping me in writing things carefully.  I would also like to thank the University of Rennes I, where most of this work was completed, as well as MSRI for their hospitality.

\newpage

	\tableofcontents

		\section{Background and conventions}\label{background}

				\subsection{Notation and Conventions} We use $\cong$ for homeomorphic, $\simeq$ for homotopic and $\overset{\text{iso}}\simeq$ for properly isotopic. By $S\hookrightarrow M$ we denote an embedding of $S$ into $M$ while $S\imm M$ denotes an immersion. By a proper embedding $(S,\partial S)\hookrightarrow (M,\partial M)$ we mean an embedding of $S$ in $M$ mapping boundary to boundary, we allow $\partial S=\emp$, and such that the preimage of compact sets is compact. All appearing 3-manifolds are assumed to be aspherical and orientable.

				 By $\Sigma_{g,n}$ we denote an orientable surface of genus $g$ with $n$ boundary components. We say that a manifold is closed if it is compact and without boundary. Unless otherwise stated we use $I=[0,1]$ to denote the closed unit interval and generally by $\mathbb A$ we denote an annulus $\mathbb A\eqdef \mathbb S^1\times I$ and we use $\mathbb T^2\eqdef \mathbb S^1\times \mathbb S^1$ for a torus. By $\pi_0(M)$ we denote the set of connected components of $M$.

	 	By \textit{wing} of a solid torus $V\subset (M,\partial M)$ we mean a 3-dimensional essential thickened annulus $\mathbb A\times I$ with $\mathbb A\times \set 0$ in $\partial M$ and $\mathbb A\times\set 1$ in $\partial V$ such that $\mathbb A\times\set 1$ winds at least once along a longitude of $V$. Note that the topological type of the solid torus and its wings is still a solid torus. If a solid torus $V$ has $n$ pairwise disjoint wings $w_1,\dotsc, w_n$ then, they decompose $\partial V$ into $2n$ parallel annuli $A_1,\dotsc ,A_{2n}$ such that the wings $w_i\cong\mathbb  A\times I$ are attached to $A_{2i}$ and every pair of subsequent wings $w_i,w_{i+1}$ is separated in $\partial V$ by the annulus $A_{2i+1}$.

									\begin{center}\begin{figure}[h!]
	 										\centering
	 										\def\svgwidth{350pt}
	 										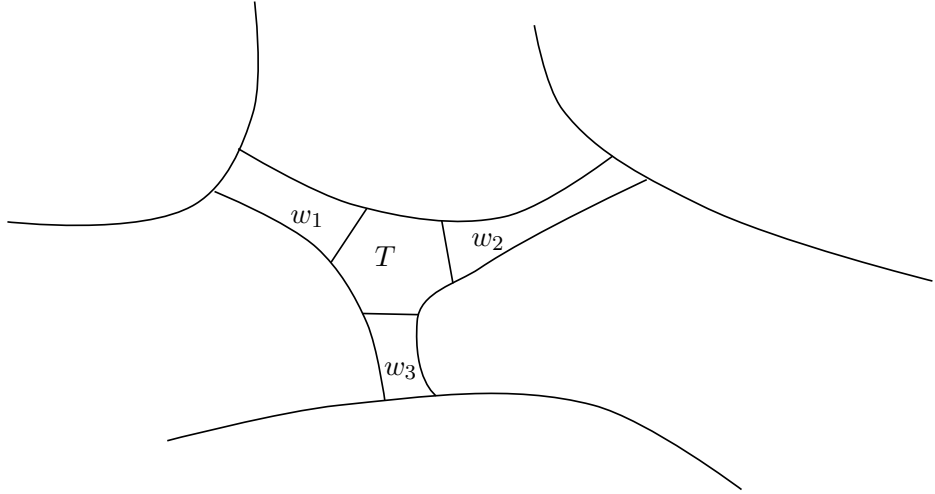
											\caption{Cross section of a solid torus $T$ with three wings $w_1,w_2,w_3$.}
\end{figure}\end{center}			

Let $M$ be an open manifold, by an \emph{exhaustion} $\set{M_i}_{i\in\N}$ we mean a nested collection of compact submanifolds $M_i\subset\text{int}(M_{i+1})$ with $\cup_{i\in\N}M_i=M$. 					
		By \emph{gaps} of an exhaustion $\set{M_i}_{i\in\N}$ we mean the connected components of $\overline {M_i\setminus M_{i-1}}$. Given a manifold with non-empty incompressible boundary $(M,\partial M)$ and a $\pi_1$-injective embedding $\iota:N\hookrightarrow M$ we have a decomposition of $\iota(\partial N)$ into two submanifolds meeting along simple closed curves.\ These two complementary submanifolds are the \emph{outer boundary}: $\partial_{\text{out}} N\eqdef \iota(\partial N)\cap \partial M$ and the closure of the complement: $\partial_{\text{int}} N\eqdef \overline{\iota(\partial N)\setminus \partial_{\text{out}}N}$ which we call the \emph{interior boundary}.
		
		\medskip
						
\subsection{Some 3-manifold topology} 

We now recall some facts and definitions about 3-manifold topology. For more details on the topology of 3-manifolds some references are \cite{He1976,Ha2007,Ja1980}.

Let $M$ be an orientable 3-manifold, then $M$ is said to be \emph{irreducible} if every embedded sphere $\mathbb S^2$ bounds a 3-ball $\mathbb B^3$. Given a connected properly immersed surface $S\imm M$ we say it is \emph{$\pi_1$-injective} if the induced map on the fundamental groups is injective. Furthermore, if $S\hookrightarrow M$ is embedded and $\pi_1$-injective we say that the surface $S$ is \emph{incompressible} in $M$. By the Loop Theorem \cite{He1976,Ja1980} if $S\hookrightarrow M$ is a two-sided surface that is not incompressible we have that there is an embedded disk $D\subset M$ such that $\partial D=D\cap S$ and $\partial D$ is non-trivial in $\pi_1(S)$. Such a disk is called a \emph{compressing} disk. 

An irreducible 3-manifold with boundary $(M,\partial M)$ is said to have \emph{incompressible} \emph{boundary} if every map of a disk: $(\mathbb D^2,\partial\mathbb D^2)\hookrightarrow (M,\partial M)$ is homotopic via maps of pairs into $\partial M$. Therefore, a manifold $(M,\partial M)$ has incompressible boundary if and only if each component $S$ of $\partial M$ is incompressible.

\bdefi We say that an open 3-manifold $M$ is \emph{tame} if it is homeomorphic to the interior of a compact 3-manifold $\overline M$. \edefi

\bdefi Given an irreducible, open 3-manifold $M$ we say that a codimension-zero submanifold $N\overset{\iota}{\hookrightarrow} M$ forms a \textit{Scott core} if the inclusion map $\iota$ is a homotopy equivalence. \edefi 
		By \cite{Sc1973,HS1996,RS1990} if $M$ is an orientable irreducible 3-manifold such that $\pi_1(M)$ is finitely generated we have that a Scott core exists and is unique up to homeomorphism.
		
		\bdefi Given an open 3-manifold $M$ with Scott core $N\hookrightarrow M$ we say that an end $E\subset\overline{ M\setminus N}$ is \textit{tame} if it is homeomorphic to $S\times [0,\infty)$ for $S=\partial E$.\edefi
		
		 A Scott core $C\subset M$ gives us a bijective correspondence between the ends of $M$ and the components of $\partial C$. For the core $C$ we say that a surface $S\in\pi_0(\partial C)$ \textit{faces the end $E$} if $E$ is the component of $\overline{M\setminus C}$ with boundary $S$. It is a known fact that if an end $E$ is exhausted by submanifolds homeomorphic to $S\times I$ then $E$ is a tame end.

\vspace{0.3cm}

Let $M$ be a tame 3-manifold with compactification $\overline M$. Let $C\hookrightarrow M\subset \overline M$ be a Scott core with incompressible boundary then, by Waldhausen's cobordism Theorem \cite[5.1]{Wa1968} every component of $\overline M\setminus \text{int}(C)$ is a product submanifold homeomorphic to $S\times I$ for $S\in\pi_0(\partial C)$. Similarly, we have:

\blem\label{scottcore} Let $\iota:(N,A)\hookrightarrow (M,R)$ be a Scott core for a non-compact irreducible 3-manifold $(M,R)$ that admits a manifold compactification $\overline M$ with $R\subset\partial \overline M$. If $\partial N$ is incompressible in $M$ rel $A$ then $M\cong \text{int} (N)$ and $N\cong \overline M$.
\elem 
\bpf
Consider a component $U$ in $\overline{\overline M\setminus N}$. Then $U$ corresponds to an end of $M$ and since $(N,A)\hookrightarrow (M,R)$ is a homotopy equivalence we have that there exists $S\in \pi_0(\partial N\setminus A)$ facing $U$. Since $\pi_1(S)$ surjects onto $\pi_1(U)$ and $S$ is incompressible in $M$ by Waldhausen's cobordism Theorem \cite[5.1]{Wa1968} we get that $U\cong S\times I$ and the result follows.
\epf

Finally we say that a properly embedded annulus $(\mathbb A,\partial\mathbb  A)$ in a 3-manifold $(M,\partial M)$ is \emph{essential} if $A$ is $\pi_1$-injective and it is not boundary parallel, i.e.\ not isotopic into the boundary. Moreover, a loop $\gamma$ in a surface $(S,\partial S)$ is similarly said \emph{essential} if it is not homotopic into the boundary and $[\gamma]$ is non zero in $\pi_1(S)$.

\subsubsection{JSJ Decomposition}

\bdefi
 A Seifert-fibered 3-manifold $M$ is a compact, orientable, irreducible 3-manifold that has a fibration by circles.
\edefi

\bdefi\label{jsjdefin}
Given a compact 3-manifold $(M,\partial M)$ with incompressible boundary a \emph{characteristic submanifold} for $M$ is a codimension-zero submanifold $(N,R)\hookrightarrow ( M,\partial  M)$ satisfying the following properties:
\begin{enumerate}
\item[(i)] every $(\Sigma,\partial \Sigma)\in\pi_0(N)$ is an essential $I$-bundle or a Seifert-fibered manifold;
\item[(ii)] $\partial N\cap \partial M=R$;
\item[(iii)] all essential maps of a Seifert-fibered manifold $S$ into $ ( M,\partial  M)$ are homotopic as maps of pairs into $(N,R)$;
\item[(iv)] $N$ is minimal, that is no component $P$ of $N$ is homotopic into a component $Q$ of $N$.

\end{enumerate}

\edefi
By work of Johannson \cite{Jo1979} and Jaco-Shalen \cite{JS1978} we have such a submanifold for compact, irreducible 3-manifolds with incompressible boundary:

\begin{Theorem*}[Existence and Uniqueness] Let $(M,\partial M)$ be a compact, irreducible 3-manifold with incompressible boundary. Then there exists a characteristic submanifold $(N,R)\hookrightarrow (M,\partial M)$ and any two characteristic submanifolds are isotopic. 
\end{Theorem*}

This is also called the JSJ or annulus-torus decomposition \cite{Jo1979,JS1978}. Of importance to us will be the fact that if $(F\times I,F\times \partial I)$ is essential and embedded in $(M,\partial M)$ then it is isotopic into the characteristic submanifold $(N,R)$. 

\bdefi

A \emph{window} in a compact irreducible 3-manifold $(M,\partial M)$ with incompressible boundary is an essential $I$-subbundle of the characteristic submanifold.
\edefi

One of the main application of the JSJ decomposition \cite{Jo1979,JS1978} is that:

\bthm
Let $(M,\partial M)$ and $(M',\partial M')$ be compact irreducible 3-manifolds with incompressible boundary and denote by $(N,R)$, $(N',R')$ respectively their characteristic submanifolds. Given a homotopy equivalence $f:M\rar M'$ then we have a homotopy $f\simeq \phi$ such that:
\begin{enumerate}
\item[(i)] $\phi:\overline{M\setminus N}\diffeo \overline{M'\setminus N'}$ is a homeomorphism;
\item[(ii)] $\phi:N\rar N'$ is a homotopy equivalence.
\end{enumerate}
\ethm

In particular if $M$ is acylindrical, that is $M$ has no essential annuli, we have that $M$ has no characteristic submanifold. Therefore, any homotopy equivalence $f:M\rar N$ is homotopic to a homeomorphism.

	\section{Examples and Basic results}\label{motivatingexample}

	 	In this section we show that manifolds in $\M$ are \emph{locally hyperbolic}, that is every cover $N\twoheadrightarrow M$ with $\pi_1(N)$ finitely generated is hyperbolizable.\ For more details on hyperbolic manifolds see \cite{Kap2001,BP1992,MT1998}.

		\blem\label{lochyp}
		
		Let $M$ be orientable and  exhausted by $\set{ M_i}_{i\in\N}$. If each $M_i$ has incompressible boundary in $M$ then $M$ is locally hyperbolic if and only if each $M_i$ is hyperbolizable.
		
		\elem
		\bpf
		$(\Rar)$ Let $M$ be locally hyperbolic and define $N_i$ to be the cover of $M$ corresponding to $\pi_1(M_i)$. Since $\pi_1(M_i)$ is finitely generated and $M$ is locally hyperbolic the cover $N_i$ is hyperbolizable. By the Lifting criterion \cite{Ha2002} for all $i\in\N$ the compact manifolds $M_i$ lift homeomorphically to the cover $N_i$ and the lift has incompressible boundary. Moreover, since all manifolds $M_i, N_i$ are aspherical by Whitehead's Theorem \cite{Ha2002} the lift $\tilde\iota:M_i\hookrightarrow N_i$ forms a Scott core for $N_i$. By the uniqueness of Scott cores \cite{HS1996}, Tameness \cite{AG2004,CG2006} and the fact that $\partial M_i$ is incompressible by applying Lemma \ref{scottcore} we have that $\text{int}(M_i)\cong N_i$ hence  $M_i$ is hyperbolizable.
		
$(\Lar)$ Assume that each $M_k$ is hyperbolizable, since we have that: $\pi_1(M)=\varinjlim_k \pi_1(M_k)$ for every finitely generated $H\leqslant \pi_1(M)$ we can find $i$ such that $\pi_1(M_k)$ contains a generating set of $H$. Therefore, $H$ is a subgroup of $\pi_1( M_k)$ and we denote by $M(H)$ the cover of $M$ corresponding to $H$. Since the cover $M(H)$ factors through the cover $M(k)\eqdef M(\pi_1(M_k))$ we have the following commutative diagram:

		\be		\xymatrix{ M(H) \ar[dr]_{\pi_H}\ar[r]^{\pi'} & M(k)\ar[d]^{\pi}\\ & M }
	 					\ee
		Hence to show that $M(H)$ is hyperbolizable it suffices to show that for all $k\in\N$ the covers $M(k)$ corresponding to $\pi_1(M_k)$ are hyperbolizable.  
		
		Fix $M_k$ and pick a basepoint $x_0\in M_k$ and for all $i\geq k$ consider the covers $\pi_i:(M_i(k),x_i)\rar (M_i,x_0)$ corresponding to $\pi_1(M_k,x_0)$ and denote by $\pi:(M(k),\tilde x_0)\rar (M,x_0)$ the cover of $M$ corresponding to $\pi_1(M_k)$. For all $i\geq k$ the manifolds $M_i(k)$ are hyperbolizable with finitely generated fundamental group, hence by the Tameness Theorem \cite{AG2004,CG2006} they are tame. A core for $M_i(k)$ is given by the homeomorphic lift of $M_k$. Thus, since $\partial M_k$ is incompressible by Lemma \ref{scottcore} all of the $M_i(k)$ are homeomorphic to $\text{int}(M_k)$ and compactify to a manifold homeomorphic to $M_k$. 
		
We now want to show that the cover $M(k)$ is homeomorphic to $\text{int}(M_k)$, hence hyperbolizable. Since $M_k$ lifts to $M(k)$ and induces an isomorphism on $\pi_1$ and the manifolds are aspherical by Whitehead's Theorem \cite{Ha2002} the homemorphic lift $\tilde \iota:(M_k,x_0)\hookrightarrow (M(k),\tilde x_0)$ is a homotopy equivalence, hence $\widetilde M_k\eqdef\tilde\iota(M_k)$ forms a Scott core for $M(k)$.

		  In order to prove that $M(H)\cong \text{int}(M_k)$ it suffices to show that all ends of $M(k)$ are tame. Let $E\in\pi_0\left(\overline{M(k)\setminus \widetilde M_k}\right)$ be an end and denote by $S\in\pi_0\left(\partial \widetilde M_k\right)$ the surface facing $E$. By the Lifting criterion we have  the following commutative diagram:

		\be
		\xymatrix{ (M(k),\tilde x_0)\ar[r]^\pi & (M,x_0) \\
		(M_i(k),\tilde x_0')\ar@{^{(}..>}[u]^{f_i}\ar[r]^{\pi_i} & (M_i,x_0)\ar@{^{(}->}[u]_\iota }
		\ee 
		where the cover: $\pi_i:M_i(k)\rar M_i$ lifts to an embedding $f_i$ into $M(k)$. We denote by $\widetilde M_i(k)$ the images of the $f_i$ and note that their union is $M(k)=\cup_{i\in\N} \widetilde M_i(k)$. Since the $\widetilde M_i(k)$ are tame in each $\overline{\widetilde M_i(k)\setminus \widetilde M_{i-1}(k)}$ we can find a surface $S_i$ homotopic to $S$. Since $\cup_{i\in\N} \tilde M_i(k)$ forms an exhaustion the $\set{S_i}_{i\in\N}$ form a sequence of homotopic embedded surfaces exiting the end $E$ hence by Waldhausen's Cobordism Theorem \cite[5.1]{Wa1968} we have $E\cong S\times [0,\infty)$.
		\epf
		\bdefi
An open 3-manifold $M$ lies in the class $ \M$ if it is irreducible, orientable and satisfies the following properties:
\begin{itemize}
\item[(i)] $M=\cup_{i\in\N} M_i$ where each $M_i$ is a compact, orientable and hyperbolizable 3-manifold;
\item[(ii)] for all $i:\partial M_i$ is incompressible in $M$.
\end{itemize}

Moreover, we say that $M\in\M^B$ if $M\in\M$ and for all $i\in\N$ all components of $\partial M_i$ have genus at most $g=g(M)$
\edefi
	By Lemma \ref{lochyp} and its proof we obtain:
	\bcor
	Given $M=\cup_{i\in\N} M_i\in\M$ then $M$ is locally hyperbolic. Moreover, the cover of $M$ corresponding to $\pi_1(M_i)$ is homeomorphic to $\text{int}(M_i)$.
	\ecor

	\subsection{Locally hyperbolic not Hyperbolic 3-manifolds}\label{lochyp}

	To motivate the construction of the bordification we construct an example, similar to the one constructed in \cite{C20171}, of a non hyperbolizable manifold $N\in\M$. 
	
	We first note that the condition that manifolds in $\M$ are exhausted by hyperbolizable manifolds is a necessary condition. It is not hard to construct manifolds where the gaps $\overline{M_i\setminus M_{i-1}}$ are hyperbolizable but $\cup_{i\in\N} M_i$ has divisible elements (for an example see Appendix \ref{appendix B}). On the other hand, by using our main result it can be shown that manifolds in $\M^B$ do not have divisible elements (see Remark \ref{catzero}).
	
		Consider a compact hyperbolizable 3-manifold $X$ such that $\partial X$ is formed by two incompressible surfaces $S^\pm$ of genus two. Moreover, assume that $X$ has a unique essential cylinder $C$ connecting $S^-$ to $S^+$ along separating curves (see Appendix \ref{appendix A} for the construction of the manifold $X$). 
		
	\begin{center}\begin{figure}[h!]
	 									 	\centering
	 										\def\svgwidth{150pt}
	 										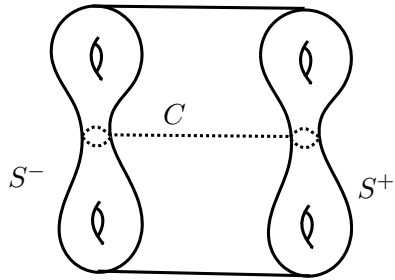\caption{The manifold $X$.}
											\end{figure}\end{center}
	
		We now build the example $N$ of a non-hyperbolizable manifold in $\M$. Let $N_i$, $i\in\Z$, be copies of $X$. By gluing $S_i^+$ to $S_{i+1}^-$ via the identity, so that the boundaries of the cylinders match up we obtain a 3-manifold $N= \quotient{\left(\coprod_{i\in\Z} N_i\right)}{\sim}$ that has a properly embedded cylinder $A=\cup_{i\in\Z}C_i$ with the property that the ends of $A$ are in bijection with the ends of $N$.

		The manifold $N$ has an exhaustion $\set{M_j}_{j\in\N}$ given by taking the portion of the manifold co-bounded by two genus two surfaces $S_{\pm j}^\pm$. Since $\partial M_j$ are all incompressible and each $M_j$ is hyperbolizable we have that $N\in\M^B$. The manifold $N$ looks as:
	
	\begin{center}\begin{figure}[h!]
	
	 										\def\svgwidth{360pt}
	 										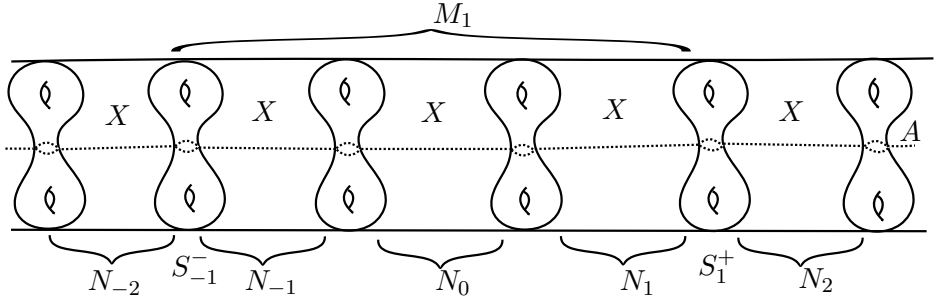\caption{The cylinder $A$ is dashed.}
											\end{figure}\end{center}
	
	The proof that $N$ is not hyperbolizable follows from the techniques and proofs of \cite{C20171}.
	
	\brem
	By \cite{C2018c} one can also show that the manifold $N$ is not homotopy equivalent to any hyperbolic 3-manifold.
	\erem
	
	We now show how to obtain a maximal bordification for $M$. Let $A:\mathbb S^1\times\R\hookrightarrow N$ be the properly embedded bi-infinite cylinder and pick disjoint end's neighbourhoods $A^+,A^-\subset A$. Let $E^+,E^-\subset N$ be open regular neighbourhoods of $A^+,A^-$ respectively. 
										
										By adding two open cylinders at 'infinity' we can bordify $N$ so that in the bordification $\overline N$ the ends $A^\pm$ of the annulus are compactified. A way of achieving this construction is to remove the open neighbourhoods $E^\pm$ from $N$.
										
										 The resulting manifold $\overline N\eqdef N\setminus E^-\coprod E^+$ satisfies: $$\partial\overline N\cong\mathbb  S^1\times (0,1)\coprod\mathbb  S^1\times (0,1)\text{ and } N\cong \text{int}(\overline N)$$
										 
\begin{center}\begin{figure}[h!]
	 										\centering
	 										\def\svgwidth{300pt}
	 										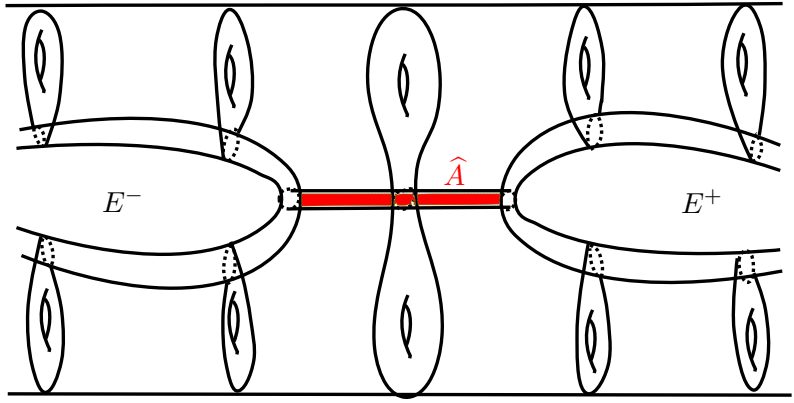
	 										
	 										\caption{The essential cylinder $\hat A$ in the bordification.}\label{figure4}
	 									\end{figure}\end{center}

										  We call such a manifold $\overline N$ a \emph{bordification} for $N$, for a picture see Figure \ref{figure4}. 
										
										This bordification is maximal since the gaps $\overline{M_{i+1}\setminus M_i}$ have a unique cylinde which is isotopic into $A$. Thus, any properly embedded product submanifold $\P:F\times[0,\infty)\hookrightarrow M$ is properly isotopic into $A$. Hence, we cannot compactify any other product to enlarge $\partial\overline N$.

										Note that $\overline N$ has an essential separating annulus $\hat A\eqdef \overline{A\setminus E^\pm}$ that splits $\overline N$ into two infinite type 3-manifolds $\overline N_1,\overline N_2$. Moreover, the manifolds $\overline N_i$, $i=1,2$, are acylindrical since they are obtained by gluing acylindrical manifolds along incompressible hyperbolic subsurfaces of their boundaries.  It seems reasonable to say that the annulus $\hat A$ is the characteristic submanifold of $\overline N$ since if we split along it we get two acylindrical 3-manifolds. The annulus $\hat A$ is an example of a doubly peripheral annulus, so one of the annuli in Theorem \ref{maintheorem}.

			\section{Topological Constructions}\label{section2}

The aim of this section is to study some topological properties of open 3-manifolds $M$ that admit a compact exhaustion $\set{M_k}_{k\in\N}$ by hyperbolizable 3-manifolds with incompressible boundary. We denote this class by $\M$. Moreover, if we assume that for all $k\in\N$ the boundary components $S$ of $\partial M_k$ have uniformly bounded complexity: $\abs{\chi(S)}\leq g\eqdef g(M)$ we say that $M\in\M_g\subset \M$ and we let $\M^B\eqdef\cup_{g\geq2} \M_g$.

Section \ref{secbordification} deals with the problem of constructing canonical bordifications for manifolds $M$ in $\M$. A bordification for $M\in\M$ is a 3-manifold with boundary $\overline M$ such that $M\cong\text{int}(\overline M)$, $\partial\overline M$ has no disk component and with the property that $\overline M$ does not embed in $\overline M'$ where two components $S_1,S_2$ of $\partial\overline M$ are essential subsurfaces of the same component of $\partial\overline M'$. In subsection \ref{maxbordboundedcase} we will show that manifolds in $\M^B$ admit a unique maximal bordification and in \ref{inftypeprod} we extend this to manifolds in $\M$. 

In subsection \ref{sectionjsj}, by studying the characteristic submanifolds $N_i$ of the $ M_i$ and how they change as we go deeper in the exhaustion we construct the characteristic submanifold $N$ for $\overline M$. In the last subsection \ref{acylindricityconditions} we show that even if $\overline M$ is not  acylindrical if $\overline M$ has no doubly peripheral essential cylinders we can find a collection of simple closed curve $ P\hookrightarrow \partial\overline M$ that make $\overline M$ acylindrical relative to $ P$.

		\subsection{Existence of maximal bordifications for manifolds in $\M$}\label{secbordification}

The aim of this section is to show that an open 3-manifold $M$ with a compact exhaustion by hyperbolizable 3-manifolds with incompressible boundary admits a ``maximal" manifold bordification $\overline M$. The boundary components of $\overline M$ are in general open surfaces and come from compactifying properly embedded \emph{product submanifolds} of the form $F\times [0,\infty)$ where $F$ is an incompressible surface. 

We will work in the following class of 3-manifolds:

		\bdefi\label{M} 
An open 3-manifold $M$ lies in the class $ \M$ if it is irreducible, orientable and satisfies the following properties:
\begin{itemize}
\item[(i)] $M=\cup_{i\in\N} M_i$ where each $M_i$ is a compact, orientable and hyperbolizable 3-manifold;
\item[(ii)] for all $i:\partial M_i$ is incompressible in $M$.
\end{itemize}

Moreover, we say that $M\in\M_g$ if $M\in\M$ and for all $i\in\N$ all components of $\partial M_i$ have genus at most $g$. We write $\M^B$ for the class $\cup_{g\geq 2}\M_g$.
\edefi

\bdefi\label{bordidefin}
Given $M\in\M$ we say that a pair $(\overline M,\iota)$, for $\overline M$ a 3-manifold with boundary and $\iota:M\rar \text{int}(\overline M)$ a marking homeomorphism, is a \emph{bordification} for $M$ if the following properties are satisfied:
\begin{enumerate}
\item[(i)] $\partial\overline M$ has no disk components and every component of $\partial\overline M$ is incompressible;
\item[(ii)] there is no properly embedded manifold $$(\mathbb A\times [0,\infty),\partial\mathbb  A\times [0,\infty))\hookrightarrow (\overline M,\partial\overline M)$$ 
\end{enumerate}
Moreover, we say that two bordifications $(\overline M,f),(\overline M',f')$ are equivalent, $(\overline M,f)\sim (\overline M',f')$, if 
we have a homeomorphism $\psi: \overline M\diffeo \overline M'$ that is compatible with the markings, that is: $\psi\vert_{\text{int}(\overline M)}\iso f'\circ f^{-1}$. We denote by $\cat{Bor} (M)$ the set of equivalence classes of bordified manifolds.
\edefi

Condition (ii) is so that $(\overline M,\partial\overline M)$ does not embed into any $(\overline M',\partial\overline M')$ so that two cusps in $\partial \overline M$ are joined by an annulus in $\partial \overline M'$. In particular, this means that we cannot keep adding annuli to $\partial\overline M$ that are isotopic into some other component of $\partial\overline M$. Condition (i) is also so that we can have \emph{maximal} bordification since it is always possible to add disk components to $\partial \overline M$ by compactifying properly embedded rays.

We will build a \emph{maximal} bordification $\overline M_m\in\cat{Bord}(M)$. The bordified manifold $\overline M_m$ has the key property that every properly embedded product submanifold of $M$ is compactified in $ \overline M$ and $M$ is homeomorphic to $\text{int}(\overline M)$. The main result of the section is:

\bthm\label{bordification}
Let $M$ be an orientable, irreducible 3-manifold such that $M=\cup_{i\in\N} M_i$ and $M\in\M$ then, there exists a unique bordification $[\overline M]\in\cat{Bor}(M)$ such that every properly embedded submanifold $F\times [0,\infty)$ in $\overline M$ is properly isotopic into a collar neighbourhood of a subsurface of $\partial \overline M$.
\ethm

We will first deal with the case in which the manifolds $M=\cup_{i\in\N} M_i$ have the property that the genus of $S\in\pi_0(M_i)$ is uniformly bounded, i.e. $M\in\M^B$. Then we will show how to generalize the main technical results to deal with manifolds that have exhaustions with arbitrarily large boundary components.

\subsubsection{Existence of maximal bordification for manifolds in $\M^B$}		\label{maxbordboundedcase}

For open manifolds we define:
		\bdefi\label{productsdefinition}
Given an open 3-manifold $M$ a \emph{product $\P$} is a proper $\pi_1$-injective embedding $\P:F\times [0,\infty)\hookrightarrow M$ for $F$ a, possibly disconnected and of infinite type, surface with no disk components. Given products $\P,\mathcal Q$ then \emph{$\P$ is a subproduct of $\mathcal Q$} if $\P$ is properly isotopic to a restriction of $\mathcal Q$ to a subbundle. Whenever $\P$ is a subproduct of $\mathcal Q$ we write $\P\subset\mathcal Q$. We say that a product $\P:F\times [0,\infty)\hookrightarrow M$ is \emph{simple} if for $\set{F_i}_{i\in\N}$ the connected components of $F$ no $\P_i\eqdef \P\vert_{F_i\times[0,\infty)}$ is a subproduct of a $\P_j$.
 
\edefi 

Note that the image of every level surface of a product $\P$ in $M$ is incompressible in $M$. With an abuse of notation we will often use $\P$ for the image of the embedding.

Given a compact, irreducible and atoroidal 3-manifold $(M,\partial M)$ with incompressible boundary and two essential properly embedded $I$-bundles $P,Q$ we can find a characteristic submanifold \cite{Jo1979,JS1978} $N(P)$ extending $P$, i.e. $P$ is contained in a component of $N(P)$. Then by JSJ Theory \cite{Jo1979,JS1978} we can isotope $Q$ into the characteristic submanifold $N(P)$. Then, up to another isotopy of $Q$ supported in $N(P)$, either $Q$ and $P$ are disjoint or they intersect in one of the following ways:

\begin{enumerate} 
\item[(i)] their union forms a larger connected $I$-bundle in $N(P)$;
\item[(ii)] both $P$ and $Q$ are products over annuli and they intersect `transversally', that is $V\in\pi_0(P\cap Q)$ is a solid torus containing a fiber of both $P$ and $Q$.
\end{enumerate}
In case (ii) we have that $P$ and $Q$ are contained in a essential solid torus component of $N(P)$. We will refer to the second type of intersection as a \emph{cross shape}.

\begin{figure}[h!]

\includegraphics[scale=0.5]{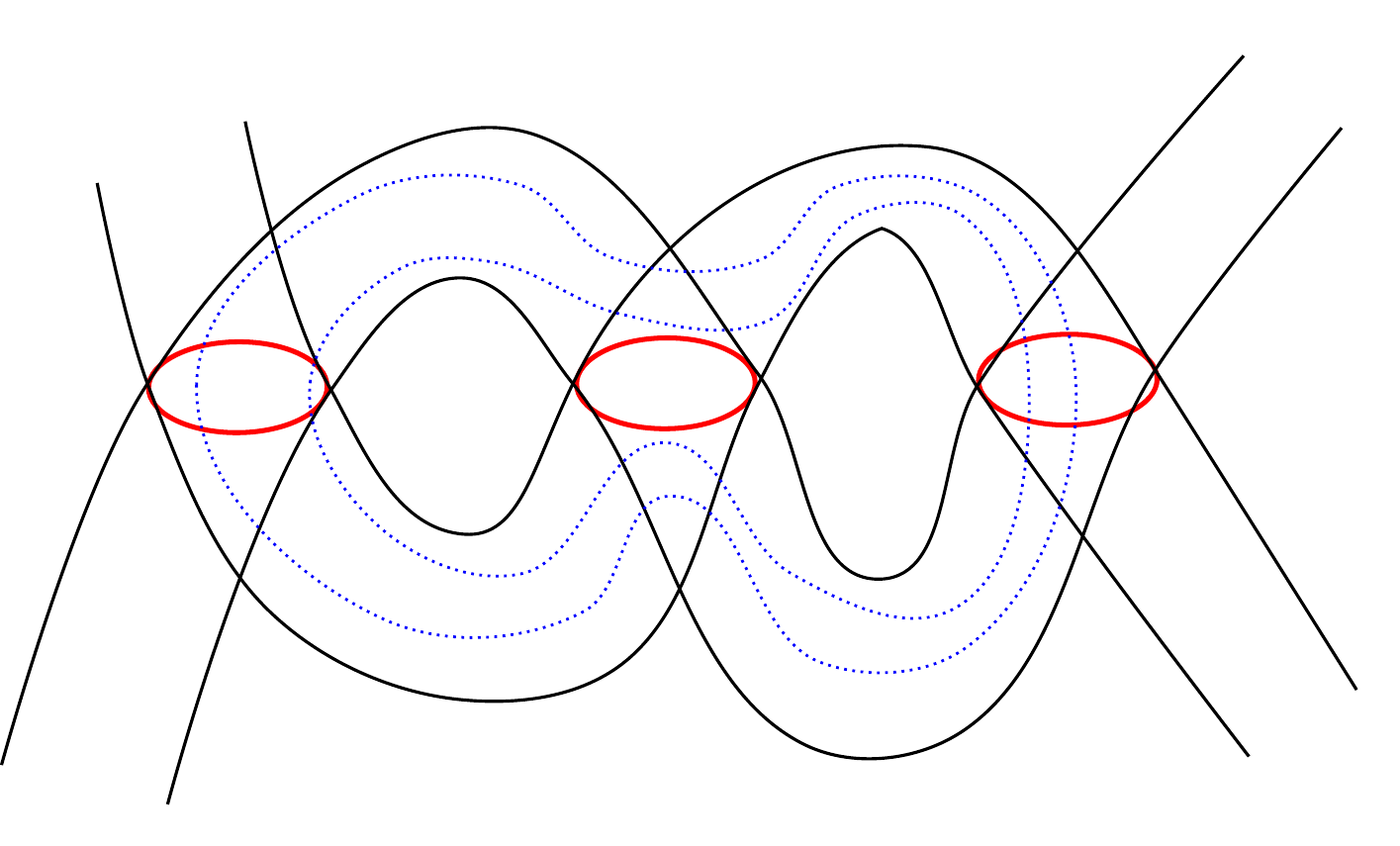}
\caption{The essential torus is dotted in blue.}\label{figure1}
\end{figure}

\blem
Let $P_1,P_2$ be essential properly embedded thickened annuli in a compact irreducible 3-manifold $M$ with incompressible boundary. If, up to isotopy, $P_1,P_2$  intersect in cross shapes more than twice then $M$ is not atoroidal. 
\elem
\bpf
By JSJ theory we can isotope them, relative to the boundary, into the characteristic submanifold.Then they are either disjoint or they intersect in a cross shape. In the latter case, by JSJ theory their union gives a solid torus piece in the JSJ decomposition. However, if the ambient manifold is atoroidal we have that the $P_i$ cannot intersect essentially more than twice. If they intersect at least twice their union contains two essential tori that are joined by an annulus. This configuration of essential tori and annuli (see Figure \ref{figure1}) contradicts the fact that $M$ is atoroidal.  

\hphantom{ada}  \epf

\brem\label{annuliintersection} 

Let $P_1,P_2$ be infinite products over annuli in $M\in\M$ whose component of intersection are cross shapes. If $P_1,P_2$ have at least two components of intersections we can find $k$ such that $ M_k$ contains an embedded copy of Figure \ref{figure1} contradicting the atoroidality of $M_k$. Thus, since $P_1\cap P_2$ intersect in fibers at most once, after a proper isotopy, we can find a compact set $K$ such that outside $K$ the products $P_1,P_2$ are either parallel or have disjoint representatives.

\erem

\bdefi
Let $M\in\M$ and let $X_k\eqdef\overline{M_k\setminus M_{k-1}}$ be the gaps of the exhaustion $\set{M_k}_{k\in\N}$. Given characteristic submanifolds $(N_k,R_k)\subset (X_k,\partial X_k)$, for $R_k\eqdef \partial N_k\cap \partial X_k$, we say that they form a \emph{normal family} $\mathcal N\eqdef\set{N_k}_{k\in\N}$ if: whenever we have essential subsurfaces $S_1, S_2$ of $R_k,R_{k+1}$ respectively that are isotopic in $\partial M_k$ we have an essential subsurface $S\subset R_k\cap R_{k+1}$ such that $S\iso S_1\iso S_2$ in $\partial M_k$. Thus, if $\set{N_k}_{k\in\N}$ forms a normal family we can assume that if a component of $R_{k+1}$ is isotopic into $R_k$ then it is contained in $R_k$, i.e. the $N_k$'s match up along the $\partial M_k$'s.
 \edefi

\blem\label{normalfamily}Given $M=\cup_{k\in\N} M_k\in\M$ there exists a normal family $\mathcal N$ of characteristic submanifolds.
\elem
\bpf
By \cite{Jo1979,JS1978} each $X_k\eqdef \overline{M_k\setminus M_{k-1}}$ has a characteristic submanifolds $N_k$ and define $R_k\eqdef \partial N_k\cap \partial X_k$. Consider $(N_2,R_2)$ and let $S\subset R_2$ be the maximal, up to isotopy, essential subsurface of $R_2$ that is isotopic in $\partial M_1$ into $R_1$. Let $\Sigma$ be a component of $S$. If $\Sigma$ is the boundary of a wing of a solid torus of $N_2$ we can isotope $N_2$ so that $\Sigma\subset R_1$. If $\Sigma$ is contained into an $I$-bundle $P$ we can isotope $P$ so that $P\cap R_1$ contains a subsurface isotopic to $\Sigma$. By doing this for all components of $S$ we obtain that $N_1,N_2$ form a normal pair. We then iterate this construction for all $X_k$ and $N_k$ to obtain the required collection of characteristic submanifolds.
\epf

\bdefi
Given $M\in\M$ and a normal family $\mathcal N=\set{N_k}_{k\in\N}$ of characteristic submanifolds for $X_k\eqdef\overline{M_k\setminus M_{k-1}}$ a product $\P:F\times[0,\infty)\hookrightarrow M$ is in \emph{standard form} if for all $ k\in\N$ every component of $\im( \P)\cap X_k$ is an essential $I$-bundle contained in $N_k$ or an essential sub-surface of $\partial M_k$.
\edefi

\bdefi
We say that a product $\P: F\times [0,\infty)\hookrightarrow M$ is of finite type if the base surface $F$ is compact and of infinite type if the base surface $F$ is of infinite type.
\edefi

In the case that $\P: F\times [0,\infty)\hookrightarrow M$ is of finite type, i.e. $F$ is a compact surface, we have that $\P$ is in standard form if and only if, up to reparametrization, we have $\im(\P)\cap \cup_{k\in\N}\partial M_k=\cup_{i\in\N} \P(F\times\set i)$ and each submanifold $\P(F\times [i,i+1])$ is an essential $I$-bundle contained in $N_{k_i}$ for some $k_i\in\N$.

From now on we will focus on the case of manifolds $M$ that are in $\M^B$, i.e. $M$ is in $\M_g$ for some $g\in\N$. This is so that every product $\P:F\times[0,\infty)\hookrightarrow M$ has the property that every component is a finite type product, see Corollary \ref{boundedgenusboundedproducts}. The next pages will be dedicated to the proof of the following Theorem:

\begin{customthm}{3.18} Consider a product $\P:\Sigma\times [0,\infty)\hookrightarrow M$ with $M=\cup_{k\in\N} M_k\in\M^B$ and $P_i\eqdef\P\vert_{\Sigma_i\times [0,\infty)}$ for $\Sigma_i$ a connected component of $\Sigma$. Given a normal family $\mathcal N=\set{N_k}_{k\in\N} $ of characteristic submanifold for $X_k\eqdef \overline {M_k\setminus M_{k-1}}$ then, there is a proper isotopy $\Psi^t$ of $\P$ such that $\Psi^1$ is in standard form.

\end{customthm}
 The proof is fairly technical and involves ideas and techniques coming from standard minimal position argument. We start by showing that we can properly isotope $\pi_1$-injective submanifolds so that the intersections with the boundaries of the exhaustion are $\pi_1$-injective surfaces.
   
\blem\label{esssubsurface}
Let $\Psi:N\hookrightarrow M\in\M$ be a $\pi_1$-injective proper embedding of an irreducible 3-manifold $N$. Then, there exists a proper isotopy $\Psi^t$ of the embedding $\Psi^0=\Psi$ such that all components of $\mathcal S\eqdef \Psi^1(N)\cap \cup_{k\in\N}\partial M_k$ are $\pi_1$-injective surfaces and no component $S$ of $\mathcal S$ is a disk.\elem
\bpf By a proper isotopy of $\Psi$, supported in $\epsilon$-neighbourhoods of the $\partial M_k$'s we can assume that $\forall k:\partial\im( \Psi)\pitchfork \partial M_k$ so that $\partial M_k\cap\im( \Psi)$ are properly embedded surfaces in $\im(\Psi)$. Thus, we only need to show:

\vspace{0.3cm} 

\paragraph{Claim:} Up to a proper isotopy of $\Psi$ we have that every component of $\mathcal S\eqdef\cup_{k\in\N}\partial M_k\cap \im(\Psi)$ is a properly embedded incompressible surface in $\im(\Psi)$ and $\mathcal S$ has no disk components.

\vspace{0.3cm}

Since every component $S$ of $\mathcal S$ is a subsurface of some $\partial M_k$ and $\partial M_k$ is incompressible in $M$ it suffices to show that up to a proper isotopy of $\Psi$ we have that every component of $\mathcal S$ is an essential subsurface of some $\partial M_k$. Therefore, we have to show that for every component $S$ of $\mathcal S$ we have that $\partial S$ is essential in $M$.

Define $B_k\eqdef \pi_0(\partial M_k\cap \Psi(\partial N))$, since $\Psi$ is proper embedding we have that: for all $k\in\N$ $ \abs{B_{k}}<\infty$. The first step is to show that up to isotopies we can remove all inessential components of $B_k$. To do the isotopies we will need \emph{good balls} for $\partial M_k\cap \partial\im( \Psi)$. These are embedded closed $3$-balls $B\subset\overline {M\setminus M_{k-1}}$ with $\partial B=D_1\cup_\partial D_2$ where $D_1,D_2$ are disks such that $D_1\subset \partial M_k$, $D_2\subset \partial \im( \Psi)$ and $\partial B\cap\partial\im( \Psi)= D_2$. Given a good ball $B$ we can push $\Psi$ through $B$ to reduce $B_k$. Pushing through a $3$-ball effectively adds/deletes a $3$-ball from $\im(\Psi)$. We now define: 
\be\mathcal D_k\eqdef\set{(D_1,D_2)\vert D_1\subset\partial M_k,D_2\subset\partial\im( \Psi)\text{ disks with: }\partial D_1=\partial D_2}\ee
Notice that  by the Loop Theorem \cite{He1976,Ja1980} and incompressibility of $\cup_{k\in\N} \partial M_k$ and  of $\partial\im( \Psi)$ if $\mathcal D_k=\emp$ it means that every component of $B_k$ is essential.

By an iterative argument the key thing to show is that if for all $n< k$ we have that $\mathcal D_n=\emp$, then if $\mathcal D_k\neq\emp$ it contains a good ball for $\partial M_{k}$.

Since $\partial M_k$ and $\partial\im(\Psi)$ are properly embedded in every compact subset we see finitely many components of intersection. Therefore, we can take an innermost component in $\partial M_k$. Thus, we have a disk $D_1\subset \partial M_k$ such that $ D_1\cap \partial \im( \Psi)=\partial D_1$ and the loop $\gamma\eqdef\partial D_1$ is contained in $\partial\im( \Psi)$. Since $\partial \im( \Psi)$ is incompressible we see that $\gamma$ bounds a disk $D_2\subset\partial \im( \Psi)$ and since $D_1$ was picked to be innermost we have that $D_2\cap D_1=\gamma$. By irreducibility of $M$ we have that the embedded 2-sphere $\mathbb S^2\eqdef D_1\cup_\gamma D_2$ bounds a 3-ball $B$. The only thing left to check is that $B\subset \overline{M\setminus M_{k-1}}$. 

The disk $D_2\subset\partial \Psi$ does not intersect any component of $\partial M_n$ with $ n<k$, otherwise by incompressibility of $\partial M_n$ and by taking an innermost disk of intersection we would have $\mathcal D_n\neq\emp$. Hence, $B$ is a good ball for $\partial M_k$.

Thus, we can push $\im(\Psi)$ through the good ball $B$ to reduce $B_k$ without changing any $B_n$ for $n<k$. This process either adds or deletes a 3-ball to $\im(\Psi)$ therefore the homeomorphism type does not change. Moreover, since $B_k$ is finite and every time we remove a good ball it goes down by at least one we have that by pushing through finitely many good 3-balls, i.e. after a proper isotopy $\Psi^t_k$ of $\Psi$ supported in $\overline{M\setminus M_{k-1}}$, we have that every component of $B_k$ is essential. The composition in $k\in\N$ of all the isotopies $\Psi_k^t$ is still proper since the support of $\Psi_k^t$ is contained in $\overline{M\setminus M_{k-1}}$. Thus we obtain a proper isotopy $\Psi^t\eqdef\varinjlim_{k\in\N}\Psi^t_k$ of $\Psi$ such that for all $k\in\N$ $D_k=\emp$. \epf
 
 Before we can prove Theorem \ref{prodstandardform} we will need some technical Lemmas about isotopies of annuli and $I$-bundles in 3-manifolds.

\blem\label{outermost} Let $(M,\partial M)$ be an irreducible 3-manifold with a collection of properly embedded pairwise disjoint boundary parallel annuli $A_1,\dotsc,A_n$ with $\partial A_i\subset \partial M$, for $1\leq i\leq n$. Then there exists pairwise disjoint solid tori $V_1,\dotsc, V_\ell$ in $M$ such that $\cup_{i=1}^n A_i\subset\cup_{k=1}^\ell V_k$ and for all $1\leq k\leq \ell$ we have that $\partial V_k=C_k\cup_\partial A_{i_k}$ for $C_k$'s pairwise disjoint annuli in $\partial M$.
\elem
\bpf Every annulus $A_i$ is properly isotopic rel $\partial A_i$ to an annulus $C_i\subset\partial M$ and for all $1\leq i\leq n$ each pair $A_i,C_i$ co-bounds a solid torus $V_i\subset M$.

\vspace{0.3cm} 

\paragraph{Claim: } For $i\neq j$ the annuli $C_i,C_j\subset\partial M$ are either disjoint or $C_i\subsetneq C_j$.

\vspace{0.3cm} 

\bpfc Since $\partial C_i\cap\partial C_j=\emp$ if $C_i\cap C_j\neq\emp$ we have that at least one component of $\partial C_i\eqdef \alpha_1\cup \alpha_2$ is contained in $\text{int}(C_j)$. Then, if we look at the solid torus $V_j$ we see that $V_j\cap A_i\neq \emp$. Thus, either $A_i\subset V_j$, which gives us that $C_i\subsetneq C_j$, or it escapes. In the latter case we have that $\partial V_j\cap \text{int}( A_i)\neq\emp$ but $\partial V_j=C_j\cup_\partial A_j$ and $\text{int}(A_i)\cap C_j=\emp$. Hence, we must have that $A_i\cap A_j\neq\emp$ contradicting the fact that the annuli $A_1,\dotsc, A_n$ were pairwise disjoint in $M$. 
\epfc

By taking maximal pairs $(C_{i_k},V_k)$, $1\leq k\leq\ell$, with respect to inclusions, we get a collection $\mathcal V=\set{V_1,\dotsc,V_\ell}$ of finitely many solid tori that contain all the annuli $A_1,\dotsc, A_n$. 

Moreover, for $1\leq k\leq \ell$ the $V_k$'s are pairwise disjoint. If not we would have two solid tori $V_t,V_h$ with intersecting boundary and again we contradict the fact that the $A_i$'s are pairwise disjoint or the fact that the $V_k$'s were maximal with respect to inclusion. Thus the required collection of solid tori is given by $\mathcal V$. \epf

Say we have $N\subset \text{int}(M)$ and $F\times I\subset \text{int}(M)$ where $N,M$ are irreducible manifolds with incompressible boundary. If $\partial N\pitchfork F\times I$ by applying Lemma \ref{outermost} to $F\times I$ we can remove $\partial$-parallel annuli of $\partial N\cap F\times I$ by pushing $F\times I$ through the solid tori. Therefore, we have:

\bcor\label{pushanninter} 
Given $F\times I\subset \text{int}(M)$, with $N\subset\text{int}(M)$ and $\partial$-parallel annuli $A_1,\dotsc, A_n\subset \partial N\cap F\times I$ then there is an isotopy $\Psi_t$ of $F\times I$ that is the identity outside neighbourhoods of the $V_k$'s such that for all $t\in[0,1]$ $\Psi_t(F\times I)\subset F\times I$ and $A_i\cap\Psi_1(F\times I)=\emp$.
\ecor

\blem\label{pushinesseibundles}
Let $V$ be a solid torus with $\partial V\eqdef C_1\cup_\partial C_2$, for $C_1,C_2$ annuli and let $A_1,\dotsc,A_n$ be properly embedded $\pi_1$-injective annuli in $V$ such that for $1\leq  i\leq n$ $\partial A_i\subset C_1$. Given a properly embedded annulus $S\subset V$ with $\partial S=\partial C_1=\partial C_2$ there exists an isotopy $\Psi_t$ of $V$ that is constant on $\partial V$ such that $\Psi_1(\cup_{i=1}^nA_i)\cap S=\emp$.\elem
\bpf The annulus $S$ splits the solid torus $V$ into two solid tori $V_1$ and $V_2$ such that $\partial V_k=S\cup_\partial C_k$ for $k=1,2$. Since there are finitely many $A_i$'s and they all have boundary in $C_1\subset\partial V$ we can find an annulus $S'\subset V_2$ with $\partial S'=\partial S$ such that all $A_i$ are contained in the component of $V\setminus S'$ containing $C_1$. By pushing $S'$ to $S$ we obtain the required isotopy. \epf

\blem\label{wald}
Let $\P:\Sigma\times [0,\infty)\hookrightarrow M$ be a product, for $M=\cup_{k\in\N}M_k\in\M$ and let  $\P_i: \Sigma_i\times [0,\infty)\rar M$ be the restriction of $\P$ to the connected components of $\Sigma\times [0,\infty)$ in which we assume that the $\Sigma_i$'s are compact and let $A\subset \Sigma$ be the essential subsurface containing all annular components. Assume that every component of $\mathcal S\eqdef\cup_{k\in\N}\partial M_k\cap \im(\P)$ is properly embedded and incompressible and no component $S$ of $\mathcal S$ is a boundary parallel annulus in $\im(\P)\setminus \P(A\times\set 0)$. Then, there exists a proper isotopy $\Psi^t$ of the embedding $\P$ such that for all $t\in[0,1]$ $\Psi^t(\Sigma\times [0,\infty))\subset\im(\P)$ and for every component $S$ of $\cup_{k\in\N}\partial M_k\cap \im(\Psi^1)$ the surface $(\Psi^1)^{-1}(S)$ is a horizontal fiber in some component of $\Sigma\times[0,\infty)$. Moreover, we have that $\Psi^1:\Sigma\times[0,\infty)\hookrightarrow M$ maps $\Sigma\times\set 0$ into $\mathcal S$ and if $\P(\Sigma\times\set 0)\subset \mathcal S$ then the isotopy can be assumed to be constant on $\P(\Sigma\times\set 0)$.
\elem
\bpf Since we will do proper isotopies supported in $\P_i$ with image in $\P_i$ we can work connected component by connected component. Therefore, it suffices to prove the proposition in the case that $\Sigma\times [0,\infty)$ is connected. We define $ \hat{\mathcal S}$ be the collection of surfaces $S\subset \mathcal S$ such that $\P^{-1}(S)\cap\Sigma\times\set 0=\emp$. Since $\P$ is a proper embedding and $\Sigma\times\set 0$ is compact the set $\hat{\mathcal S}$ is not empty.

\vspace{0.3cm} 

\paragraph{Claim 1:} Given a component $S$ of $\hat{\mathcal S}$, then $\P^{-1}(S)$ is isotopic to a horizontal surface in $\Sigma\times[0,\infty)$.

\vspace{0.3cm} 

\bpfc Since $S$ is compact and $\P$ is a proper embedding we have $0<t_1<\infty$ such that $S\subset \P(\Sigma\times [0,t_1])$ and since $\P^{-1}(S)\cap \Sigma\times\set 0=\emp$ we have $\partial S\subset\P( \partial\Sigma\times (0,t_1])$. By \cite[3.1,3.2]{Wa1968}\footnote{The exact statement we are using here is an easy consequence of the ones cited. In particular we are are applying Waldhausen's result to a an isotopic fiber structure on the $I$-bundle. } we have an isotopy $\phi_t$ of $\P^{-1}(S)$ supported in $\Sigma\times[0,t_1]$ that is the identity on $\partial (\Sigma\times[0,t_1])$ such that the natural projection map: $p:\Sigma\times [0,t_1]\rar \Sigma$ is a homeomorphism on $\phi_1(\P^{-1}(S))$. Moreover, the surface $\P^{-1}(S)$ is also isotopic to a subsurface $S'$ of $\Sigma\times\set {t_1}$. Since $p$ is a homeomorphism on $\P^{-1}(S)$ and $S$ is not a $\partial$-annulus we have that the boundary components of $S$ are in bijection with a subset of boundary components of $\partial\Sigma$. Thus, we have that $S'$  is a clopen subset of $\Sigma$, hence $S$ must be homeomorphic to $\Sigma$.
\epfc

All surfaces $\hat{\mathcal S}$ are pairwise disjoint and isotopic to a fiber. Thus, we can label them by $\set{S_n}_{n\in\N}$ such that for $n<m$ we have that $S_n$ is contained in the bounded component of $\im(\P)\setminus N_\epsilon(S_m)$. 

Consider $S_1$ then we have positive real numbers $a_1,b_1$ with: $0<a_1<b_1<\infty$ such that $\P^{-1}(S_1)\subset \Sigma\times [a_1,b_1]$. By Waldhausen Cobordism Theorem \cite[5.1]{Wa1968} we can change the fibration so that $\P^{-1}(S_1)$ is a horizontal fiber in $ \Sigma_i\times [a_1,b_1]$ hence in $\P$. Then, since $\P^{-1}(S_1)$ is a horizontal fiber in $\P$ by a proper isotopy of $\P$ supported in $\im(\P)$ and with image in $\im(\P)$ we can `raise' $\P(\Sigma\times\set 0)$ to $S_1$ so that $\cup_{k\in\N} \partial M_k\cap \im(\P)=\hat{\mathcal S}$. Note that, these isotopy preserves all properties of $\mathcal S$ and $\hat{\mathcal S}$. Moreover, since the last isotopy removed all components of intersection of $\cup_{k\in\N}\partial M_k\cap\P_i$ that were not in $\hat{\mathcal S}$ we obtain that $\mathcal S=\hat{\mathcal S}$.

Also note that if $\P(\Sigma\times\set0)\subset S$ we don't have to do `raise' isotopy an we automatically have that $\mathcal S=\hat{\mathcal S}$. 

\vspace{0.3cm} 

\paragraph{Claim 2: } Assume that for $1\leq k\leq n$ the surfaces $\P^{-1}(S_k)$ are the horizontal fibers $\Sigma\times \set{k-1}$ in $\Sigma\times[0,\infty)$. Then, by a proper isotopy that is the identity on $\Sigma\times[0,n-1]$ we can make $\P^{-1}(S_{n+1})$ equal to $\Sigma\times\set n$.

\vspace{0.3cm} 

\bpfc
Since $S_{n+1}$ is compact and $\P$ is properly embedded we have $0<n-1<t_n<\infty$ such that the surface $S_{n+1}$ is contained in $\Sigma\times [n-1,t_n]$. Then, by Waldhausen Cobordism Theorem \cite[5.1]{Wa1968} the submanifolds bounded by $S_n,S_{n+1}$ and $\partial \Sigma\times [n-1,t_n]$ is homeomorphic to $\Sigma\times [0,1]$ with $\Sigma\times\set 0=S_n$ and $\Sigma\times \set 1=S_{n+1}$ thus by changing the fiber structure we get that $S_{n+1}$ is also horizontal and equal to $\Sigma\times\set n$. \epfc

By iterating Claim 2 we get that all components of $\mathcal S$ are horizontal in $\P$. Moreover, since $\P(\Sigma\times\set 0)\subset \mathcal S$ we complete the proof. The last statement of the Lemma holds by the observation before Claim 2 and the fact that the isotopies in Claim 2 are constant on $\Sigma\times \set 0$.\epf

From now on we will use \emph{annular product} to indicate a product $\mathcal A: \mathbb A\times [0,\infty)\hookrightarrow M$ where $\mathbb A = \mathbb S^1\times I$ is an annulus.

\bdefi Given a connected product $\P:\Sigma\times[0,\infty)\hookrightarrow M$ such that for every component $S$ of $\cup_{k\in\N}\partial M_k\cap \im(\P)$ the surface $\P^{-1}(S)$ is a horizontal fiber of a component of $\Sigma\times[0,\infty)$ we say that $Q\eqdef \P(\Sigma\times [a,b])$ is a \emph{compact region of $\P$ at $\partial M_k$} if $Q\cap\partial M_k=\P(\Sigma\times\set{a,b})$. Whenever the product and the level is clear we will just write compact region.
\edefi

\bprop\label{inessannuli}
Let $M=\cup_{k\in\N}M_k\in\M$ and $\P:\Sigma\times[0,\infty)\hookrightarrow M$ be a product such that for each component $S$ of $\mathcal S\eqdef \cup_{k\in\N}\partial M_k\cap\im(\P)$ then $\P^{-1}(S)$ is a horizontal surface in some component of $\Sigma\times[0,\infty)$. Consider the subproduct $\mathcal A\subset \P$ consisting of all annular products of $\P$. If for $k< m$ all compact regions of $\mathcal A$ at $\partial M_k$ are essential in either $M_k$ or $\overline{M\setminus M_k}$ then, there is a proper isotopy $\Psi^t_m$ of $\P$ supported in $\overline{M\setminus M_{m-1}}$ such that all compact regions of $\mathcal A$ at $\partial M_m$ are essential.
\eprop
\bpf
Since $\im(\P)$ and $\cup_{k\in\N}\partial M_k$ are properly embedded there are finitely many annular products of $\mathcal A$ that intersect $\partial M_m$. Consider a compact region $Q\eqdef \P(\mathbb A\times [a,b])$ of $\mathcal A$ at $\partial M_m$, since $Q\cap \partial M_m=\P(\mathbb  A\times\set{a,b})$ we have that $Q$ is either contained in $M_m$ or in $\overline {M\setminus M_m}$. Let $\mathscr A_m$ be the collection of all compact regions of $\mathcal A$ at $\partial M_m$ that are boundary parallel in either $M_m$ or $\overline{M\setminus M_m}$. We have that $\abs{\mathscr A_m}<\infty$ and is bounded by $b_m\eqdef \abs{\pi_0(\partial M_m\cap\partial \im(\P))}$ which is finite by properness of the embedding.

\vspace{0.3cm} 

\paragraph{Claim:} Let $\P(\mathbb A\times [a,b])$ be a compact region in $ \mathscr A_n$, for $n\in\N$ such that: 
\be \P(\mathbb A\times I)\cap\cup_{k=1}^n \partial M_k=\P(A\times\set{a,b})\subset\partial M_n\ee
and $\P(A\times[a,b])$ is inessential. Then, there is a solid torus $V\subset \overline{M\setminus M_{n-1}}$ containing $\P(\mathbb A\times [a,b])$ such that all components of $\im(\P)\cap V$ are $\partial$-parallel $I$-bundles contained in $\im(\mathcal A)$.

\vspace{0.3cm} 

\bpfc Consider $\P(\partial\mathbb  A\times [a,b])$ then these are embedded annuli $C_1,C_2$ in either $X_n$ or $\overline {M\setminus M_n}$. If $\P(\mathbb A\times [a,b])$ is $\partial$-parallel so are $C_1,C_2$, hence we have that one of them co-bounds with an annulus $C\subset \partial M_n$ a solid torus $V$ containing $\P(\mathbb A\times[a,b])$. Without loss of generality we can assume that $\partial V=C_1\cup_\partial C$. Since every component of $\cup_{k\in\N} \partial M_k\cap\im (\P)$ is a horizontal fiber in some component of $\im(\P)$ we have that no component of $\partial M_n\cap\im(\P)$ is a boundary parallel annulus or a disk. Since every properly embedded $\pi_1$-injective surface in a solid torus $V$ is either a disk or an annulus we see that $\im(\P)\cap V$ have to be subbundles $Q_1,\dotsc, Q_n$ of annular products in $\im(\P)$. Moreover, all $Q_i$, for $1\leq i\leq n$, are inessential $I$-bundles. \epfc

Let $Q=\P(\mathbb A\times [h,\ell])$ in $\A_k$ be a compact region. Assume that $Q\subset M_m$ and that $Q\cap \partial M_{m-1}\neq\emp$. Since components of intersection of $Q\cap\cup_{k\in\N}\partial M_k$ are horizontal fibers of $Q$ let $\P(\mathbb A\times\set a),\P(\mathbb A\times\set b)$ with $h<a<b<\ell$ be the first and last component of intersections in $Q$ of $Q\cap\partial M_{m-1}$ and let $Q'\eqdef \P( \mathbb A\times [a,b])\subset Q$. Since $\P(A\times [h,a])\subset Q$ and $\P(\mathbb A\times [b,\ell])\subset Q$ have boundaries on distinct components of $X_m\eqdef\overline{M_m\setminus M_{m-1}}$ we get that they are essential $I$-bundles. Therefore, since $Q$  is inessential in $M_m$ we have some $k< m$ such that $Q'\cap (X_k\coprod X_{k+1})$ has a component $T\eqdef \P( \mathbb A\times[t_1,t_2])$ that is inessential and $T$ is a thickened annulus intersecting $\partial M_k$ in $\P(\mathbb A\times\set{t_1,t_2})$. Since $T\subset Q'$ we have that $Q'$ has a compact region that is $\partial$-parallel in either $ M_{k-1}$ or $\overline{M\setminus M_{k-1}}$ contradicting the hypothesis that for all $k<m$ all compact regions of $\mathcal A\cap X_k$ were essential. 

Therefore, we have that $Q$ is boundary parallel in either $ X_m$ or in $\overline{M\setminus M_{m}}$. By the Claim we have a solid torus $V$ such that $V\cap \im(\P)=\im(\mathcal A)\cap V$. By Lemma \ref{pushinesseibundles} we have a proper isotopy of $\P$ supported in a solid torus $N_\epsilon(V)$ contained in $\overline{M\setminus M_{m-1}}$ that removes $Q$ from $\mathcal A_m$ and reduces $b_m$ by at least two.

Thus, we obtain a proper isotopy of $\P$ supported in $\overline{M\setminus M_{m-1}}$ that removes $Q$ from $\mathcal A_m$. Finally, since $\mathcal A_m$ has finitely many elements the composition of these isotopies gives us a proper isotopy $\Psi^t_m$ of $\P$ that makes all sub-bundles of $\mathcal A\cap X_k$ for $k\leq m$ minimal. Moreover, since all the isotopies are supported in $\overline {M\setminus M_{m-1}}$ we get that $\Psi^t_m$ is also supported outside $M_{m-1}$.\epf

The last thing we need to prove Theorem \ref{prodstandardform} is:

\bprop\label{standardisotopy} Let $\P:\Sigma\times [0,\infty)\hookrightarrow M$ be a product with $M=\cup_{k\in\N} M_k\in\M$ and let $A\subset \Sigma$ be the collection of components of $\Sigma $ that are homeomorphic to annuli. Then, there is a proper isotopy of $\P\iso\mathcal Q$ such that all components of $\mathcal S\eqdef\im( \mathcal Q)\cap \cup_{k\in\N} M_k$ are properly embedded $\pi_1$-injective surfaces in $\im(\P)$ such that no $S\in\mathcal S$ is a $\partial$-parallel annulus in $\im(\mathcal Q)\setminus \mathcal Q(A\times\set 0)$ or a disk.

\eprop
\bpf Since products are $\pi_1$-injective by Lemma \ref{esssubsurface} we have that:

\vspace{0.3cm}

\paragraph{Step 1:} Up to a proper isotopy of $\P$ we have that every component of $\mathcal S\eqdef\cup_{k\in\N}\partial M_k\cap \im(\P)$ is a properly embedded incompressible surface in $\im(\P)$ and $\mathcal S$ has no disk components.

\vspace{0.3cm} 

Let $\mathcal A: A\times [0,\infty)\hookrightarrow M$ be the restriction of $\P$ to $A\subset\Sigma$. We first isotope $\mathcal A(A\times\set 0)$ so that every component of $\mathcal A(A\times\set 0)$ is an essential annulus in some $\partial M_k$. Let $\mathcal A_1
\eqdef \mathcal A(A_1\times[0,\infty))$, $A_1\in\pi_0(A)$, be a component of $\im(\mathcal A)$ then not all $\pi_1$-injective annuli $\mathcal S\cap \mathcal A_1$ can have boundary on a component of $\partial \mathcal A_1\setminus\mathcal A (A\times\set 0)$ since otherwise by a proper isotopy of $\mathcal A_1$ supported in $\mathcal A_1$ we would have that $\mathcal A_1\cap \cup_{k\in\N}\partial M_k$ would be compact which contradicts the fact that $\mathcal A_1$ is a proper embedding. Therefore, we must have an essential annulus $S\subset\partial M_k$ of $\mathcal S\cap \mathcal A_1$ whose boundaries are on distinct components of $\partial \mathcal A_1\setminus \mathcal A(A\times\set0)$. Therefore, since $S\cap \mathcal A_1$ and $\mathcal A(A_1\times\set 0)$ are isotopic in $\mathcal A_1$ we can isotope $\mathcal A_1$ so that $\mathcal A(A_1\times\set 0)$ is mapped to $S\subset\partial M_k$. By doing this for all components of $\A$ we can assume that $\mathcal A(A\times\set0)\subset\cup_{k\in\N}\partial M_k$.

\vspace{0.3cm}

\paragraph{Step 2:} Up to a proper isotopy of $\P$ supported in $\im(\P)$ we have that no component $S$ of $\mathcal S$ is a boundary parallel annulus in $\im(\P)\setminus\mathcal A(A\times\set 0)$.

\vspace{0.3cm}

Let $\mathcal A_k$ be the collection of annuli of $\mathcal S_k\eqdef \mathcal S^0\cap\partial M_k$ that are $\partial$-parallel in $\im(\P)\setminus\mathcal A(A\times\set 0)$. Since $\P$ is a proper embedding we have that for all $k\in\N$ $\abs{\pi_0(\mathcal A_k)}<\infty$. By an iterative argument it suffices to show the following:

\vspace{0.3cm} 

\paragraph{Claim 1:} If for $1\leq n< k$ $\mathcal A_n=\emp$ then via an isotopy $\phi_k^t$ of $\P$ supported in $\overline{M\setminus M_{k-1}}\cap\im(\P)$ we can make $\mathcal A_{k}=\emp$.

\vspace{0.3cm} 

\bpfc For all $k\in\N$ we have $0<a_k<b_k<\infty$ such that $\mathcal A_k\subset \P(F_k\times[a_k,b_k])$ for $F_k\subset\Sigma$ a finite collection of connected components of $\Sigma$. 

Denote by $A_1,\dotsc, A_n$ the $\partial$-parallel annuli in $\mathcal A_k$. By applying Corollary \ref{pushanninter} to each component of $\P(F_k\times[a_k,b_k])$ we have a local isotopy $\phi_k^t$ of $\P$ that removes all these intersections. The isotopy $\phi_k^t$ is supported in a collection of solid tori $\mathcal V_k\subset F_k\times[a_k,b_k]$ thus it can be extended to the whole of $\P$. Moreover, if we consider for $n<k$ a component of intersection of $\partial M_n\cap \P(\mathcal V_k)$ then it is either a boundary parallel annulus or a disk. However, we assumed that for $n<k$ $\mathcal A_n=\emp$ and by \textbf{Step 1} no component of $\cup_{k\in\N}\partial M_k\cap\im(\P)$ is a disk thus, the solid tori $\mathcal V_k$ that we push along are contained in $\im(\P)\cap \overline{M\setminus M_{k-1}}$. Therefore, we get a collection of solid tori $\mathcal V_k\subset\im( \P)\cap\overline{M\setminus M_{k-1}}$ such that pushing through them gives us an isotopy $\phi_k^t$ of $\P$ that makes $\mathcal A_k=\emp$. \epfc

Since for all $k\in\N$ $\text{supp}(\phi_k^t)=\mathcal V_k$ is contained in $\overline{M\setminus M_{k-1}}$ the limit $\phi^t$ of the $\phi^t_k$ gives us a proper isotopy of $\P$ such that for all $k\in\N$ $\mathcal A_k=\emp$.

 This concludes the proof of \textbf{Step 2} and the Lemma follows. \epf

We can now show that products whose components are of finite type can be put in standard form.

\bthm\label{prodstandardform} Consider a product $\P:\Sigma\times [0,\infty)\hookrightarrow M$ with $M=\cup_{k\in\N} M_k\in\M$ where the $\Sigma_i$ are the connected component of $\Sigma$ and are compact. Given a normal family $\mathcal N=\set{N_k}_{k\in\N} $ of characteristic submanifold for $X_k\eqdef \overline {M_k\setminus M_{k-1}}$ there is a proper isotopy $\Psi^t$ of $\P$ such that $\Psi^1:\Sigma\times [0,\infty)\hookrightarrow M$ is in standard form.

\ethm
\bpf
From now on we denote the gaps of the exhaustion $\set{M_k}_{k\in\N}$ of $M$ by $X_k\eqdef\overline{M_k\setminus M_{k-1}}$, by definition we have that $N_k\subset X_k$. With an abuse of notation we will often confuse a product $\P$ with its image $\im(\P)$ and we define $\P_i\eqdef\P\vert_{\Sigma_i\times [0,\infty)}$.  By Lemma \ref{standardisotopy} up to a proper isotopy of $\P$ all components of $\mathcal S$ are properly embedded $\pi_1$-injective surfaces in $\im(\P)$ such that no component $S$ of $\mathcal S$ is a $\partial$-parallel annulus or a disk. Then, we are in the setting of Lemma \ref{wald}, thus by a proper isotopy of $\P$ and a reparametrization we have that $\mathcal S=\set{S_i^n}_{i,n\in\N}$ where $S_i^n\eqdef \P(\Sigma_i\times \set n)$ and we let $I^n_i\eqdef \P(\Sigma_i\times [n,n+1])$.

\vspace{0.3cm}

\paragraph{Step 1:} Up to a proper isotopy of $\P$ we can make for all $k\in\N$ all $I$-bundle components $I_i^n$, $i,n\in\N$, of $\im(\P)\cap X_k$ essential.

\vspace{0.3cm}

Since every $I$-bundle $I^n_i$ over a surface $\Sigma_i$ with $\chi(\Sigma_i)<0$ is automatically essential we only need to deal with annular components of $\P$, i.e. products $\P_i:\Sigma_i\times [0,\infty)\hookrightarrow M$ where $\Sigma_i\cong \mathbb A$. We denote by $\mathcal A\subset\P$ the collection of all annular products. By Proposition \ref{inessannuli} and an iterative argument we will show that by isotopies supported in $\overline{M\setminus M_k}$ we can make $\mathcal A\cap X_k$ essential.

Let $\mathscr A_k$ be the $\partial$-parallel compact regions of $\mathcal A$ at $\partial M_k$. Since $\mathcal A\subset \P$ is properly embedded we have that for all $k$ each $\mathscr A_k$ has finitely many components each of which is a compact region over an annulus. By applying Proposition \ref{inessannuli} to $\mathscr A_1\subset\P$ we obtain a proper isotopy $\Psi_1^t$ that makes all compact regions $Q$ in at $\partial M_1$ are essential in either $ M_1$ or $\overline{M\setminus M_1}$. In particular this gives us that every $I$-bundle in $\im(\P)\cap X_1=\im(\P)\cap M_1$ is essential.

We now proceed iteratively. Assume that we made for $1\leq n< k$ all compact regions $Q\in \mathscr A_n$ essential at $\partial M_n$. Then, by applying Proposition \ref{inessannuli} to $\mathscr A_{k}\subset\P$ we get a proper isotopy $\Psi_{k}^t$ supported in $\overline {M\setminus M_{k-1}}$ that makes all compact regions $Q\in \mathscr A_{k}$ essential at $\partial M_{k}$. In particular we get that for all $1\leq n\leq k$ all components of $\im(\P)\cap X_n$ are essential $I$-bundles or essential subsurfaces of $\partial X_n$. 

Since the isotopies $\Psi_k^t$ are supported in $\overline{M\setminus M_{k-1}}$ their composition yields a proper isotopy $\Psi^t\eqdef\varinjlim_{k\in\N}\Psi^t_k$ of $\P$ such that for all $k\in\N$ every component of $\im(\P)\cap X_k$ is an essential $I$-bundle or an essential subsurfaces of $\partial X_n$ given by $\P(\Sigma\times\set0)$.

\vspace{0.3cm}

\paragraph{Step 2:} By a proper isotopy of $\P$ we have that $\im(\P)\subset\cup_{k\in\N} N_k$.  

\vspace{0.3cm}

By Step 1 we have that for all $k\in\N$ the $I$-bundle components of $ \im(\P)\cap X_k$ are essential and pairwise disjoint. Consider $X_1=M_1$ then by JSJ theory we can isotope $\im(\P)\cap X_1$ so that $\im(\P)\cap X_1\subset N_1$. Moreover, since $\P(\Sigma\times (0,\infty))\cap \partial M_1$ is, up to isotopy, contained in both $R_1$ and $R_2$\footnote{We remind the reader that $R_i\eqdef\partial N_i\cap \partial X_i$.} by definition of normal family we can assume that $\P(\Sigma\times (0,\infty))\cap \partial M_1$ is contained in $R_{1,2}\eqdef R_1\cap R_ 2$. This isotopy is supported in a neighbourhood of $X_1$, hence it can be extended to a proper isotopy $\Psi_1^t$ of $\P$. Noting that each component of $\P(\Sigma\times\set0)$ is isotoped at most once to obtain the required proper isotopy it suffices to work iteratively by doing isotopies relative $R_{k,k+1}\eqdef R_k\cap R_ {k+1}$. 

Assume that we isotoped $\P$ such that for all $1\leq n\leq k$ we have that $\P\cap X_n\subset N_n$ and such that $\P(\Sigma\times (0,\infty))\cap \partial M_n$ is contained $R_{n,n+1}$. Since the components of $\im(\P)\cap X_{k+1}$ are essential $I$-bundles of $X_{k+1}$ with some boundary components contained in $R_{k,k+1}$ we can isotope them rel $R_{k,k+1}$ inside $N_{k+1}$ so that their boundaries are contained in $R_{k,k+1}\coprod R_{k+1,k+2}$. This can be extended to an isotopy $\Psi_{k+1}^t$ of $\P$ whose support is contained in $\overline{M\setminus M_k}$, hence the composition of these isotopies gives a proper isotopy of $\P$ such that $\forall k\in\N: \im(\P)\cap X_k\subset N_k$, thus completing the proof. \epf

As a consequence of the Theorem we have:

\bcor\label{boundedgenusboundedproducts}
If $\P:\coprod_{i=1}^\infty F_i\times [0,\infty)\hookrightarrow M$ is a product in $M\in\M_g$ then, every $F_i$ is of finite type and $\abs{\chi(F_i)}$ is uniformly bounded by $2g-2$.
\ecor
\bpf It suffices to show that the statement holds for connected products. Assume that we have a connected product $\P: F\times [0,\infty)\hookrightarrow M$ with $\abs{\chi(F)}\not\leq 2g-2$. Without loss of generality we can assume that: 

$$\abs{\chi(F)}=n>2g-2$$

 since even if $\P$ is not a product of finite type we can find a subproduct $\P_n\eqdef \P\vert_{ F_n\times [0,\infty)}$ where $F_n$ is an essential connected finite type subsurface of $\Sigma$ with $\abs{\chi(F_n)}=n$. 
 
By Theorem \ref{prodstandardform} up to a proper isotopy of $\P$ we can assume $\P$ to be in standard form. Then, the surface $F$ is an essential subsurface of $\Sigma_h\in\pi_0(\partial M_k)$ for some $k\in\N$. Since $h\leq g$ we have that $n=\abs{\chi(F)}\leq 2g-2$, which gives us a contradiction. \epf

Thus, we have that:

\bcor\label{productshavestandardform}
Given a product $\P:\Sigma\times[0,\infty)\hookrightarrow M$ if $M$ is in $\M^B=\cup_{g\geq 2}\M_g$, then there is a proper isotopy $\Psi^t$ of $\P$ such that $\Psi^1$ is in standard form.
\ecor

By isotopying surfaces in general position we have:

\blem\label{matchingfabrition}
Let $\phi_i:(F_i\times I,F_i\times\partial I)\hookrightarrow (F\times I,F\times\partial I)$, $i=1,2$ be essential embeddings in which $\phi_i(F_i\times\set 0)$. Then by a proper isotopy of $\phi_1,\phi_2$ we have that $\im(\phi_1)\cup\im(\phi_2)=\im(\phi_3)$ where $\phi_3:(F_3\times I,F_3\times\partial I)\hookrightarrow (F\times I,F\times\partial I)$ is an essential embedding. Moreover, if $\phi_2(F_2\times\set0)\subset\phi_1( F_1\times\set0)$ we can do the isotopy rel $\phi_1(F_1\times\set0)$.
\elem

\bdefi
Given a normal family of characteristic submanifolds $\mathcal N=\set{N_k}_{k\in\N}$ for $M=\cup_{k\in\N}M_k\in \M$ and a $\pi_1$-injective subbundle $w\eqdef F\times I\hookrightarrow N_k$ with $N_k\subset X_k\eqdef\overline{M_k\setminus M_{k-1}}$ we say that \emph{$w$ goes to infinity} if it can be extended via $I$-bundles $w_i\hookrightarrow N_{k_i}$, $\set{k_i}_{i\in\N}\subset\N$ and $w_0=w$, to a product $F\times [0,\infty)\hookrightarrow M$.
\edefi

Note that each $ F\times I\cong w\subset\pi_0(N_k)$ with $\chi (F)<0$ has at most two extensions to infinity since these $I$-bundles do not branch in any $N_k$. On the other hand annular products can branch off in solid tori in $N_k$ and thus may have infinitely many extensions to infinity.
\begin{rem}\label{locintwind}
Let $w_1,w_2$ be subbundles of $w\subset N_k$ going to infinity. Say that $w\overset\phi\cong F\times I$ and $w_i\overset{\phi_i}\cong F_i\times I$, then if $\phi_1(F_1\times \set 0)\cap \phi_2(F_2\times \set 0)$ is an essential subsurface $F_{1,2}$, $\pi_1$-injective and not $\partial$-parallel, then $w_3\eqdef w_1\cup w_2$ gives a product going to infinity containing the ones given by $w_1,w_2$ as subproducts.
\end{rem}

\bdefi
For $M=\cup_{k\in\N} M_k\in\M$ we say that a product $\P: F\times [0,\infty)\hookrightarrow M$ \emph{starts at $X_k\eqdef \overline{ M_k\setminus M_{k-1}}$} if $\im(\P)\cap X_k$ contains a component homeomorphic to $F\times I$ and $k$ is minimal with respect to this property.
\edefi
Recall that a simple product $\P$ is a product such that no component of $\im(\P)$ is properly isotopic into any other one, see Definition \ref{productsdefinition}.

\blem\label{maxwindows}
Let $M=\cup_{k=1}^\infty M_k\in\M$ and $\mathcal N=\set{N_k}_{k\in\N}$ be a normal family of characteristic submanifolds for $X_k\eqdef\overline{M_k\setminus M_{k-1}}$. Then, for all $k\in\N$ there exists a simple product $\P_k$, in standard form, starting at $X_k$ that contains, up to proper isotopy, all products at $X_k$ generated by sub-bundles over hyperbolic surfaces of windows of $N_k$.
\elem
\begin{proof}
Let $W\subset N_k$ be the collection of $I$-bundles over hyperbolic surfaces of $N_k$. Then, $W$ is homeomorphic, via a map $\phi$, to $F\times I$. If no sub-bundle $F'\times I$ of $F\times I$ goes to infinity there is nothing to do and $\P_k$ is just the empty product.

 Otherwise let $S\times I\subset F\times I$ be a sub-bundle in which $S$ has maximal Euler characteristic and fewest number of boundary components going to infinity through $ F\times \set 1$ such that the product $\mathcal Q:S\times[0,\infty)\hookrightarrow M$ it generates is simple and $\mathcal Q(S\times\set 0)\subset F\times\set 1$. By definition we see that $\mathcal Q$ is also in standard form. We now need to show that $\mathcal Q$ contains all subbundles going to infinity. Let $w'\eqdef \phi(S'\times I)$ be a product going to infinity not properly isotopic into a product given by some components of $\mathcal Q$. Via an isotopy of $S'$ we can assume that $S'$ is in general position with respect to $S$. If $S'\subset S$ we are done. Otherwise since $S'$ and $S$ are in general position no component of $S'\setminus S$ is a disk $D$ such that $\partial D=\alpha\cup_\partial \beta$ with $\alpha\subset \partial  S$ and $\beta\subset \partial S'$. Say we have a disk $D$ component in $S'\setminus S$ then $\partial D$ is decomposed into arcs $\alpha_1,\dotsc,\alpha_{2n}$ such that the odd ones are in $\partial S\cap S'$ and the even ones are in $\partial S$ and $n\geq 2$. Thus, by adding $D$ to $S$ we get that $\chi (S\cup D)=\chi (S) - n+1<\chi (S)$ contradicting the maximality of $\abs{\chi(S)}$. Thus, all components $\Sigma$ of $S'\setminus S$ have $\chi (\Sigma)\leq 0$. 
 
 If we have one $\Sigma\in\pi_0(S'\setminus S)$ that is not an annulus with a boundary component in $\partial S$ and one in $\partial S'$ we would also get that by adding it to $S$ we would get $\abs{\chi(\Sigma\cup_\partial S)}>\abs{\chi(S)}$. Thus, we must have that all components of $S'\setminus S$ are annuli $A$ with one boundary component in $\partial S$ and the other in $\partial S'$ or with both boundary components in $\partial S$. The latter case cannot happen since then by adding $A\times I$ to $S\times I$ we would have gotten a new sub-bundle $\Sigma\times I$ going to infinity through $F\times \set 1$ such that $\chi(\Sigma)=\chi(S)$ but $\abs{\pi_0(\partial\Sigma)}<\abs{\pi_0(\partial S)}$. Hence, $S'$ is isotopic to a subsurface of $S$.
 
 Therefore, we obtain a product $\P_1$ containing all windows going to infinity going through $F\times\set 1$. By doing the same proof for $F\times\set 0$ we obtain another product $\P_0$, hence the required product $\P_k$ is $\P_0\coprod \P_1$. \end{proof}

\bese
For example for the manifold constructed in Section \ref{lochyp} the two annular products start at $X_1= M_1$.
\eese

We will now define maximal products which are the products that we will compactify to construct the maximal bordification.

\bdefi
A simple product $\P: \coprod_{i=1}^\infty \Sigma_i\times[0,\infty)\hookrightarrow M\in\M$ is \emph{maximal} if given any other product $\mathcal Q$ in $M$ then $\mathcal Q$ is properly isotopic to a subproduct of $\P$.\edefi

\bthm\label{evilproposition}
Given $M=\cup_{k\in\N}M_k\in\M^B$ with $\mathcal N$ a normal family of characteristic submanifolds there exists a product in standard form $\P_{max}:F\times [0,\infty)\hookrightarrow M$ such that any other product $\mathcal Q$ is properly isotopic to a sub-product of $\P_{max}$.
\ethm
\bpf
 Let $\mathcal N\eqdef \set{N_k}_{k\in\N}$ be the normal family of characteristic submanifolds for $X_k\eqdef\overline{M_k\setminus M_{k-1}}$. With an abuse of notation in the proof we will often confuse a product with its image. By taking a maximal collection of product starting at $X_i$ we will build collections of pairwise disjoint, disconnected products $P_i$ such that:
 \begin{enumerate}
 \item[(i)] for all $i\in\N:P_i\subset \cup_{k\in\N} N_k$ are essential $I$-bundles;
 \item[(ii)] for all $n\in\N: \cup_{i=1}^n P_i$ is a simple product;
 \item[(iii)] for all $k\in\N: \cup_{i=1}^\infty P_i\cap X_k$ is closed.
 \end{enumerate}
 
 Then by defining $\P_{max}\eqdef \cup_{i=1}^\infty P_i$ we obtain a product that we will show to be maximal by our choice of $P_i$. The fact that $\P_{max}$ is a simple product follows by (ii) and (iii) while (i) gives us that $\P_{max}$ is in standard form. 
 
 With an abuse of notation we will also use $P_i$ to denote the image of the product.
 
 \paragraph{Existence:} Consider $N_1\subset X_1=M_1$. We add to $P_1$ the products coming from Lemma \ref{maxwindows} applied to the windows of $X_1=M_1$. Note that $P_1$ contains finitely many products since all base surfaces are isotopically distinct subsurfaces of $\partial M_1$ and each such submanifold generates at most two non-properly isotopic products. Also note that all such products are necessarily pairwise disjoint and in standard form since they are in every $N_k$, $k\geq 1$.

Next, consider submanifolds of the form $\mathbb A\times I\subset N_1$ that go to infinity. Potentially every such manifold has countably many extensions. If this is the case we choose one representative $A_h^1$, $h\in\N$, for each extension not isotopic into a subproduct of $P_1$ and we add it to $P_1$. However, there is no reason why two such products are not intersecting. So far $P_1$ satisfies (i) and the only obstruction to (ii) is that annular products may intersect in cross shapes inside solid tori components of $N_k$. Moreover, since in each $X_k$ there are finitely many distinct isotopy classes of pairwise disjoint annular products we can choose the representatives $\set{A_h^i}_{h\in\N}$ so that $\forall k:\abs{\set{h\in\N\vert A_h^1\cap M_k\neq\emp}}<\infty$. Therefore, we can also assume that $P_1$ satisfies (iii).

Let $\set{ Q_n}_{n\in\N}\subset P_1$ be all annular subproducts that are not pairwise disjoint. By Remark \ref{annuliintersection} we have that all the intersections of $Q_n$ with $Q_j$ are contained in some compact set (each annular product can intersect another one at most twice). Therefore, by flowing each $Q_j$ in the ``time" direction so that $Q_j$ is disjoint from $Q_n$ with $1\leq n\leq j$ we get a proper isotopy of the $Q_j$'s so that in the image they are pairwise disjoint. Moreover, we can also assume that the new $\set{A_h^1}_{h\in\N}$ still satisfy $\forall k:\abs{\set{h\in\N\vert A_h^1\cap M_k\neq\emp}}<\infty$. By construction $P_1$ satisfies condition (i) and since the annular products do not accumulate it is a product. Moreover, $P_1$ is simple since by construction no annular product is isotopic into a product over a hyperbolic surface and by Lemma \ref{maxwindows} the subproduct of $P_1$ given by products over hyperbolic surfaces is simple. Thus, $P_1$ satisfies conditions (i), (ii) and (iii).

We now proceed inductively. Assume we defined $P_j$, $1\leq j\leq n$, satisfying (i)-(iii) and so that we have representatives of products that start at $X_j$, $k\leq n$, and the annular products $\set{A_h^j}_{h\in\N}$ of $P_j$ intersecting any given $M_k$ are finite.

Consider $X_{n+1}$ and add to $P_{n+1}$, as for $P_1$, the collection of products going to infinity coming from Lemma \ref{maxwindows} applied to $X_{n+1}$ that are not properly isotopic into subproducts of $P_j$, for $1\leq j\leq n$. Every product in $P_{n+1}$ is, by construction, a sub-bundle of $N_k$ for all $k\in\N$. Thus, $P_{n+1}$ satisfies condition (i). We now need to make sure that $P_{n+1}$ satisfies condition (ii), i.e. that $\cup_{i=1}^{n+1} P_i$ is a simple product. Condition (iii) follows from the fact that no product of $P_{n+1}$ intersects $M_n$, otherwise it would have been included in $P_n$. The problem is that the union $\cup_{i=1}^{n+1} P_i$ might not be an embedding, however up to an isotopy of $P_{n+1}$ we can make it an embedding so that $P_1,\dotsc, P_{n+1}$ satisfy (i), (ii) and (iii). 

Let $\mathcal Q\subset P_{n+1}$ be a product over a connected compact surface $F$ with $\chi(F)<0$. For all $k\geq n+1$ $Q\cap N_k$ is a sub-bundle of a, not necessarily connected, window $w_k\in\pi_0( N_k)$. Therefore, $\mathcal Q$ can only intersect products $T\subset\cup_{i=1}^n P_i$ that are also sub-bundles of the same window $w_k$. Let $\set{w_{k_n}}_{n\in\N}\subset\set{w_k}_{k\in\N}$ be the windows containing the intersections of $T$ and $\mathcal Q$. Then, by Lemma \ref{matchingfabrition} we have an isotopy of $Q$ supported in $w_{n_1}$ such that $(\mathcal Q\cup T)\cap w_{n_1}$ is a sub-bundle of $w_{n_1}$. By an iterative argument using Lemma \ref{matchingfabrition} on $w_{n_{k+1}}$ and doing isotopies rel $w_{n_k}\cap w_{n_{k+1}}$ we get a proper isotopy of $\mathcal Q$ such that now $\mathcal Q\cup\cup_{i=1}^n P_i$ is a product. 

By repeating this for the finitely many such $Q$'s in $ P_{n+1}$ we obtain a collection of products $P_{n+1}$ such that all products over surfaces $F$ with $\chi(F)<0$ can be added to $\cup_{i=1}^n P_i$ to define a, possibly disconnected, product $\cup_{i=1}^{n+1} P_i$. Moreover, this product is still simple by Lemma \ref{maxwindows}.

Finally, we add, as for $P_1$, one representative $A^{n+1}_h$, $h\in\N$, for each extension to infinity of annular products starting in $X_{n+1}$ and not properly isotopic into any subproduct of $ P_i$, $1\leq i\leq n+1$. Note that this condition necessarily implies that for $k\leq n:A^{n+1}_h\cap X_k=\emp$, otherwise it would have been added in some $P_k$ with $k\leq n$. Therefore, we can assume that:
\begin{itemize}
\item $\forall k>n:\abs{\set{h\in\N\vert A_h^{n+1}\cap M_k\neq\emp}}<\infty$;
\item $\forall k\leq n:\set{h\in\N\vert A_h^{n+1}\cap M_k\neq\emp}=\emp$
\end{itemize}

For the same reasons as before condition (ii) might still fail, however by doing the same isotopies as for $\set{A^{n+1}_h}_{h\in\N}$ as for $P_1$ so that they become pairwise disjoint and are also disjoint from all annular products $\set{A_h^i}_{i\in\N}$ in $P_i$, $1\leq i\leq n$. Moreover, since no connected subproduct of $P_{n+1}$ intersects $M_n$, otherwise it would have defined a product starting at $X_{n}$ and so it would have been added to $P_n$, we have that:
\be\forall k\leq n: \cup_{i=1}^{n+1} P_i\cap X_k=\cup_{i=1}^{n} P_i\cap X_k \qquad \ee
\noindent and that $\cup_{i=1}^{n+1} P_i\cap X_{n+1}$ is compact, hence it also satisfies (iii).

We then define $\P_{max}\eqdef \cup_{i=1}^\infty P_i$, and $\P_{max}$ satisfies (i) and (ii). Thus, $\P_{max}\eqdef \cup_{i=1}^\infty P_i$ is homeomorphic to $F\times [0,\infty)$ for $F$, in general, some disconnected surface $F=\coprod_{n\in\N} F_n$ where the $F_n$ are all essential subsurface of a fixed genus $g=g(M)$ surface. Property (iii) follows from the previous remark since:
\be\forall k\in\N: \P_{max}\cap X_k=\cup_{i=1}^{\infty} P_i\cap X_{k}=\cup_{i=1}^k P_i\cap X_k\ee
which is compact by (ii). Therefore, by construction $\P$ is a simple a product.

\paragraph{Maximality:} Let $\mathcal Q$ be a product in $M\in\M^B$. Since we are only interested in $\mathcal Q$ up to proper isotopy by Corollary \ref{productshavestandardform} we can assume that it is in standard form with respect to $\mathcal N$. Let $Q_i\cong F_i\times [0,\infty)$ be a connected subproduct of $\mathcal Q$. This means that there is a minimal $k_i$ such that $Q_i\cap X_{k_i}$ is a collection of essential $I$-bundles each one homeomorphic to $F\times I$. Hence it is, up to proper isotopy, contained in a component of $\P_{max}$. Therefore, we get that each connected finite type product $Q_i$ is properly isotopic into $\P_{max}$.

 Let $\P\subset \P_{max}$ be a connected subproduct and let $\mathcal Q_\P$ be all the connected subproducts of $\mathcal Q$ isotopic into subproducts of $\P$. Since $\P$ and $\mathcal Q_\P$ are in standard form, up to a proper isotopy of $Q_\P$ flowing in the `time' direction, they are contained in the same collection of components $\set{ w_n}_{n\in\N}$ for $w_n\subset N_{k_n}$. Then, by doing isotopies in each $w_n$ rel $w_{n-1}$ we can properly isotope $\mathcal Q_\P$ into $\P$. By doing this for all $\P\subset \P_{max}$ we complete the proof. \epf

By the maximality condition we get:

\bcor\label{maxareiso}
If $\P$ and $\mathcal Q$ are both maximal products in $M\in\M^B$ then they are properly isotopic.

\ecor

\bdefi
Given an irreducible 3-manifold $(M,\partial M)$ and a product $\P:F\times [0,\infty)\hookrightarrow M\in\M$ we say that it is $\partial$-\emph{parallel} if $\im(\P)$ is properly isotopic into a collar neighbourhood of a subsurface of $\partial M$. If $\P$ is not $\partial$-parallel we say that it is \emph{essential}. 
\edefi

\bese
Given $(M,\partial M)$ with $S\in \pi_0(\partial M)$ a punctured surface we can build a $\partial$-parallel product $\P$ by taking a collar neighbourhood of a puncture of $S$ and pushing it via a proper isotopy inside $\text{int}(M)$.
\eese

For convenience we recall the definition of a bordification:
\bdefi
Given $M\in\M$ we say that a pair $(\overline M,\iota)$, for $\overline M$ a 3-manifold with boundary and $\iota:M\rar \text{int}(\overline M)$ a marking homeomorphism, is a \emph{bordification} for $M$ if the following properties are satisfied:
\begin{enumerate}
\item[(i)] $\partial\overline M$ has no disk components and every component of $\partial\overline M$ is incompressible;
\item[(ii)] there is no properly embedded manifold $$(\mathbb A\times [0,\infty),\partial\mathbb  A\times [0,\infty))\hookrightarrow (\overline M,\partial\overline M)$$
\end{enumerate}
Moreover, we say that two bordifications $(\overline M,f),(\overline M',f')$ are equivalent $(\overline M,f)\sim (\overline M',f')$ if 
we have a homeomorphism $\psi: \overline M\diffeo \overline M'$ that is compatible with the markings, that is: $\psi\vert_{\text{int}(\overline M)}\iso f'\circ f^{-1}$. We denote by $\cat{Bor} (M)$ the set of equivalence classes of bordified manifolds
\edefi

Condition (ii) is so that $(\overline M,\partial\overline M)$ does not embed into any $(\overline M',\partial\overline M')$ in a way that two cusps in $\partial \overline M$ are joined by an annulus in $\partial \overline M'$. Condition (i) is so that we can have 'maximal' bordification since it is always possible to add disk components to $\partial \overline M$ by compactifying properly embedded rays and so that collar neighbourhoods of $\partial\overline M$ correspond to products in $M$. 
\bdefi
We say that a bordication $[(\overline M,f)]\in\cat{Bor}(M)$ is \emph{maximal} if $\overline M$ has no essential products.\edefi

\blem\label{maprodequivmaxbord}
A bordification $[(\overline M,f)]\in\cat{Bord}(M)$ is maximal if and only if the preimage of a collar of $\partial\overline M$ in $M$ via $f$ is a maximal product.
\elem
\bpf
Let $(\overline M,f)$ be a maximal bordification and let $\partial \overline M= \coprod_{i=1}^\infty S_i$ and $\P_i\eqdef f^{-1} (N_\epsilon(S_i))$. Then, we get a product $\P\eqdef\coprod \P_i$ in $M$. Moreover, $\P$ is simple since otherwise we would have two component $S_1,S_2$ of $\partial \overline M$ that can be joined by a submanifold homeomorphic to $\mathbb A\times [0,\infty)$, contradicting property (ii) of the definition of a bordification. Finally, by property (i) of a bordification we see that $\P$ is $\pi_1$-injective and by maximality of the bordification every product $\mathcal Q$ in $M$ is isotopic into $\P$ and hence $\P$ is a maximal product.

Similarly if for $(\overline M,f)$ a bordification we have that $\P\eqdef f^{-1}(N_\epsilon(\partial\overline M)$ is a maximal product then $\overline M$ is maximal since if not we would have another bordification: $(\overline M',f')$ and an embedding:

$$\psi: (\overline M,\partial\overline M)\hookrightarrow (\overline M',\partial\overline M')$$

such that $\partial\overline M'\setminus \psi(\partial\overline M)$ contains a non-annular component contradicting the maximality of $\P$. \epf

We can now prove the main result of the section:

\bthm\label{bordification}
Let $M\in\M^B$ then there exists a unique maximal bordification $[(\overline M,\iota)]\in\cat{Bor}(M)$.
\ethm

\bpf Since $M\in\M^B$ by Theorem \ref{evilproposition} we have a maximal product $\P_{max}: F\times [0,\infty)\hookrightarrow M$. We now want to compactify $\P_{max}$ by adding $\text{int}(F)\times\set\infty$ to $M$.
Topologically the subproduct $ \P\eqdef\P_{max}\vert_{\text{int}(F)\times [0,\infty)}$ can be naturally compactified to $\overline \P: \text{int}(F)\times [0,\infty]\hookrightarrow M\cup \text{int}(F)\times\set\infty$ by adding the boundary at infinity $\text{int}(F)\times\set{\infty}$ to $M$.

Define $\overline M\eqdef M\coprod \text{int}(F)\times\set\infty$ with the topology that makes $\overline P:\text{int}(F)\times [0,\infty]\hookrightarrow\overline M$ into a homeomorphism onto its image. To see that $\overline M$ is a 3-manifold it suffices to show that $F\times [0,\infty)\cup \text{int}(F)\times\set{\infty}$ is a 3-manifold. This follows from the fact that $F\times [0,\infty)\cup \text{int}(F)\times\set{\infty}$ is naturally an open submanifold of $F\times [0,1]$ and so the smooth structures agree. Moreover, we have that the inclusion: $\id: M\hookrightarrow \overline M$ is an embedding. Since products have no disk components we have that $[(\overline M,\id)]\in\cat{Bord}(M)$ and by Lemma \ref{maprodequivmaxbord} we get that $[(\overline M,\id)]$ is a maximal bordification.

\paragraph{Uniqueness:} Say we have another maximal bordification $[(\overline M',\iota')]$, then by Lemma \ref{maprodequivmaxbord} we obtain a maximal product $\P'$. Since the products $\P$ and $\P'$ are maximal by Corollary \ref{maxareiso} we have a proper isotopy $H_t$ from $\P$ to $\P'$.

We can then extend this proper isotopy to a proper isotopy $\widehat H_t:M\rar M$. The diffeomorphism $H_1:M\rar M$ mapping $\im(\P_m)$ to $\im(\P')$ extends to a diffeomorphism $\psi:\overline M\rar\overline M'$ mapping $\partial \overline M$ to $\partial \overline M'$. By construction we have that this gives an equivalence of bordifications concluding the proof.\epf

\subsubsection{Extension to Manifolds in $\M$}\label{inftypeprod} Note that in Theorem \ref{bordification} we used the fact that the manifold was in $\M^B$ just to say that we had maximal products via Theorem \ref{evilproposition}. Thus,  the aim of this subsection is to show how one can extend Theorem \ref{evilproposition} to deal with infinite type products. To do so it suffices to show that infinite type product can be put in standard form, i.e. extending Theorem \ref{prodstandardform}.

	We will use the term \emph{window} to denote an essential $I$-subbundle of a component of the characteristic submanifold.

\bdefi
We say that an $I$-bundle $F\times I$ embedded in an irreducible 3-manifold $ (M,\partial M)$ with incompressible boundary, not necessarily boundary to boundary, is \emph{mixed} if it is $\pi_1$-injective and it contains a window of $M$.
\edefi
\bese
Let $M$ be a compact, irreducible 3-manifold with incompressible boundary and let $w\overset\phi\cong F\times I$ be a window in $M$ with $\phi(F\times\set i)\subset S_i\in\pi_0(\partial M)$, for $i=0,1$, and $S_0\neq S_1$. If we denote by $N_1$ a collar neighbourhood of $S_1$ we have that $w\cup N_1$ is a mixed $I$-bundle with fiber structure $S_1\times I$.
\eese
We now extend Theorem \ref{prodstandardform}:

	\blem\label{infstandardform}
	Given $\P:\Sigma\times[0,\infty)\hookrightarrow M$ a product in $M\in\M$ and a normal family of characteristic submanifolds $\mathcal N=\set{N_k}_{k\in\N}$ with $N_k\subset\overline{M_k\setminus M_{k-1}}$, then there is a proper isotopy $\Psi^t$ of $\P$ such that $\Psi^1$ is in standard form.
	\elem
	\bpf By Lemma \ref{standardisotopy} we can assume that after a proper isotopy of $\P$ all components of $\mathcal S\eqdef \cup_{k\in\N}\partial M_k\cap\im(\P)$ are properly embedded $\pi_1$-injective surfaces in $\im(\P)$ and the components $S$ of $\mathcal S$ are neither disk nor $\partial$-parallel annuli in $\im(\P)\setminus\P(A\times\set 0)$ for $A$ the collection of annular components of $\Sigma$.
	
	\vspace{0.3cm}
	
	\paragraph{Step 1} Up to a proper isotopy of $\P$ for every component $S$ of $\mathcal S$ we have that $\P^{-1}(S)$ is isotopic to an essential subsurface of $\Sigma\times\set 0$.
	
	\vspace{0.3cm}
	
	By applying \textbf{Claim 1} of Lemma \ref{wald} to a finite type sub-product $\P(\Sigma_S\times[0,\infty))$ containing $\P^{-1}(S)$ we get that every $\P^{-1}(S)$ is isotopic rel $\partial$ to an essential subsurface $F$ of $\Sigma$. 
	
	\vspace{0.3cm}

Since $\Sigma=\cup_{i=1}^\infty\Sigma_i\cup \cup_{i=1}^\infty \Delta_i$ where the $\Delta_i$'s are components of infinite type we have an essential subsurface $T\subset\Sigma$ such that $T\eqdef\cup_{i=1}^\infty\Sigma_i\cup \cup_{i=1}^\infty T_i$ where the $T_i
\neq\emp$ are finite type hyperbolic essential subsurfaces of the $\Delta_i$'s. We denote by $\P_{T}\eqdef \P\vert_{T\times[0, \infty)}\subset \P$ the subproduct it generates. We denote by $\mathcal S_T\subset \mathcal S$ the subcollection of components of $\mathcal S$ that do not intersect $\P(T\times\set 0)$. By properness of the embedding $\P$ we see that for all $i\in\N$: 
$$\P(\Sigma_i\times [0,\infty))\cap \mathcal S\setminus \mathcal S_T\coprod\P(\Delta_i\times [0,\infty))\cap \mathcal S\setminus \mathcal S_T$$
 has finitely many components.

\begin{center}
 \begin{figure}[h!]\centering	\def\svgwidth{350pt}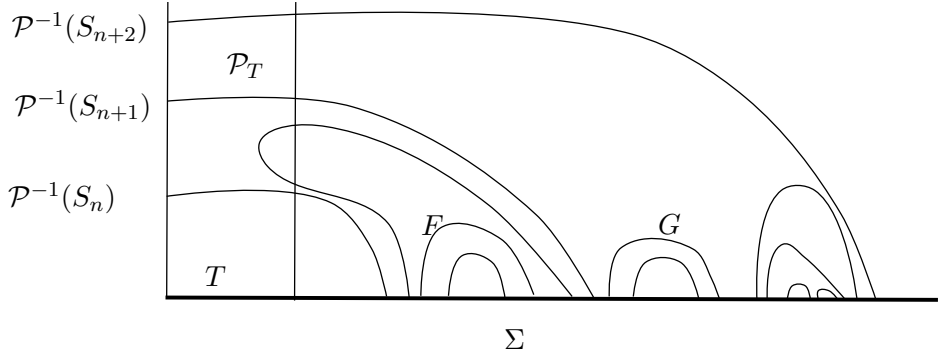
\caption{Schematic of the intersection of $\cup_{k\in\N}\P^{-1}(\partial M_k)$ in $\Sigma\times[0,\infty)$ where $F,G$ are surfaces in $\P^{-1}(\mathcal S\setminus\mathcal S_T)$ and we assume that $\Sigma$ is a connected infinite type surface.}\label{fig6}
 	 \end{figure}
 \end{center}

Note that not all surfaces of $\P(\Delta_i\times [0,\infty))\cap\mathcal S_T$ can be $\partial$-parallel annuli or disks in $\P_T$ since then by Lemma \ref{outermost} and an iterative argument we have a proper isotopy of $\P\vert_{T_i\times[0\infty)}$ such that $\im(\P\vert_{T_i\times[0\infty)})$ does not intersect $\cup_{k\in\N}\partial M_k$. Therefore, in each $\P(\Delta_i\times[0,\infty))$ we have a sequence $\set{S_n^i}_{n\in\N}\subset \mathcal S_T$ such that for all $n\in\N: S_n^i\cap \im(\P_T)$ is an essential properly embedded surface in $\im(\P_T)$. Thus, in $\Sigma\times[0,\infty)$ we have the configuration depicted in Figure \ref{fig6}.

We denote by $F_n^i$ the essential subsurface of $\Sigma$ that $\P^{-1}(S_n^i)$ is isotopic to rel $\partial$. Since Lemma \ref{wald} does isotopies supported in the image of $\P$ we can apply it connected component by connected component and we can assume that for all $n$ $\P^{-1}(S_n^i)\cap T_i\times[0,\infty)=T_i\times\set{a_n^i}$. Thus, we have that for all $n\in\N$ the surfaces $\P^{-1}(S_n^i)$ co-bound with $\Sigma$ $I$-bundles $H_n^i$. Moreover, since for all $n\in\N: T_i\times[0,\infty)\cap H_n^i=T_i\times[0,a_n]$ and $S_n^i\cap S_{n=1}^i=\emp$ we get that $H_n^i\subsetneq H_{n+1}^i$. Moreover, by properness of $\P$ we have that $\cup_{n\in\N} H_n$ is an open and closed submanifold of $\Sigma\times[0,\infty)$ thus, we have that $\cup_{n=1}^\infty H_n=\Sigma\times[0,\infty)$.

Moreover, since each component of $\overline{\im(\P)\setminus\cup_{n,i\in\N} S_{n}^i}$ is essential it is either an $I$-bundle, if $S_{n}^i\iso S_{n+1}^i$ or a mixed $I$-bundle homeomorphic to $S_{n+1}^i\times I$.
		
We now want to remove all components of intersection of $\im(\P)\cap\cup_{k\in\N}\partial M_k$ that are not $S_n^i$ for some $i,n$. By Lemma \ref{wald} we have an isotopy supported in the image of the products of finite type so that they are in standard form. Therefore, we only need to worry about the components of $\P$ that are of infinite type.
			
\vspace{0.3cm} 

	\paragraph{Claim:}Let $\Delta_i$ be a component of infinite type, then up to a proper isotopy of $\P$ supported in $\im(\P\vert_{\Delta_i\times[0,\infty)})$ we have that $\cup_{k\in\N} \partial M_k\cap\im( \P\vert_{\Delta_i\times[0,\infty)})=\cup_{n\in\N}S_n^i$.
	
\vspace{0.3cm}

\bpfc		
	 Via a proper isotopy of $\P$ supported in $\im(\P\vert_{\Delta_i\times[0,\infty)})$ we can assume that 
	 $$\P^{-1}(S_1)\subset \Sigma\times\set 0$$
	 so that $\P(\Delta_i\times [0,\infty))\cap \mathcal S=\P(\Delta_i\times [0,\infty))\cap \mathcal S_T$. We will now do isotopies of the $H_n^i$ relative to $S_n^i,S_{n+1}^i$. All components of $\mathcal S^i\eqdef \P(\Delta_i\times [0,\infty))\cap \mathcal S_T$ that are not $\set{S_n^i}_{n\in\N}$ are contained in a $H_n^i\cong  S_{n+1}\times I$ and are essential. We denote this collection of components $\mathcal S_n^i$. By properness of $\P$ and the fact that $H_n^i$ is compact we get that $\mathcal S_n^i$ has finitely many components $L_1,\dotsc , L_k$ and $\P^{-1}(L_j)$ is isotopic to a subsurface $F_j\subset \Sigma\times\set 0$. Moreover, since the $\P^{-1}(L_j)$ are pairwise disjoint and separating in $H_n^i$ we can find an innermost one. That is if $F_j\subset\Sigma$ is the surface that $\P^{-1}(L_j)$ is properly isotopic to then there are no other components $L_h$ of $\mathcal S^i_n$ such that $\P^{-1}(L_h)$ is contained in the submanifold $J$ bounded by $\P^{-1}(L_j)\cup F_j$. Then, by a proper isotopy of $\P\vert_{H_n^i}$ supported in $\im(\P\vert_{\Delta_i\times[0,\infty)})$ that is the identity on $S_n^i,S_{n+1}^i$ we can push $F_j$ to $\P^{-1}(L_j)$ to reduce $\abs{\pi_0(\mathcal S^i_n)}$. Thus, by concatenating these finitely many isotopies we obtain a proper isotopy $\psi_n^i$ of $\P$ supported in $\im(\P\vert_{H_i^n})$ that is constant on $S_n^i,S_{n+1}^i$ such that $\P^{-1}(\mathcal S^i_n)=S^i_n\coprod S^i_{n+1}$. Since the $\psi_n^i$ are constant on $S_n^i,S_{n+1}^i$ they can be glued together to obtain a proper isotopy of $\P$ supported in $\im(\P\vert_{\Delta_i\times[0,\infty)})$ so that $\cup_{k\in\N} \partial M_k\cap\im( \P\vert_{\Delta_i\times[0,\infty)})=\cup_{n\in\N}S_n^i$.\epfc
	 
So far we have $\Sigma\times[0,\infty)=\cup_{n\in\N} J_n$ and for all $n$ $\P(J_n)$ are mixed $I$-bundles such that there exists $k_n$ with $\P(J_n)\subset X_{k_n}$. Since $\P(J_n)$ is a mixed $I$-bundle we can decompose it in the window $w_n$ and $Q_n$ the non-window part, with an abuse of notation we denote their preimages in $J_n$ by the same name. 

Thus for $\Sigma_i\eqdef \P^{-1}(S_i)$ in $\Sigma\times[0,\infty)$ we have the following configuration:

\begin{center}
 \begin{figure}[h!]\centering	\def\svgwidth{360pt}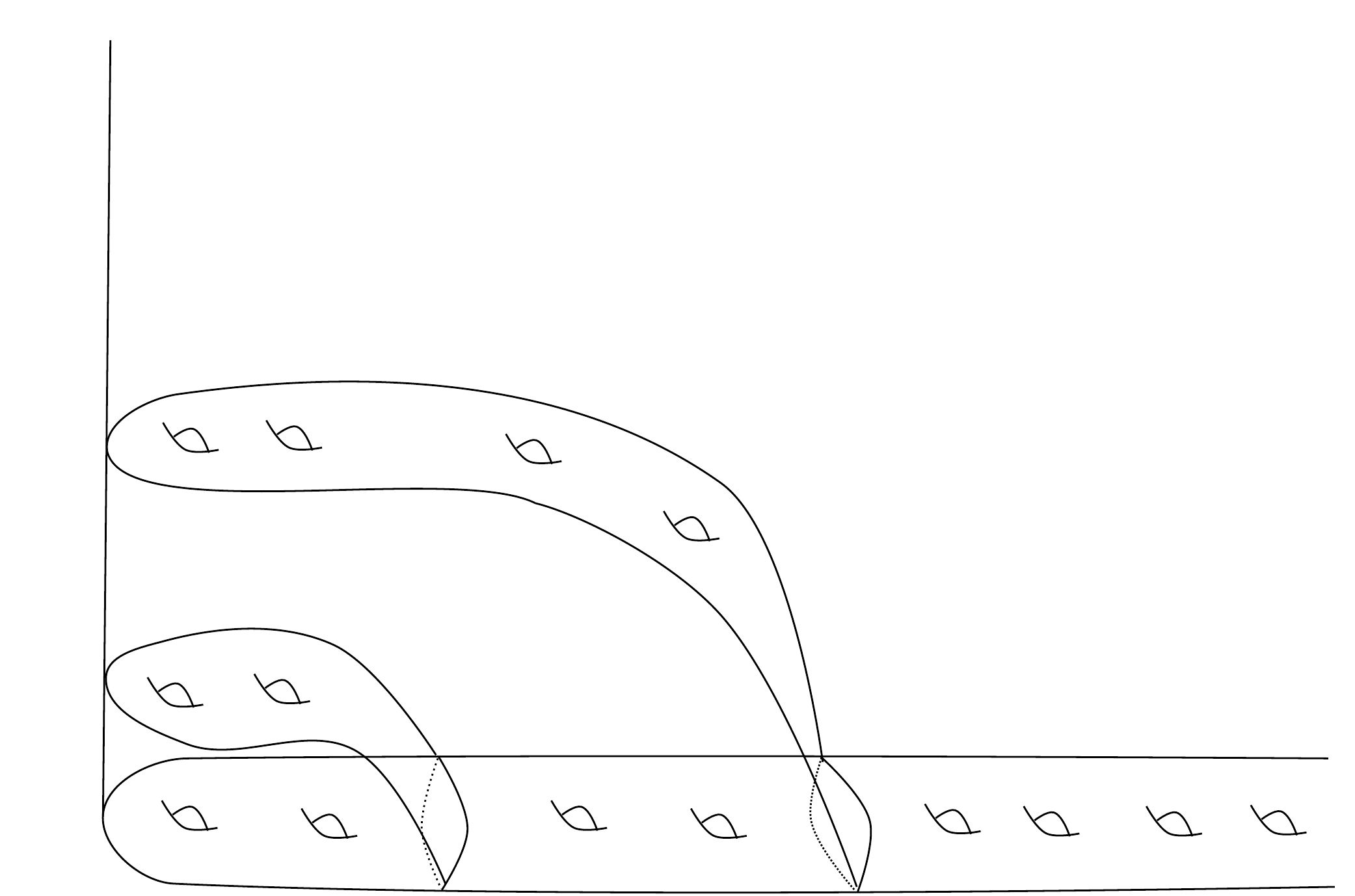

 	 \end{figure}
 \end{center}
 
 \vspace{0.3cm}
 
 \paragraph{Step 2:} Up to a proper isotopy of $\P$ we have that for all $k$ all components of $\im(\P)\cap X_k$ are essential $I$-bundles.
 
 \vspace{0.3cm}
 We will again do isotopies supported in $\im(\P)$ so we can assume that all products are connected. We see that $S_n$ is the boundary of the window $w_n$ and the non-window part $Q_{n-1}$ of the mixed $I$-bundle $I_{n-1}$ is isotopic into $w_n$. Via an isotopy of $\P$ supported in $X_{k_n}\cup X_{k_{n-1}}$ that is the identity on $\Sigma_{n-1}$ and $S_{n+1}$ we can isotope $Q_{n-1}$ into $w_{n}\subset \P(J_{n})$. The $n$-th isotopy is supported in neighbourhoods of the non-window part $Q_n$ and the image is contained in $w_{n+1}$. Therefore, the support of the $n+1$ isotopy does not intersect $w_{n-1}$ thus the composition of these isotopies yields a proper isotopy of $\P$.
 
  Thus one gets a ``staircase picture" in which every step is isotopic to a window in $N_k$, hence it is an essential $I$-bundle, and $\Sigma\times\set 0$ is the boundary of the ``stairs".

\begin{center}
 \begin{figure}[h!]\centering	\def\svgwidth{350pt}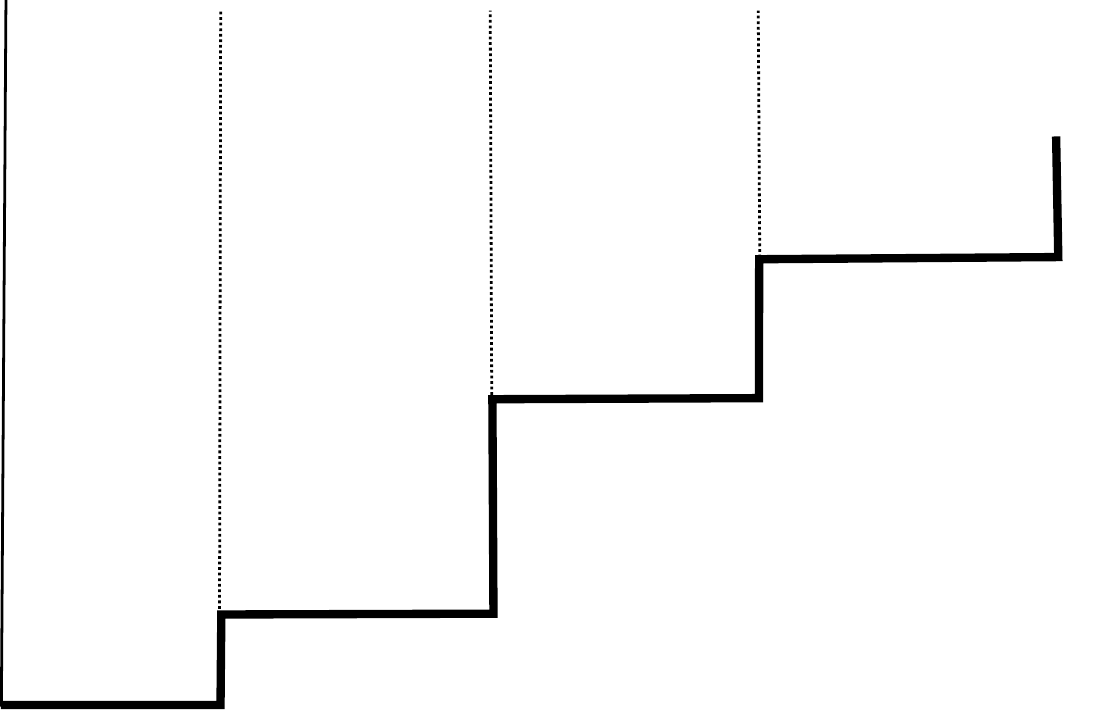
\caption{Where the surfaces $S_i$ coming from the exhaustion induce products $P_i$ whose union is $P$.}\label{staircase}
 	 \end{figure}
 \end{center}
 
 	\paragraph{Step 3:} Up to proper isotopy $\im(\P)\cap X_k\subset N_k$.
	
This follows from the corresponding \textbf{Step} in Theorem \ref{prodstandardform}. \epf

Thus we obtain:

\bthm\label{maxbordification}
Let $M\in\M$ be an open 3-manifold. Then, there exists a unique maximal bordification $[(\overline M,\iota)]\in\cat{Bor}(M)$.
\ethm

\subsubsection{Minimal Exhaustions}\label{2.1.1}

If $M\in\M$ we have two types of $I$-bundles between the boundaries of the gaps in the compact exhaustion: \emph{type I} products have as bases compact surfaces while \emph{type II} products have as bases closed surfaces. 

 \begin{center}
 \begin{figure}[h!]\centering	\def\svgwidth{250pt}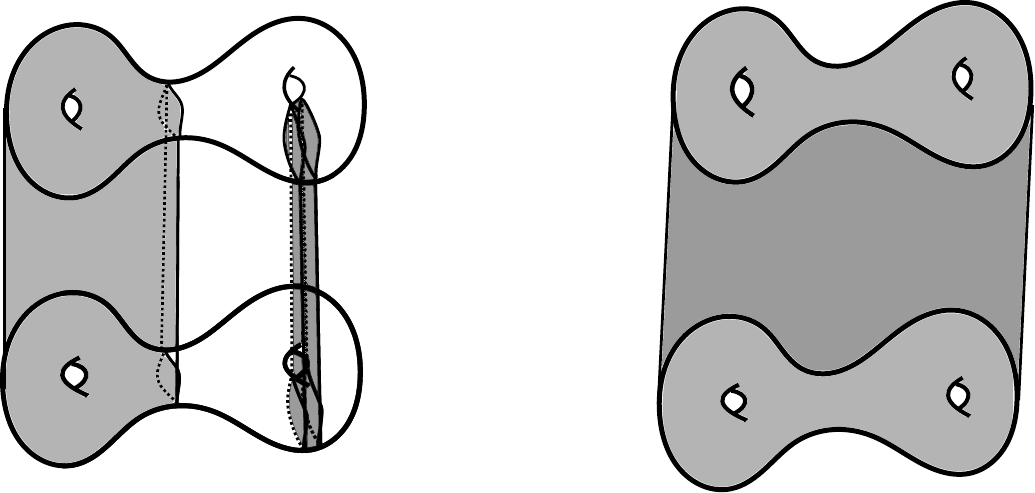

 	\caption{Type I and type II products.}\label{Figure2}
 \end{figure}
 \end{center}

We now want to show that type II products either correspond to tame ends of $M$ or can be thrown away by modifying the exhaustion. 

\bdefi\label{minimalexhdef} We say that a compact exhaustion $\set{M_i}_{i\in\N}$ of $M$ is \emph{minimal} if the two following conditions hold:
\begin{itemize}
	\item[(i)] for all $i<j$: there are no pairs of closed orientable surfaces $F\in \pi_0(\partial M_i)$, $F'\in \pi_0(\partial M_j)$ such that $F\simeq F'$, unless they bound neighbourhoods $U,V$ of the same tame end $E$ of $M_i$;
	\item[(ii)] for all $i$ no component of $\overline{M\setminus M_i}$ is compact.
	\end{itemize}
 \edefi

 These conditions are so that the exhaustion has minimal redundancy. 
 
 \blem\label{minexh}
 Let $M$ be an irreducible 3-manifold with a compact exhaustion $\set{ M_i}_{i\in\N}$ where each $M_i$ has incompressible boundary then, $M$ has a minimal exhaustion.
 \elem
 \bpf For the second condition of a minimal exhaustion we just look at the various $\overline{M\setminus M_i}$ and whenever we see a compact component we add it to all $M_j$ with $j\geq i$. By repeating this process for all components of every $\partial M_k$ we obtain an exhaustion that satisfies the first condition of a minimal exhaustion.\ With an abuse of notation we still denote this new exhaustion by $\set{M_i}_{i\in\N}$.
 
 We now deal with the second condition. Since $M\setminus \text{int}(M_i)$ has no compact component no component of $\partial M_i$ is homotopic to a component of $\partial M_i$ in $M\setminus M_i$ since then by  Lemma \cite[5.1]{Wa1968} we would get that it is homeomorphic to an $I$-bundle. Assume that for $i<j$ we have two distinct closed incompressible surfaces $F_i,F_j$ in $\partial M_i,\partial M_j$, respectively, that are homotopic in $ M_j$. By replacing $F_i$ and $i$, if needed, we can assume that $i$ is minimal. By Lemma \cite[5.1]{Wa1968} the surfaces $F_i,F_j$ bound an $I$-bundle $J$ in $M_j$. Up to an isotopy of $J$ rel $\partial J$ we can assume that $J\cap \partial M_k$, for $i\leq k\leq j$, are level surfaces in $J$. Then, either $J\subset \overline{M_j\setminus M_i}$ or by Lemma \cite[5.1]{Wa1968} $M_i\cong F_i\times I$ and we have $J'\subset  \overline{M_j\setminus M_i}$ given an isotopy from a component of $\partial M_i$ to $F_j$.

 Consider the connected component $U$ of the gap $\overline{M_j\setminus M_i}$ containing the two surfaces. By Lemma \cite[5.1]{Wa1968} we have that: $U\cong F_i\times I$. Then there are two cases: 
 
 \begin{itemize} \item[(i)] either there is $k>j$ such that $F_j$ is not homotopic to any other $F_k\in \pi_0(\partial M_k)$; \item[(ii)] $\forall k>j$ there is $F_k\in\partial M_k$ with $F_k\simeq F_i$. \end{itemize}
 
 In the first case we have a minimal $k\in\N$ with $k>j>i$ such that $F_i$ is not homotopic to any $F_k\in\pi_0(\partial M_k)$. Then, by Lemma \cite[5.1]{Wa1968} the connected component $U$ of $\overline{M_{k-1}\setminus M_i}$ containing $F_i$ and the surface $F_{k-1}$ that it is homotopic to is an $I$-bundle over $F_i$. Hence, we can modify our exhaustion by adding $U$ to all $M_s$ with $i\leq s<k-1$ and leave the other elements of the exhaustion unchanged.
 
In the latter case for all $ k>i$ there is a boundary component $ F_k\in\pi_0(\partial M_k)$ homotopic to $ F_i$. Therefore, $\forall k>i:$ the connected component $U_{k,i}\in\pi_0( \overline{M_k\setminus M_i})$ containing $F_i,F_k$ is an $I$-bundle over $F_i$. Hence, we obtain an exhaustion by submanifolds homeomorphic to $F_i\times I$ of a connected component $E$ of $ \overline{M\setminus  M_i}$. Thus $E\cong F_i\times [0,\infty)$ and $E$ is a tame end of $M$.   \epf

\subsection{Characteristic submanifold for bordifications of manifolds in $\M$}\label{sectionjsj}
In this section we construct the \emph{characteristic submanifold} $(N,R)$ of $\overline M$. Specifically we will prove the following Theorem:

\bthm\label{exisjsj1}
The maximal bordification $(\overline M,\partial\overline  M) $ of $M\in \M$ admits a characteristic submanifold $(N,R)$ and any two characteristic submanifolds are properly isotopic.
\ethm

We will first define characteristic submanifolds for bordifications of manifolds in $\M$ and postponing the proof of existence we prove some general facts about characteristic submanifolds and construct families of characteristic submanifolds for the exhaustion. The proof of Theorem \ref{exisjsj1} is divided into two sections, in which we first prove existence of characteristic submanifolds and then uniqueness.

\subsubsection{Characteristic submanifolds}

 In this subsection we  define characteristic submanifolds for the bordifications of manifolds in $\M$ and describe their components.

\bdefi
Given 3-manifolds $M, N$ and a $\pi_1$-injective submanifold $R\subset \partial N$ a continuous map $f:(N,R)\rar (M,\partial M)$ is \emph{essential} if $f$ is not homotopic via map of pairs to a map $g$ such that $g(N)\subset\partial M$. Similarly we say that a submanifold $N$ is \emph{essential} in $M$ if by taking $R\eqdef N\cap \partial M$ then the embedding is essential.
\edefi

In Definition \ref{jsjdefin} we defined a characteristic submanifold $N$ for a compact irreducible 3-manifold with incompressible boundary $M$. In this setting characteristic submanifolds exists and are unique, up to isotopy, by work of Johannson \cite{Jo1979} and Jaco-Shalen \cite{JS1978}. In the case that $M$ is atoroidal, see \cite[2.10.2]{CM2006}, we get that all components of $N$ fall into the following types:

\begin{itemize}
\item[(1)] $I$-bundles over compact surfaces;
\item[(2)] solid tori $V\cong\mathbb S^1\times\mathbb D^2$ such that $V\cap\partial M$ is a collection of finitely many annuli;
\item[(3)] thickened tori $T\cong\mathbb T^2\times I$ such that $T\cap\partial M$ is a collection of annuli contained in $\mathbb T^2\times\set0 $ and the torus $\mathbb T^2\times\set 1$.
\end{itemize}

We now show that the maximal bordification of manifolds in $\M$ is atoroidal.

\blem
Let $\overline M\in\cat{Bord}(M)$, for $M\in\M$, be the maximal bordification for $M\in\M$ then $\overline M$ is atoroidal.
\elem
\bpf Let $\mathcal T:\mathbb T^2\rar \overline M$ be an essential torus and $\set{M_i}_{i\in\N}$ the exhaustion of $M$. By compactness of $\mathcal T(\mathbb T^2)$ we have that, up to a homotopy pushing $\mathcal T(\mathbb T^2)$ off of $\partial\overline M$, $\mathcal T:\mathbb T^2\rar \overline M$ factors through some $ M_i$. Since $M_i$ is atoroidal and $\mathcal T(\mathbb T^2)\subset  M_i$ is essential we have that $\mathcal T(\mathbb T^2)$ is homotopic into a torus component $T $ of $\partial  M_i$. For all $j>i$: by Waldhausen Cobordism's Theorem \cite{Wa1968} $T$ is isotopic in $M_j$ to a torus component $T_j$ of $\partial M_j$ and so $T,T_j$ cobound an $I$-bundle $I_j$. By the arguments of Theorem \ref{prodstandardform} up to an isotopy of $I_j$ we can assume that $I_j\cap \partial M_k$ are level surfaces of $I_j$ for $i\leq k\leq j$. Thus, either $I_j\cap M_i=T$ or $M_i\cong \mathbb T^2\times I$ and then for all $j>i$ we have that $\overline {M_j\setminus M_i}\cong \mathbb T^2\times I\coprod \mathbb T^2\times I$. In either case, we get that the component of $\overline {M_j\setminus M_i}$ containing $T,T_j$ is homeomorphic to $\mathbb T^2\times I$ and since $j$ was arbitary $\mathcal T(\mathbb T^2)$ is homotopic into $\partial\overline M$. Therefore, every $\pi_1$-injective torus in $\overline M$ is homotopic into $\partial \overline M$ and so inessential. \epf

In our setting we have $\overline M\in\cat{Bord}(M)$, for $M\in\M$, with $\text{int}(\overline M)$ exhausted by compact hyperbolizable 3-manifolds $M_i$ with incompressible boundary. Therefore, we have a collection $(N_i,R_i)\hookrightarrow (M_i,\partial M_i)$ of characteristic submanifolds whose components are of the form (1)-(3). Thus, since $\overline M$ is atoroidal for a characteristic submanifold $N$ of $\overline M$ we expect the components of $N$ to be of the following types:
\begin{itemize}
\item[(i)] $I$-bundles over compact incompressible surfaces;
\item[(ii)] solid tori $V$ with finitely many wings\footnote{Recall that a \emph{wing} is a thickening of an essential annulus $A$ with one boundary component on the solid torus $V$ and one on $\partial\overline M$.  };
\item[(iii)] thickened essential tori $\mathbb T^2\times I$ corresponding to a torus component of $\partial\overline M$ possibly with finitely many wings;
\item[(iv)] limit of nested solid tori or of nested thickened essential tori.
\end{itemize}
Except for (iv) these are the same components that one finds in the usual JSJ decomposition of compact atoroidal 3-manifolds with incompressible boundary. 

The difference for manifolds in $\M$ is that we can have a countable family of nested solid tori or thickened essential tori having no parallel wings. One can think of these as \textit{infinitely winged solid tori} (IWSD) or \textit{infinitely winged essential tori} (IWET). These are solid tori $V$, or thickened tori $T$, with infinitely many wings, specifically in each $ M_i$ we have that $V\cap  M_i$, or $T\cap M_i$, has a component that is isotopic to an essential solid torus $V_i\subset N_i$ or an essential thickened torus $T_i\subset N_i$ with $a_i$-wings and $a_i\nearrow\infty$. 

We will now build such an example.

\bese[A 3-manifold with an infinitely winged solid torus.]\label{IWSD}
												 					
	 										 									Let $(N,\partial N)$ be an acylindrical and atoroidal compact 3-manifold with boundary an incompressible genus two surface (for example see \cite[3.3.12]{Th1997} or Appendix \ref{appendix A}).

																				 Let $T$ be a solid torus with three wings winding once around the soul of the solid torus. The boundary of $T$ is decomposed into $6$ annuli, one for each wing and one between each pair of wings. 
																				 
																				 Consider the manifold obtained by gluing the annular end of the wings of $T$ to three copies of $\Sigma_2\times I$\footnote{By $\Sigma_2$ we mean a genus two surface.} along a neighbourhood of a curve $\gamma\subset\Sigma_2\times\set 0$ separating $\Sigma_2\times \set 0$ into two punctured tori. The resulting 3-manifold has for boundary six copies of $\Sigma_2$. Three boundary components are coming from the three copies of $\Sigma_2\times \set 1$ and the other three are coming from gluing two punctured tori in the $\Sigma_2\times\set 0$'s along an annulus in the boundary of the solid torus. 
																				 
																				 By gluing 3 copies of $N$ along the second type of $\Sigma_2$ we obtain a 3-manifold $X$ as in the picture:
	 									
	 									\begin{center}\begin{figure}[h!]
	 											\centering
	 											\def\svgwidth{150pt}
	 											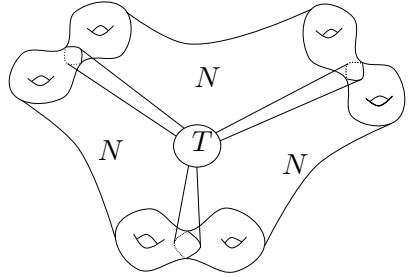
												\caption{An $X$-piece.}
	 											
	 										\end{figure}\end{center}
	 										
	 										The 3-manifold $X$ is hyperbolizable with incompressible boundary and has the property that its characteristic submanifold is given by the solid torus with three wings $T$. We now construct a 3-manifold $M$ by gluing together countably many copies $\set{X_i}_{i=1}^\infty$ of $X$ and product manifolds $P\eqdef \Sigma_2\times [0,\infty)$. We denote by $T_i$  the three winged solid torus in $X_i$. The gluing is given by the following tree pattern in which the gluing maps are just the identity:
	 										
	 										\begin{center}\begin{figure}[h!]
	 												\centering
	 												\def\svgwidth{200pt}
	 												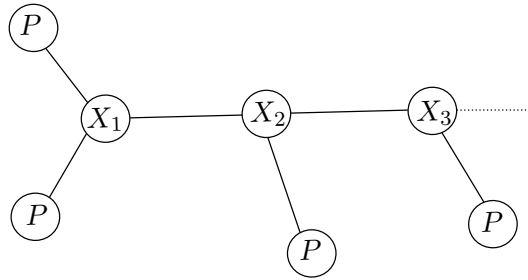
													\caption{The gluing pattern for $M$.}	 											\end{figure}\end{center}

	 											The manifold $M$ has a compact exhaustion given by taking $M_i$ to be the manifold up to the $i$-th copy of $X$ and the compact submanifolds of the product ends given by $\Sigma_2\times [0,i]$. Hence, $\partial M_i$ is formed by $2+i$ copies of $\Sigma_2$ all of which are incompressible. At each $M_i$ the characteristic submanifold is a solid torus $\tau_i$ with $2+i$ wings. Moreover, since the $M_i$ are atoroidal Haken 3-manifold by the Hyperbolization theorem \cite{Kap2001} they are hyperbolizable. Therefore, since all boundaries of the $M_i$ are incompressible and of genus two $M$ is a manifold in $\M$.
												
										 Moreover, the JSJ submanifold of $ M_1$ is given by the solid torus with 3-wings $T_1$ and the JSJ submanifold of the component of $\overline{ M_{j}\setminus M_{j-1}}$ that is not an $I$-bundle is also given by the solid torus $T_j$. Let $T_\infty$ be the submanifold of $M$ obtained by taking all the $\set{T_j}_{j\geq 1}$ and adding to it cylinders going to infinity in all the tame ends. Then $T_\infty$ is an example of an \emph{infinitely winged solid torus} since it is an open 3-manifold that compactifies to a solid torus $V$ and is homeomorphic to $V\setminus L$ where $L\subset\partial V$ is a collection of pairwise disjoint isotopic simple closed curves forming a closed subset of $\partial V$. Namely $L=\set{L_{\frac 1 n}}_{n\in\N}\cup L_0$ is in bijection with the ends of $M$ where $L_0$ is the non-tame end and the $\set{L_{\frac 1 n}}_{n\in\N}$ correspond to the tame ends. Moreover, it is topologised so that $L_{\frac 1n}\rar L_0$ as $n\rar\infty$.				
			
							\begin{center}\begin{figure}[h!]
	 										\centering
	 										\def\svgwidth{350pt}
	 										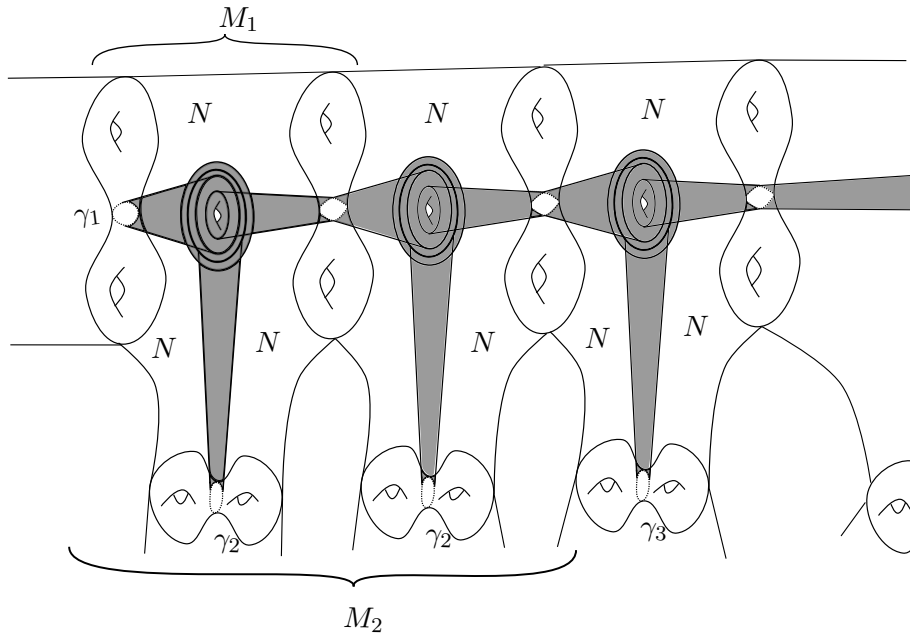
	 										
	 										\caption{The tori $T_j$'s are in grey and the tame ends are unmarked.}\label{fig:charmanifold}
	 									\end{figure}\end{center}	
The maximal bordification $\overline M$ of $M$ has for boundary components an open annulus $A$ and countably many genus two surfaces $\set{\Sigma_i}_{i\in\N}$. The annulus comes from the compactification of a product $P:\mathbb A\times[0,\infty)\hookrightarrow  T_\infty$ going out the non-tame end that is contained in the interior of $T_\infty$. The genus two surfaces $\Sigma_i$ come from compactifying all the tame ends $P_i$.
In the maximal bordification $\overline M$ we have that the characteristic submanifold $N$ is given by $N\eqdef T_\infty\cup_{i\in\N} A_i\cup A'$ where the $A_i\subset \Sigma_i$ are the annuli that $T_\infty$ limits to and $A'\subset A$ is a core annulus for $A$.

It is easy to modify the above example to obtain a 3-manifold containing an IWET by adding to any tame end $P_i$, along the boundary $S_i$, a compact hyperbolizable 3-manifold $Y$ with incompressible boundary $\partial Y\cong T\cup\Sigma_2$  and such that a simple closed loop $\beta$ in the boundary torus $ T$ is isotopic, in $Y$, to the separating curve of the genus two boundary $\Sigma_2$ glued to the separating loop of $P_i$.

															\eese 
															
															\vspace{0.5cm}

Thus, we define:
\bdefi
Given a 3-manifold $M\in\M$ let $ (\overline M,\partial \overline M)$ be the maximal bordification, which could be $M$ itself, then we define the \emph{characteristic submanifold } $(N,R)\hookrightarrow (\overline M,\partial \overline M)$ to be a codimension-zero submanifold satisfying the following properties:
\begin{enumerate}
\item[(i)] every $\Sigma\in\pi_0(N)$ is homeomorphic to either: 
\begin{itemize}
\item an essential $I$-bundle over a compact surface;
\item  an essential solid torus $V\cong \mathbb S^1\times\mathbb D^2$ with $V\cap\partial\overline M$ a collection of finitely many parallel annuli or a non-compact submanifold $V'$ that compactifies to a solid torus such that $V'\cap \partial \overline M$ are infinitely many annuli;
\item an essential thickened torus $T\cong\mathbb T^2\times [0,1]$ such that $T\cap\partial \overline M$ is an essential torus and a, possibly empty, collection of parallel annuli in $\partial T$ or a non-compact manifold $T'$ that compactifies to a thickened torus such that $T'\cap\partial\overline M$ is an essential torus and infinitely many annuli;
\end{itemize}
\item[(ii)] $\partial N\cap \partial\overline M=R$;
\item[(iii)] all essential maps of an annulus $(\mathbb S^1\times I,\mathbb S^1\times\partial I)$ or a torus $\mathbb T^2$ into $ (\overline M,\partial\overline  M)$ are homotopic as maps of pairs into $(N,R)$;
\item[(iv)] $N$ is minimal i.e. there are no two components of $N$ such that one is homotopic into the other.

\end{enumerate}
\edefi

\blem\label{charlimitsolidtori}
Let $M=\cup_{i\in\N} M_i\in\M$ and $\set{T_i}_{i\in\N}$ be a collection of essential solid tori $T_i\subset M_i$ such that for $j>i$ $T_j\cap M_i$ is compact, contains $T_i$ and $\overline{T_{i+1}\setminus T_i}$ are essential solid tori. Moreover, assume that the inclusion maps $\iota_i:T_i\hookrightarrow T_{i+1}$ induce isomorphisms on $\pi_1$. Then the direct limit $T\eqdef \varinjlim_i T_i$ is a properly embedded submanifold of $M$ such that $T\cong \mathbb S^1\times\mathbb D^2\setminus L$ for $L$ a closed subset of $\partial (\mathbb S^1\times\mathbb D^2)$ consisting of parallel simple closed curves.
\elem
\bpf Since all inclusions induce isomorphism on the fundamental groups and $\pi_1(T)=\varinjlim_i \pi_1(T_i)$ we have that $\pi_1(T)\cong\Z$. A non-compact manifold $N$, possibly with boundary, is a \emph{missing boundary manifold} if $N\cong \overline N\setminus L$ where $\overline N$ is a manifold compactification of $N$ and $L$ is a closed subset of $\partial\overline N$. By Tucker's Theorem \cite{Tu1974} the manifold $T$ is a missing boundary manifold if the complement of every compact submanifold has finitely generated fundamental group. Since $T=\cup_{i=1}^\infty T_i$ it suffices to check the above condition for the $T_i$. 

Let $Q\in\pi_0(\overline{T\setminus T_i})$, then since $\overline{T_j\setminus T_i}$ are solid tori $Q$ is either a solid torus or another direct limit of nested solid tori in which the inclusions induce isomorphism in $\pi_1$. In either case $\pi_1(Q)\cong \Z$. Therefore $T$ compactifies to $\widehat T$ and $T$ is homeomorphic to $\widehat T\setminus L$ where $L\subset \partial \hat T$ is a closed set. Since $\hat T$ is compact, irreducible and $\pi_1(\hat T)\cong\Z$ we have that $\hat T\overset\psi\cong V\eqdef \mathbb S^1\times\mathbb D^2$, see \cite[Theorem 5.2]{He1976}. 

\vspace{0.3cm}

\paragraph{Claim:} Up to a homeomorphism of $V$ the set $L$ is a union of of parallel curves.

\vspace{0.3cm}

\bpfc Given the homeomorphism $\psi: T\rightarrow V\setminus L$ for $V=\mathbb S^1\times\mathbb D^2$ we see that every $T_i\subset T$ is mapped to a solid torus $V_i\subset V\setminus L$ such that $\partial V_i=S_i\cup A_1^i\cup\dotsc\cup A_{n_i}^i$ where $S_i\subset \partial V\setminus L$ and the $A_j^i$'s are compact properly embedded annuli in $V$. Moreover, since the $A^i_j$'s are $\pi_1$-injective embedded annuli in $\partial V_i$ they are isotopic annuli in $\partial V_i$ and are the images of the annuli of $T_i$ contained in $\partial M_i$ to which new solid tori get glued in $\overline{M_{i+1}\setminus M_i}$ to obtain $T_{i+1}$. Moreover, since $V\setminus V_i$ is a collection of solid tori we have that each $A_{j_k}^i$ is $\partial$-parallel in $V$. We define $\mathcal A_i\eqdef\hat A_1^i\cup\dots \cup\hat A_{m_i}^i\subset \partial V$ for $\hat A_{j_k}^i$ the annulus in $\overline{\partial V\setminus V_i}$ co-bounded by $\partial A_{j_k}^i$. Every component $\gamma$ of $L$ is given as a countable intersection of a sequence of the annuli $\hat A_j^i$. Moreover, we fix a fiber structure on $V$ such that it is fibered by circles. Then, for each $V_i$ we want to construct a homeomorphism $\phi_i$ of $V$ such that:

\begin{itemize}
\item $\phi_i\vert _{V_{i-1}}=\phi_{i-1}$;
\item $\phi_i(\mathcal A_i)$ are fibered annuli contained in $\phi_{i-1}(\A_{i-1})$.
\end{itemize} 

Assume we defined such a $\phi_j$ for all $j\leq i$. To define $\phi_{i+1}$ we only need to change $\phi_i$ in the solid tori co-bounded by $\phi_i(\mathcal A_i)$ and $\phi_i(\cup_{j=1}^{n_i} A_j^i)$ which lie in the complement of $\phi_i(V_i)$. Each such solid torus $W_k$ has boundary given by $\phi_i(\hat A_k^i\cup A_{j_k}^i)$ and in $\phi_i(\hat A_{j_k}^i)$ contains some of the annuli $\set{\phi_i(\hat A_k^{i+1})}_{1\leq k\leq m_{i+1}}$. Then by an isotopy $\psi_{i,k}^t$ supported in $\phi_i(\hat A_k^i)$ that is the identity on $\partial \phi_i(\hat A_k^i)$ we can make the $\phi_i(\hat A_k^{i+1})$ fibered in $\phi_i(\hat A_j^{i})$. By extending the isotopy $\psi_{i,k}^t$ to the solid torus and taking the time one map we obtain the required homeomorphism $\phi_{i+1}\eqdef \psi_{i,k}^1$. 

	Finally, the map $\phi\eqdef\varinjlim \phi_i$ is a homeomorphism from $V\setminus L$ to $V\setminus L$ such that now $L=\cap_i \phi_i(\cup_j\hat A_j^i)$ where the $\phi_i(\hat A_j^i)$ are now compatibly fibered annuli. Thus every component of $L$ is also fibered and if we assume that the $\phi_i$ are strictly contracting on the annuli we get that the components of $L$ do not contain any annuli and are indeed parallel loops.
\epfc

By taking the homeomorphism $\phi\circ\psi :T \rightarrow V\setminus L$ we obtain the required conclusion. 
\epf

Similarly we obtain:

\bcor\label{charesstori}
Let $M=\cup_{i\in\N} M_i\in\M$ and $\set{T_i}_{i\in\N}$ be a collection of essential tori such that for $j>i$ $T_j\cap M_i$ is compact, contains $T_i$ and $\overline{T_{i+1}\setminus T_i}$ are essential solid tori or thickened essential tori. Moreover, assume that the inclusion maps $\iota_i:T_i\hookrightarrow T_{i+1}$ induce isomorphisms on $\pi_1$. Then the direct limit $T\eqdef \varinjlim_i T_i$ is a properly embedded submanifold of $M$ such that $T\cong\mathbb T^2\times (0,1]\setminus L$ for $L$ a closed subset of $\mathbb T^2\times \set 1$ consisting of parallel simple closed curves.

\ecor

The two above Lemma deal with two of the types of components we expect to have. For $I$-bundles we have: 

\blem
Let $\overline M\in\cat{Bord}(M)$ be the maximal bordification of $M\in\M$ and let $\iota:(F\times I,F\times\partial I)\hookrightarrow (\overline M,\partial\overline M)$ be an essential $I$-bundle over a connected surface $F$. Then, the surface $F$ is compact.
\elem
	\bpf Let $\P$ denote the proper embedding $\iota\vert_{F\times (0,1)}$ and let $\set{M_i}_{i\in\N}$ be the exhaustion of $M\in\M$. Since $\P$ is a proper $\pi_1$-injective embedding by Lemma \ref{esssubsurface} we have that:
	
	\vspace{0.3cm}
	
\paragraph{Step 1:} Up to a proper isotopy of $\P$ we can make all components of $\cup_{i\in\N}\partial M_i\cap\P$ be $\pi_1$-injective subsurfaces of $\P$ and no component is a disk.

\vspace{0.3cm}

We now claim that:

	\vspace{0.3cm}
	
\paragraph{Step 2:} Up to a proper isotopy of $\P$, supported in $\im(\P)$,  no component $S$ of $\mathcal S\eqdef \cup_{i\in\N}\partial M_i\cap\im(\P)$ is a $\partial$-parallel annulus.

\vspace{0.3cm}

Let $\mathcal A_i$ be the collection of annuli of $\mathcal S_i\eqdef\partial M_i\cap\im(\P)$ that are $\partial$-parallel in $\im(\P)$. Since $\P$ is a proper embedding we have that for all $i\in\N$ $\abs{\pi_0(\mathcal A_i)}<\infty$. By an iterative argument it suffices to show the following:

	\vspace{0.3cm}
	
\paragraph{Claim } If for $1\leq n< i$ $\mathcal A_n=\emp$ then via an isotopy $\phi_i^t$ of $\P$ supported in $\overline{M\setminus M_{i-1}}\cap\im(\P)$ we can make $\mathcal A_{i}=\emp$.

\vspace{0.3cm}

\bpfc For all $i\in\N$ we have $0<a_i<b_i<\infty$ such that $\mathcal A_i\subset \P(F_i\times[a_i,b_i])$ for $F_i\subset F$ a compact essential subsurface of $F$. 

Denote by $A_1,\dotsc, A_n$ the $\partial$-parallel annuli in $\mathcal A_i$. By applying Corollary \ref{pushanninter} to $\P(F_i\times[a_i,b_i])$ we have a local isotopy $\phi_i^t$ of $\P$ that removes all these intersections. The isotopy $\phi_i^t$ is supported in a collection of solid tori $\mathcal V_i\subset F_i\times[a_i,b_i]$ such that $\partial\mathcal V_i\cap \partial F_i\times[a_i,b_i]\subset\partial F\times[a_i,b_i]$ and $\phi_i^t$ is the identity outside a neighbourhood of $\partial \mathcal V_i\setminus\partial F\times\R$ thus it can be extended to the whole of $\P$. Moreover, if we consider for $n<i$ a component of intersection of $\partial M_n\cap \P(\mathcal V_i)$ then it is either a boundary parallel annulus or a disk. However, we assumed that for $n<i$ $\mathcal A_n=\emp$ and by \textbf{Step 1} no component of $\cup_{k\in\N}\partial M_k\cap\im(\P)$ is a disk thus, the solid tori $\mathcal V_i$ that we push along are contained in $\im(\P)\cap \overline{M\setminus M_{i-1}}$. Therefore, we get a collection of solid tori $\mathcal V_i\subset\im( \P)\cap\overline{M\setminus M_{i-1}}$ such that pushing through them gives us an isotopy $\phi_i^t$ of $\P$ that makes $\mathcal A_i=\emp$. \epfc

Since for all $i\in\N$ $\text{supp}(\phi_i^t)=N_\epsilon(\mathcal V_i)$ is contained in $\overline{M\setminus M_{i-1}}$ the limit $\phi^t$ of the $\phi^t_i$ gives us a proper isotopy of $\P$ such that for all $i\in\N$ $\mathcal A_i=\emp$.

\vspace{0.3cm}

By \textbf{Step 2} every component $S$ of $\mathcal S$ the surface $\iota^{-1}(S)$ is an essential surface in $F\times \R$. In particular, since $\P$ is properly embedded we have a component $\Sigma$ of some $\partial M_i$ such that $S\eqdef \Sigma\cap\im(\P)\neq\emp$. If $S=\Sigma$ we get a contradiction since then we have an incompressible closed surface $\iota^{-1}(S)$ in the $I$-bundle $F\times I$. Therefore, we must have that $S$ is an essential proper subsurface of $\Sigma$ with $\partial S\subset \P(\partial F\times \R)$.

	However, since the surface $S$ is compact and properly embedded in $\P$ we can find $0<t_1<t_2<1$ such that $\iota^{-1}(S)\subset F'\times (t_1,t_2)$ for $F'\subset F$ a compact surface such that not all boundary components of $F'$ are boundary components of $F$. Since $\iota^{-1}(S)$ is an essential properly embedded subsurface of $F'\times [t_1,t_2]$ by \cite[3.1,3.2]{Wa1968} we have that $\iota^{-1}(S)$ is isotopic to an essential sub-surface $S'$ of $F'\times\set{t_1}$. However, $\partial \iota^{-1}(S)\subset\partial F\times I$ which are a subset of $\partial F'\times I$ and so we get a contradiction since then $\partial S'\subsetneq \partial F'$ are zero in $H_1(F')$.\epf

\subsubsection{Existence of characteristic submanifolds} In this section we prove that if $M\in\M$ then the maximal bordification $\overline M$ admits a characteristic submanifold. Before proving this existence statement we need to show that given $M\in\M$ and a normal family of characteristic submanifolds for the gaps we can find a family of characteristic submanifolds $N_i$ of the $M_i$ that are compatible with each others:

\bprop

Given $M=\cup_{i\in\N}M_i\in\M$ and a normal family of characteristic submanifolds $\set{N_{i-1,i}}_{i\in\N}$ for $X_i\eqdef \overline{M_i\setminus M_{i-1}}$, then we have characteristic submanifolds $C_i\subset M_i$ such that for all $i\in\N$ we have that for $j\geq i$: $C_{j}\cap M_i\subset C_i$ and $C_n\subset C_{n-1}\cup N_{n-1,n}$.
\eprop

In the next series of Lemmas we will construct a family $N_i$'s of characteristic submanifolds for the $M_i$'s such that $N_i\cap M_j\subset N_j$ whenever $i>j$.

\blem \label{intersection}Let $M_1\subset\text{int}(M_2)$ be hyperbolizable 3-manifolds with incompressible boundary and let $N_1,N_2$ be their characteristic submanifolds. Given distinct components $P,Q\in\pi_0(N_2)$ if every component $P\cap M_1$ is an essential submanifold of $M_1$ and $P\cap M_1$ has a component isotopic into $Q\cap M_1$ then one of $P$ or $Q$ is an $I$-bundle over a surface $F$ with $\chi(F)<0$ and the other is either a solid torus or a thickened essential torus.
\elem
\bpf

Let $S\in\pi_0(P\cap M_i)$ be a component isotopic into $S'\in\pi_0(Q\cap M_1)$. If $S\cong F\times I$, with $\chi (F)<0$, then $S$ and $S'$ are isotopic into an $I$-bundle component of $N_1$, thus $P$ (or $Q$) is a sub-bundle of $Q$ (or $P$) and we reach a contradiction since then they are not distinct components of $N_2$. 

If $S\cong \mathbb S^1\times\mathbb D^2$ is a solid torus we have that either $S'$ is an $I$-bundle $F'\times I$, with $\chi(F')<0$, and $S$ is homotopic into $\partial F'\times I$ or $S'$ is either a solid torus or an essential thickened torus. In the first case we have that $P$ is either a solid torus component or a thickened essential torus component of $N_2$ while $Q$ is an $I$-bundle over a surface of negative Euler characteristic and we are done. 

In the second case we have that $S'$ is homeomorphic to a solid torus. Since $S'$ is homotopic into $S$ we can find an embedded annulus $A$ in $M_1\setminus S\cup S'$ connecting $\partial S$ to $\partial S'$ and denote by $A'$ a regular neighbourhood $A$ intersecting $P,Q$ only in neighbourhoods of $A\cap\partial P\cup\partial Q$. Since both $P,Q$ are solid tori we get that $P\cup A'\cup Q$ is homeomorphic to a solid torus $V$. Thus we get an essential map: $f:V\rar M_2$ whose image is $P\cup A'\cup Q$. By properties of characteristic submanifolds we have that $V$ is homotopic into a component $T$ of $N_2\setminus P\cup Q$. However, this contradicts the minimality properties of $N_2$ since then $N_2\setminus P\cup Q$ would also be characteristic.

Finally if $S\cong\mathbb T^2\times I$ is isotopic into $S'$ we have that $S'$ is also homeomorphic to $\mathbb T^2\times I$ and since they contain the same $\Z^2$ subgroup of $\pi_1(M_2)$ they are the same component of $N_2$. \epf 

We now prove the iterative step of constructing a compatible family of characteristic submanifold.

\blem\label{char-submanifold}
Let $ M_1\subset\text{int}(M_2)$ be hyperbolizable 3-manifolds with incompressible boundary and let $(N_1,R_1)$, $(N_2,R_2)$  and $(N_{12}, R_{12})$ be characteristic submanifold of $M_1, M_2$ and $\overline {M_2\setminus M_1}$ respectively. Moreover, assume that $N_1,N_{12}$ form a normal family, then we can isotope $N_2$ in $ M_2$ such that $N_2\subset N_1\coprod N_{12}$. \elem 
\bpf
If, up to isotopy, $N_2\cap M_1=\emp $ then we can isotope $N_2$ so that $N_2\subset N_{12}$ and there is nothing else to do. So can we assume that the intersection, up to isotopy, is not empty thus, some component of $N_2\cap M_1$ is essential in $M_1$. 

\vspace{0.3cm}

\paragraph{Step 1:} Up to an isotopy of $N_2$ we have that every component of $N_2$ intersects $M_1$ and $\overline{M_2\setminus M_1}$ in essential $I$-bundles, essential solid tori or thickened essential tori.

\vspace{0.3cm}

By an isotopy of $N_2$ and a general position argument we can minimise $\abs{\pi_0(\partial M_1\cap N_2)}$ and have that $\partial M_1\cap N_2$ are $\pi_1$-injective surfaces, see Lemma \ref{esssubsurface}. 

Let $P\in\pi_0(N_2\cap  M_1)$ be a component of intersection coming from an $I$-bundle component $P'\cong F\times I$, with $\chi(F)<0$, of $N_2$. Since the components $S$ of $P'\cap \partial M_1$ are essential and with boundary in the side boundary of the $I$-bundle $P'$ by \cite[3.1,3.2]{Wa1968} they are isotopic to subsurfaces of the lids of the $I$-bundle region. Therefore, we have that $P\cong F\times I$ is an $I$-bundle and since it is $\pi_1$-injective it is essential.

If $P'\cong \mathbb S^1\times\mathbb D^2$ is a solid torus component of $N_2$ then $A\eqdef P'\cap \partial M_1$ is a collection of $\partial $-parallel annuli in $P'$. The annuli $A$ decompose $P'$ into a collection of solid tori each of which is contained in either $M_1$ or $\overline {M_2\setminus M_1}$. If a solid torus component $T$ of $P'\cap M_1$ is inessential, i.e. it either is $\partial$-parallel or it has, at least, two wings $w_1,w_2$ in $M_1$ that are parallel, then by an isotopy of $N_2$ that either pushes $P'$ outside of $M_1$ or pushes $w_2$ along $w_1$ outside of $M_1$ we can decrease $\abs{\pi_0(\partial M_1\cap N_2)}$ contradicting the assumption that it was minimal.

Similarly if $P\cong\mathbb T^2\times I$ we have that $\partial M_1$ decomposes $P$ into one essential thickened essential torus in $M_1$ and essential solid tori contained in $M_1$ and $\overline {M_2\setminus M_1}$.

\vspace{0.3cm}

Moreover, by properties of a normal family we can assume that up to another isotopy supported in a neighbourhood $U_1$ of $\partial M_1$ we have that $N_2\cap\partial M_1\subset R_1\cap N_{12}$.

\vspace{0.3cm}

\vspace{0.3cm}

\paragraph{Step 2:} Up to an isotopy of $N_2$ we have that $N_2\cap M_1\subset N_1$.

\vspace{0.3cm}

Since $N_2\cap M_1$ is a collection of essential Seifert-fibered 3-manifolds and $I$-bundles such that $N_2\cap \partial M_1\subset R_1$ by JSJ theory we can isotope them rel $R_1$ into $N_1$.

 Let $P\cong F\times I$ be an $I$-bundle component of $N_2$ with $\chi(F)<0$. By Step 1 we can assume that $P\pitchfork \partial M_1$ and $K_P\eqdef \abs{\pi_0(P\cap\partial M_1)}$ is minimal. Moreover, we can assume that every component of $P\cap\partial M_1$ is in $R_1=\partial N_1\cap \partial M_1$. Then $P\cap M_1=P_1\coprod P_2\coprod\dotsc\coprod P_n$ are essential $I$-bundles in $M_1$. Thus, we can isotope the $\coprod_{i=1}^nP_i$ rel $\partial P_i\cap R_1$ into $N_1$. We repeat this for all $I$-bundles of $N_2$ and by Lemma \ref{intersection} we do not need to worry of them being parallel in $M_1$. We denote by $N_2'$ the resulting submanifold. The submanifold $N_2'$ is isotopic to $N_2$ hence characteristic for $M_2$.

Let $P\in \pi_0(N_2)$ be a solid torus component. By Step 1 and the fact that $N_1,N_{12}$ form a normal family we have that each component of $P\cap\partial M_1$ is in $R_1$. Then, $P$ is decomposed by $\partial M_1$ into solid tori and annuli that are contained in $M_1$ and $\overline{M_2\setminus M_1}$. Moreover, each such component is essential,  thus every component of $P\cap M_1$ is either an essential solid torus with $k\geq 3$ wings in $N_1$ or a thickened cylinder. Each solid torus component is then isotopic into a solid torus component of $N_1$ and each annular component is isotopic into a solid torus or an $I$-bundle. Say that an annular component $A$ of $P\cap M_1$ is isotopic into the side boundary of an $I$-bundle component $Q\overset \psi\cong F\times I$ of $N_2\cap M_1$. Then, up to a further isotopy of $Q$ we can assume that both $A$ and $Q$ are contained in $N_1$.

The same process applies when $P\in\pi_0(N_2)$ is a thickened essential torus, the only difference is that if the boundary torus $T$ is in $M_1$ we also have an essential thickened torus component in $N_2\cap M_1$.

\vspace{0.3cm}

\paragraph{Step 3:} Up to an isotopy supported in $\overline{M_2\setminus M_1}$ of $N_2$ we have that $N_2\cap\overline{M_2\setminus M_1}\subset N_{12}$.

\vspace{0.3cm}

By Step 1 every component of $N_2\cap \overline{M_2\setminus M_1}$ is an essential $I$-bundle, a solid torus or a thickened essential torus. Then, by JSJ theory we can isotope $N_2\cap \overline{M_2\setminus M_1}$ into $N_{12}$. Moreover, since the components of $N_2\cap \partial M_1$ isotopic into $N_{12}\cap\partial M_1$ are already contained in $N_{12}\cap \partial M_1$, by properties of a normal family, we can assume that the isotopy is the identity on $\partial M_1$.

\vspace{0.3cm}

The composition of the isotopies yields the required characteristic submanifold $N_2\subset N_1\cup N_{12}$. \epf

Thus, we have:

\bprop\label{infjsj}
Given $M=\cup_{i\in\N}M_i\in\M$ and a normal family of characteristic submanifolds $\set{N_{i-1,i}}_{i\in\N}$ for $X_i\eqdef \overline{M_i\setminus M_{i-1}}$, then we have characteristic submanifolds $C_i\subset M_i$ such that for all $i\in\N$ we have that for $j\geq i$: $C_{j}\cap M_i\subset C_i$ and $C_n\subset C_{n-1}\cup N_{n-1,n}$.
\eprop
\bpf
We start by defining $C_1=N_1$ which obviously satisfies all the properties. Now suppose that we constructed the required collection up to level $n-1$. Let $\hat C_{n}$ be any characteristic submanifold for $M_n$ and apply Lemma \ref{char-submanifold} to $C_{n-1},\hat C_n$ and $N_{n-1,n}$ to obtain a new characteristic submanifold $C_n$ of $M_n$ such that $C_n\subset C_{n-1}\cup N_{n-1,n}$. Then, for all $j<n:$

$$C_n\cap M_j\subset (C_{n-1}\cup N_{n-1,n})\cap M_j=C_{n-1}\cap M_j\subset C_j$$

By iterating this step the result follows. \epf

We construct the characteristic submanifold $(N,R)$ of $(\overline M,\partial \overline M)$ by picking specific components of the various $(N_i,R_i)$.\ Precisely, we want to pick the components that remain essential throughout the exhaustion, we call these components \textit{admissible}. These will be components $P$ of $N_i$ with enough components of $P\cap R_i$ that generate a product in $M\setminus\text{int}(P)$. That is, if $S\cong\Sigma_{g,n}$ is a component of $P\cap R_i$ then we have a product $\P:\Sigma_{g,n}\times[0,\infty)\hookrightarrow M\setminus\text{int}(P)$ such that $\P(\Sigma_{g,n}\times\set 0)=S\subset\partial M\setminus\text{int}(P)$.

\bdefi
Let $ M=\cup_{i\in\N}  M_i\in\M$ and let $(N_i,R_i)$ be characteristic submanifolds of the $ (M_i,\partial M_i)$'s. We say that an essential submanifold $(P,Q)$ of $(N_i,R_i)$ homeomorphic to a sub-bundle, solid torus or a thickened torus is \emph{admissible} if one of the following holds:
\begin{enumerate}
\item[(i)] $Q$ has two components $A_1,A_2$ that generate in $ M\setminus \text{int}(P)$ a product $\mathcal A$;
\item[(ii)]  $P$ is homeomorphic to an essential solid torus and $Q$ has one component $A$ that generates a product in $ M\setminus \text{int}(P)$ and another component $B$ such that $B$ is the boundary of solid torus $V\subset M\setminus \text{int}(P)$ whose wings wrap $n>1$ times around the soul of $V$;
\item[(iii)] $P$ is a solid torus whose wings wrap $n>1$ times around the soul of $P$ and a component of $Q$ generates in $ M\setminus \text{int}(P)$ a product $\P$;
\item[(iv)] $P$ is homeomorphic to an essential thickened torus.
\end{enumerate}

\edefi

\blem\label{maxadmissible}
Let $M=\cup_{i\in\N}  M_i\in\M$ and let $(N_i,R_i)$ be characteristic submanifolds of the $ (M_i,\partial M_i)$'s. Then, for each $i$ there exists an admissible submanifold $(P_i,Q_i)$ of $(N_i,R_i)$ such that any admissible submanifold of $(N_i,R_i)$ is isotopic into $P_i$.
\elem
\bpf Since $N_i$ has finitely many component and any admissible submanifold of $M_i$ is isotopic into $N_i$ it suffices to work component by component. By Lemma \ref{maxwindows} for every window $W\overset\psi\cong F\times I$ over a hyperbolic surface $F$ we get a maximal submanifold $Q_W\subset W$ such that, up to isotopy, $Q_W$ contains all sub-bundles of $W$ going to infinity. The manifold $Q_W$, up to isotopy, is homeomorphic via $\psi$ to $F_1\times[0,\f1 3]\coprod F_2\times[\f 23,1]$ where $F_1,F_2$ are essential subsurfaces of $F$, which we can assume to be in general position. Then, $\Sigma\eqdef F_1\cap F_2$ is an essential sub-surface of $F$ such that we have a proper embedding $\iota:\Sigma\times\R\hookrightarrow M$ in which $\iota(\Sigma\times[0,1])\subset N_i$, thus it is an admissible submanifold and we denote it by $Q_W$. We then add to $Q_W$ a maximal collection of pairwise disjoint admissible solid tori contained in $W\setminus Q_W$. Note that an essential torus $V\in\pi_0(Q_W)$ can be isotopic into a side bundle of a window $w\in\pi_0(Q_W)$ over a hyperbolic surface. Then to define $P_i$ we do this construction for every window component and take all thickened tori components of $N_i$ and all admissible essential solid tori component. We now show that any admissible submanifolds is isotopic into $P_i$.

Let $P\overset\psi\cong S\times I\subset F\times I$, $\abs{\chi(S)}<0$, be any admissible submanifold a window component $W$, then by Lemma \ref{maxwindows}, up to isotopy, we have that  $S\times\set 1\subset F_2\times\set 1$ and $S\times\set 0\subset F_1\times\set 0$, thus $S$ is isotopic into $F_1\cap F_2$ and hence $P$ is isotopic into $Q_W$. 

If $P$ is an admissible solid torus or thickened torus then by JSJ theory is isotopic in a component of $N_i$ that is either a component of $P_i$ or $P\cong \mathbb S^1\times\mathbb D^2$ is isotopic into an $I$-bundle component $W$ of $N_i$. Since all other cases are contained in $P_i$ by construction we only need to show it for the latter case. Thus, we can assume that we have $P$ isotopic to a vertical thickened annulus $P'$ in a window $W\cong F\times I$ of $N_i$. We need to show that $P$ is isotopic into $Q_W$.  If $P'\cap Q_W=\emp$ and is not isotopic into an essential torus component of $Q_W$ we contradict the maximality of $Q_W$. Therefore, we have that $P'\cap Q_W\neq \emp$. Thus, in the window $W\overset\phi\cong F\times I$ of $N_i$ we have that $P'\overset\phi\cong A\times I$ for $A$ an annulus and we have a component $w\overset\phi\cong S\times I$ of $Q_W$  in which up to isotopy $ A\cap S\neq \emp$ and are in minimal position with respect to each other so that $S\cup A$ is an essential subsurface of $F$. We now need to deal with various cases.

Say that $w$ is admissible and of type (i) so that $\phi(S\times\partial I)$ generate a product $\P: (S\times\partial I)\times[0,\infty)\hookrightarrow M\setminus\text{int}(w)$. If $P$ is also of type (i) we get that it also generates a product $ \mathcal Q: (A\times\partial I)\times[0,\infty)\hookrightarrow M\setminus\text{int}(P)$. Thus, by adjoining $\mathcal Q$ to $\P$ we can enlarge $\P$ to a new product $\P': (S'\times \partial I)\times[0,\infty)\hookrightarrow M\setminus\text{int}(W')$ for $S'$ the essential sub-surface of $F$ filled by $S\cup A$ and $W'\eqdef \phi(S'\times I)$ contradicting the maximality of $Q$.

Now assume that $P$ is of type (ii) so that we have only one component $A_1$ of $\phi(A\times\partial I)$ generating a product while $A_2$ has a root in $M\setminus \text{int} (P)$ which is contained in some solid torus $V$. Since $V$ is compact let $k>i$ be such that $V\subset M_k$ and let $L\subset\partial M_i$ be the component containing $A_2$. Let $M_k'\eqdef \overline{M_k\setminus N_r(L)}$, then $M_k'$ is irreducible, with incompressible boundary and atoroidal thus it has a characteristic submanifold $N$. Moreover, we have that $V\subset N$ and we also have an $I$-bundle induced $J\cong S\times I$ by $\P$ such that they intersect essentially in a component of $\partial M_k'$. Then, we have a component of $N$ containing both an $I$-bundle and a root of its boundary which is impossible. Similarly this takes care of the case in which $W$ is of type (ii) and $P$ of type (i). Thus, we are only left with the case in which both $w$ and $P$ are of type (ii). 

Let $k>i$ be such that $M_k$ contains both roots of the elements of $W\cap \partial M_i$ and $P\cap \partial M_i$ and consider as before $M_k'$ and its characteristic submanifold $N$. Moreover, let $S_1,S_2$ be the surfaces induced by the regular neighbourhood of the component $L$ of $\partial M_i$. Then, we either have two simple closed loops $\alpha,\beta$ both having a root in $M_k'$ such that $\iota(\alpha,\beta)>0$ which cannot happen or we have a component of $N$ containing an $I$-bundle and a root of its boundary which also cannot happen.

Therefore, we get that every admissible submanifold of $N_i$ is indeed isotopic into $Q$ completing the proof.\epf

The following two Lemmas say that essential annuli in $\overline M$ are eventually essential in some $M_i$.

\blem\label{essentialannuli}
Let $M\in\M$ and $C:(\mathbb A,\partial\mathbb A)\hookrightarrow (\overline M,\partial\overline M)$ for $\overline M\in\cat{Bord}(M)$ the maximal bordification. If $C$ is essential in $\overline M$ there exists a minimal $n$ and a proper isotopy of $C$ such that all compact annuli of $\im(C)\cap M_n$ and $\im(C)\cap M\setminus \text{int}(M_n)$ are essential.
\elem
\bpf With an abuse of notation we will use $C$ to denote $\im(C)$. Up to a proper isotopy of $C$ that is the identity on $\partial\overline M$ we can assume that $C\pitchfork\partial M_i$ for all $i\in\N$. Now consider the minimal $i$ such that $M_i\cap C\neq \emp$ and look at the components of $C\cap M_i$. If we have a component $H$ of $C\cap M_n$ that is essential in $M_i$ up to another proper isotopy of $C$ we can push outside $M_i$ all inessential components. Then, by looking at $C\cap M\setminus \text{int}(M_i)$ by a proper isotopy we can push inside $M_i$ all inessential components. Note that by pushing components of $C\cap M\setminus \text{int}(M_i)$ into $M_i$ we might change $H$ to a component $H'$ which is however isotopic to it hence still essential. Since all these isotopies decrease the number of components of $C\cap\partial M_i$ eventually we terminate and all compact components of $C\cap M_i$ and $C\cap M\setminus \text{int}(M_i)$ are essential. Therefore, by picking $n=i$ we are done.

If not it means that all components of $C\cap M_i$ are inessential and via an isotopy $H_i$ of $C$ we can push $C$ outside of $M_i$ so that for all $k\leq i$ we have that $C\cap M_k=\emp$. This process either stops at some $k\geq i$ and by picking $n=k$ we are done by the above case or we obtain a collection of isotopies $\set{H_t^k}_{k\geq i}$ that push $C$ outside every compact subset of $M$ and are the identity on $\partial\overline M$. We will denote by $\hat C$ the properly embedded annulus $C(\mathbb S^1\times (0,1))\subset M=\text{int}(\overline M)$ and without loss of generality we assume that $M_1$ is disjoint from $\hat C$.

\vspace{0.3cm}

\paragraph{Claim:} The annulus $\hat C$ is separating in $ M$.

\vspace{0.3cm}

\bpfc If $M\setminus\hat C$ is connected there exists a loop $\alpha\subset  M$ such that $\alpha\cap \hat C\neq \emp$. Moreover, for any isotopy $H_t$ of $\hat C$ we still have that for all $t:$ $H_t(\hat C)\cap\alpha\neq\emp$. By compactness of $\alpha$ there exists $i$ such that $\alpha\subset M_i$. Then since $\hat C$ can be isotoped outside every $M_i$ we reach a contradiction and so $\hat C$ is separating in $M$. \epfc

Let $E,\hat M$ be the components of $\overline{M\setminus \hat C}$ and assume that $M_1\subset \hat M$. For all $i\in\N$ there is a proper isotopy $H_t^i$ of $\hat C$, namely the one that pushes $\hat C\cap M_i$ outside $M_i$, such that $H_1^i(\hat C)\cap M_i=\emp$. Moreover, we have that $E\iso E_i$ for $E_i$ the component of $\overline{M\setminus \im(H_1^i)}$ not containing $M_1$ and for all $i\in\N$ we have:
\begin{itemize}
\item[(i)] $E_i\cap M_i=\emp$;
\item[(ii)]  $E_{i+1}\subset E_i$;
\item[(iii)] $E_{i+1}\setminus \text{int}(E_i)$ is compact and homeomorphic to a finite collection of solid tori.
\end{itemize}

\vspace{0.3cm}

\paragraph{Claim:} The inclusion $\iota:\hat C\hookrightarrow E$ induces a homotopy equivalence.

\vspace{0.3cm}

\bpfc Since $\hat C$ and $E$ are aspherical by Whitehead Theorem \cite{Ha2002} it suffices to show that the map $\iota$ induces an isomorphism in $\pi_1$. Since $\pi_1(\hat C)$ injects in $M$ we only need to show that $\iota_*$ is a surjection. If $\iota_*$ is not surjective let $\alpha\subset E$ be a non-trivial loop that is not in the image $\iota_*(\pi_1(\hat C))$ and let $i$ be minimal such that $\alpha\subset M_i$. Then, we have a homotopy $\phi_t$ from $\alpha$ into $E_i\simeq E$ and since $H_1^i(\hat C)$ is separating we have that $\alpha$ is homotopic in $E_i$ into $\partial E_i\iso \hat C$ and so
the inclusion map is a homotopy equivalence. \epfc

\vspace{0.3cm}

\paragraph{Claim:} The submanifold $E$ is tame, hence $E\cong V\setminus L$ where $V$ is a solid torus and $L$ is a simple closed curve in $\partial V$.

\vspace{0.3cm}

\bpfc If we show that $E$ is tame, $E\cong V\setminus L$ follows by $\pi_1(E)\cong\Z$ and $\partial E=\hat C\cong \mathbb S^1\times (0,1)$. To show that $E$ is tame we will use that $E\iso E_i$, $E_{i+1}\subset E_i$ and Tucker's Theorem \cite{Tu1974}. To show that $E$ is tame we need to show that for any compact submanifold $K\subset E$ the fundamental group $\pi_1(E\setminus K)$ is finitely generated. Let $i$ be such that $E_i\cap K=\emp$ and so that $E\setminus K=E_i\cup K'$ where by (iii) $K'$ is a compact submanifold of $E$. Then by Van-Kampen's Theorem \cite{Ha2002} we have: 
$$\pi_1(E_i)*\pi_1(K')\twoheadrightarrow\pi_1(E\setminus K)$$ 
and so $\pi_1(E\setminus K)$ is finitely generated. \epfc

Since $E\cong \hat C\times [0,\infty)$ by Theorem \ref{bordification} we have a maximal bordification in which $\hat C$ compactifies to $C'$ and is $\partial$-parallel. Moreover, by uniqueness of the maximal bordification we have that $M'\overset\psi\cong \overline M$ and $\psi$ induces an isotopy from $C$ to $C'$. Contradicting the fact that $C$ was essential in $\overline M$. \epf

By Lemma \ref{maxwindows} we define:

\bdefi\label{boundatinf} Given $M=\cup_{i\in\N} M_i\in\M$ we define the \emph{boundary at infinity} of $M_i$ to be the submanifold $\partial_{\infty} M_i\subset M_i$ to be the maximal, up to isotopy, submanifold of $\partial M_i$ such that we have a simple product $\P:\partial_{\infty} M_i\times[0,\infty)\hookrightarrow M\setminus\text{int}(M_i)$ with the property that every other product $(F\times[0,\infty),F\times\set 0)\hookrightarrow (M\setminus\text{int}(M_i),\partial M_i)$ is isotopic into $\P(\partial_{\infty} M_i\times[0,\infty))$. We also define the \emph{bounded boundary} to be $\partial_{b} M_i\eqdef \overline{\partial M_i\setminus \partial_{\infty}M_i}$.\edefi

\bese
For the manifold $M$ of Example \ref{IWSD} for the elements of the exhaustion $M_i$ we have that $\partial_\infty M_i$ is given by the collection of genus two surfaces corresponding to tame ends and an annulus in the genus two surface facing the non-tame end. The sub-surface $\partial _b M_i$ is given by two punctured tori contained in the genus two surface bounding the non-tame end. 
\eese

We now extend Lemma \ref{essentialannuli} to non-embedded annuli.

\bprop\label{homessentialannuli}  Let $M\in\M$ and $C:(\mathbb A,\partial\mathbb A)\rightarrow (\overline M,\partial\overline M)$ for $\overline M\in\cat{Bord}(M)$ the maximal bordification. If $C$ is essential in $\overline M$ there exists a minimal $i$ and a proper homotopy of $C$ such that all compact components of $\im(C)\cap M_i$ and $\im(C)\cap M\setminus \text{int}(M_i)$ are essential and any $\Z^2\subset \pi_1(\im(C))$ is induced by an annulus in $\im(C)\cap M_i$.
\eprop
\bpf By compactness of the annulus we have a proper homotopy of $C$ in $\overline M$ so that we can assume that $C:(\mathbb A,\partial\mathbb A)\rightarrow (\overline M,\partial\overline M)$ is an immersion that is in general position with $\cup_{k\in\N}\partial M_k$.

\vspace{0.3cm}

\paragraph{Case 1:} Assume that $\pi_1(\im (C))$ does not contain any $\Z^2$, so that, up to homotopy, the singular locus of $C$ does not contain any essential double curve. 

\vspace{0.3cm}

Up to homotopy we can find $i\in\N$ such that $\im (C)\setminus K_i\subset \partial_\infty M_i\times [0,\infty]$ for $K_i$ a compact subset of $\im(C)$ and $\im(C)\cap\partial_\infty M_i\times [0,\infty]=\gamma_1\times[0,\infty]\coprod\gamma_2\times[0,\infty]$ for $\gamma_1,\gamma_2$ two, not necessarily simple, closed curves in $\partial_\infty M_i$.

If the $\gamma_i$ are simple then by the fact that essential annuli in $M_j$ and $\overline{M_j\setminus M_{j-1}}$ with an embedded boundary component are homotopic to embedded essential annuli we obtain a compactly supported homotopy that makes $C$ an embedding. Thus, we are done by Lemma \ref{essentialannuli}. With an abuse of notation we will use $C$ for $\im(C)$ and we will now deal with the case in which the $\gamma_i$ are not simple.

Let $F_i\eqdef Fill(\gamma_i)$ be the essential sub-surface of $\partial_\infty M_i$ filled by $\gamma_i$, for $i=1,2$. If $F_i$ is not an annulus we have that $C\cap M_i$ is essential. Moreover, every compact component of $C\cap M\setminus \text{int}(M_i)$ is also essential, it is induced by a map of an $I$-bundle over the surface $F_i$, and so we are done. Thus, we can assume that $F_i$ is homeomorphic to an annulus. 

If $F_i$ are annuli, we have that $\gamma_i\simeq \alpha_i^{n_i}$ with $\alpha_i$ simple. Then, by the previous argument we have an embedded annulus $C':(\mathbb A,\partial\mathbb A)\hookrightarrow (\overline M,\partial\overline M)$ such that $C$ is properly homotopic into $C'$. By Lemma \ref{essentialannuli} there is a proper isotopy of $C'$ and $i$ such that all compact components of $C'\cap M_i$ and $C'\cap \overline{M\setminus M_i}$ are essential. Thus, since up to a proper homotopy of $C$ it is contained in a thickening of $C'$ the result follows.

\vspace{0.3cm}

\paragraph{Case 2:} Assume that $\pi_1(\im(C))$ contains a $\Z^2$ subgroup $G$. 

\vspace{0.3cm}

Since, $\pi_1(M)=\cup_{i\in\N}\pi_1(M)$ there exists a minimal $i$ such that $\pi_1(M_i)$ contains $G$. By hyperbolicity of the $M_i$ we have that $G$ is conjugated into a subgroup of $\pi_1(T)$ for $T$ a torus in $\partial M_i$. Since, the torus $T$ is compactified in $\overline M$ we have a torus $T_\infty$ such that $T_\infty$, $T$ cobound an $I$-bundle $Q$ in $\overline M$. Moreover, up to an isotopy of $Q$ each $M_j$ intersects $Q$ into a level surface.

 If, up to homotopy, $C\subset Q$ then there is some $M_j$ such that up to homotopy $M_j\cap C$ and the compact components of $\overline{M\setminus M_j}\cap C$ are essential.

If $C$ cannot be homotoped into $Q$ we have a minimal $j\geq i$ such that $C\cap M_j\neq \emp$ and we claim it contains an essential component. If all components of $C\cap M_j$ are inessential we can homotope $C$ such that $C\cap M_j=\emp$ and since $C$ cannot be homotoped into $Q$ we have that $C\subset \overline{M\setminus M_j\cup Q}$. But $\pi_1(\im(C))$ contains $G$ and $G$ is conjugated into $\pi_1(T)$ with $T$ a torus in $\partial M_i$. By tracing the homotopy from $G$ into $\pi_1(T)$, we have a component $S$ of $\partial M_i\setminus T$ that contains a $\Z^2$ subgroup that is homotopic in $M_i$ into $\pi_1(T)$. Thus, we get that $M_i\cong \mathbb T^2\times [-i,i]$ and so $M\cong M\times\R$ in which $M_j\cong \mathbb T^2\times[-j,j]$ and the result follows. Thus, we can assume that $C\cap M_i$ has essential components. Then, as in Lemma \ref{essentialannuli}, up to a homotopy we can assume that $M_i\cap C$ and all compact components of $\overline{M\setminus M_i}\cap C$ are essential.
\epf

In the proof of existence of characteristic submanifolds for manifolds in $M$ we will need the following fact about characteristic submanifold for compact 3-manifolds with incompressible boundary.

\bcor\label{isoIbundle}
If $\psi_n:(F\times I,\partial I)\hookrightarrow (M,\partial M)$, $n=1,2$, are essential $I$-bundles and $M$ is compact, irreducible with incompressible boundary. If $\chi(F)<0$ and $\psi_1(F\times\set0)=\psi_2(F\times\set0)$ then up to isotopy we have that $\psi_1=\psi_2$. If $F$ is an annulus and we have a collection $\set{\psi_n}_{n\in\N}$ then the result is true up to sub-sequence.
\ecor

We will use the following Lemma to show that the submanifold that we build in Proposition \ref{existencechar} contains, up to homotopy, all essential cylinders.
\blem\label{admcomp}
Let $C:(\mathbb A,\partial\mathbb A)\rar (\overline M,\partial \overline M)$ be an essential cylinder such that every compact sub-annulus of $ \im(C)\cap\overline{M\setminus M_i}$ and $\im (C)\cap M_i$ is essential and if $\Z^2\subset\pi_1(\im(C))$ then it is induced by a sub-annulus contained in $M_i$. Then, every $A\subset \im(C)\cap M_i$ is homotopic into an admissible submanifold $P$ of the characteristic submanifold $N_i$ of $M_i$.
\elem
\bpf  We have $0<a<b<1$ such that $A=C(\mathbb S^1\times [a,b])$ is an essential annulus in $M_i$, hence by JSJ theory is homotopic into a component $Q$ of $N_i$. Moreover, we let $C_i$ be the collection of essential annuli induced by $C$ contained in $M_i$.

\paragraph{Case I} Assume that $Q\overset\phi\cong F\times I$ for $F$ a hyperbolic surface. 

Then, up to homotopy, we have that $A\overset\phi\cong \gamma\times I$ for $\gamma\subset F$ a $\pi_1$-injective closed curve. Thus, up to an ulterior homotopy we have that all compact components $C_c$ of $\im(C)\cap\overline {M\setminus M_i}$ are also of the form $\gamma\times I$ and so are the other components of $C_i$. Thus, in $\overline{M\setminus\phi(N_\epsilon(\phi(\gamma\times I))}$ we get two $I$-bundles $P_1,P_2$ such that $C_i\cup C_c\subset P_1\cup P_2\cup \phi(N_\epsilon(\gamma\times I))$. After this homotopy we have that $\im(C)\setminus \text{int}(M_i)$ has two unbounded components $C_1,C_2$ that have for boundary on $\partial M_i$ loops $\alpha_1,\alpha_2$ homeomorphic to $\gamma$. Then, by applying Corollary \ref{isoIbundle} to $\overline{M_n\setminus M_i}$ and a diagonal argument we obtain an embedded product $\hat C_1,\hat C_2\subset \overline{M\setminus M_i}$ such that $\hat C_i\cong \gamma\times [0,\infty)$ and so we get that $N_\epsilon(A)$ is admissible.

 \paragraph{Case II} Assume that $Q\overset\phi\cong \mathbb S^1\times \mathbb D^2$ is a solid torus of type $T_k^n$ of $N_i$ where $n$ is the number of wings and $k$ is the number of times that they wrap around the soul.

If $k>1$ we have that $T_k^n$ is admissible if a component $B$ of $Q\cap\partial M_i$ is isotopic in $M\setminus \text{int}(Q)$ to infinity. Let $A_Q$ be the collection of annuli of $C_i$ homotopic into $Q$. Since $C$ is a proper map we have finitely many such components and thus we can assume that there is $t>0$ such that $C(\mathbb S^1\times [t,1))\cap M_i$ has no components homotopic into $Q$ and $C(\mathbb S^1\times\set t)\subset Q\cap \partial M_i$. Thus we can assume that up to homotopy it is disjoint from $Q$. The loop $C(\mathbb S^1\times\set t)$ is homotopic to $\alpha^m$ for $\alpha$ the core curve of a component $B$ of $V\cap\partial M_i$. Then, by Corollary \ref{isoIbundle} and a diagonal argument we get that $N_\epsilon(\alpha)\subset \partial Q\subset\partial M_i$ generates a product in $M\setminus  \text{int}(V)$ and so $Q$ is admissible.

Similarly if $k=1$ for each $A\subset A_V$ we have that $A\cap \partial Q$ is, up to homotopy, homeomorphic to $\alpha_1^m,\alpha_2^m$ for $\alpha_1,\alpha_2$ core curves of components of $\partial Q\cap \partial M_i$. Since $A_V\subset \im(C)$ is compact we have two components $A_1,A_2\in\pi_0(A_Q)$ such that $\im(C)\setminus A_1\cup A_2$ has two unbounded component $C_1,C_2$ such that $\partial C_1,\partial C_2$ are $\alpha_1^m,\alpha_2^k$ for $\alpha_i$ core curves of components $w_1$, $w_2$ of $\partial Q\cap \partial M_i$. If the components $w_1$, $w_2$ are distinct by using Corollary \ref{isoIbundle} and a diagonal argument we get that $N_\epsilon(\alpha_1\cup \alpha_2)\subset \partial Q\subset\partial M_i$ generates a product in $M\setminus  \text{int}(Q)$ and so $Q$ is admissible.

If $w_1=w_2$ it means that either the annulus $C\setminus C_1\cup C_2$ contains a $\Z^2$ subgroup and so it is contained in $M_i$ and thus $Q$ was not a solid torus or a compact component of $C\cap \overline{M\setminus M_i}$ is contained in a $T_k^n$ torus. Thus, by using Corollary \ref{isoIbundle} and a diagonal argument we get that $w_1\subset \partial Q\subset\partial M_i$ generates a product in $M\setminus  \text{int}(V)$ and so $Q$ is admissible.

\paragraph{Case III} Assume that $Q\cong \mathbb T^2\times I$.

By definition these components are admissible and there is nothing to do. \epf
We can now prove the existence of the characteristic submanifold for bordifications of manifolds in $ \M$.

\bthm[Existence of JSJ]\label{existencechar}  Given $M\in\M$ there exists a maximal bordification $\overline M$ with a characteristic submanifold $(N,R)$.\ethm
\bpf 
 Let $\set{N_i}_{i\in\N}$ be a collection of characteristic submanifold of the $M_i$ coming from Corollary \ref{infjsj} applied to a normal family $\set{N_{i-1,i}}_{i\in\N}$ of characteristic submanifolds for the $X_i\eqdef\overline{M_{i}\setminus M_{i-1}}$. Thus, we can assume that the $(N_i,R_i)\subset ( M_i,\partial M_i)$  satisfy for all $i>j$: $N_i\cap  M_j\subset N_j$ and for all $i$ $N_i\subset N_{i-1}\cup N_{i-1,i}$.

We will construct $N$ as a bordification of a nested union of codimension-zero submanifolds $\hat N_i$. The submanifolds $\hat N_i$ will be obtained by taking admissible submanifolds of the characteristic submanifold $ N_i\subset M_i$ and the $\hat N_i$ will satisfy the following properties:
\begin{enumerate}
	\item[(i)] $\forall j\geq k:\hat N_j\cap  M_k\subset N_k$ is compact;  
	\item[(ii)] if $P\subset N_j$ is an admissible submanifold then, up to isotopy, $P\subset \hat N_j$;
	\item[(iii)] $\forall k\leq j: \hat N_j\cap  M_k=\hat N_k$.
	\end{enumerate}

	\vspace{0.3cm}
	
 Let  $\hat N_1$ be the maximal submanifold of $N_1$ containing all admissible submanifolds, see Lemma \ref{maxadmissible}. Then, $\hat N_1$ clearly satisfies (i)-(iii). We then proceed iteratively. Assume we have constructed $\hat N_i$ and start by defining $\hat N_{i+1}\eqdef \hat N_i$. 

	Let $P\overset\psi\cong F\times I$ be an $I$-bundle component, with $\chi(F)<0$, of $\hat N_i$, then $\psi(F\times\partial I)\subset\partial M_i$. Since $P$ is admissible we have that the surface $\psi(F\times \partial I)$ generates a product: 
	$$\P:(F\times\partial I)\times[0,\infty) \hookrightarrow M\setminus\text{int}(\psi(F\times I))$$
	 such that $\P(F\times\set{0,1})=\psi(F\times\partial I)$. Since $\P(F\times\set{0,1})$ are already essential sub-surfaces of $\partial M_i$ by Theorem \ref{prodstandardform} we have a proper isotopy of $\P$ rel $\P(F\times\set{0,1})$ such that $\P$ is in standard form. Thus, in $X_{i+1}$ there are finitely many essential $I$-bundles $P_1,\dotsc, P_n$ in $N_{i,i+1}$ that connect a component of $\psi(F\times\partial I)$ to either $\partial M_{i+1}$ or to another $I$-bundle $P'$ in $ N_i$. Since $P$ is admissible we have that $P'$ is also admissible and so it is contained in a component $Q$ of $\hat N_i$. Moreover, $Q$ is also homeomorphic to $F\times I$. If not, we would have that $P$ is an essential submanifold of a submanifold of $M$ homeomorphic to $F'\times\R$ where up to isotopy $P\overset\phi\cong F\times [0,1]\subset F'\times [0,1]$ is a sub-bundle. Thus, we would have that $P$ is contained in a larger admissible submanifold of $\hat N_i$ contradicting the construction.

	By adding all such $P_\ell$'s to $\hat N_{i+1}$ and repeating it for all such $I$-bundles we have that $\hat N_{i+1}$ satisfies (i) and (iii) by construction.

Let $P\overset \psi\cong \mathbb S^1\times\mathbb  D^2$ a solid torus component of $\hat N_i$. As before, we add to $\hat N_{i+1}$ all solid tori components of $N_{i,i+1}$ and thickened annuli contained in $I$-bundles of $N_{i,i+1}$ that match up with component of $P\cap\partial M_i$. Properties (i) and (iii) are still satisfied by construction. Similarly we do the case where $P\overset \psi\cong\mathbb T^2\times I$. We now claim that the only admissible submanifolds of $N_{i+1}$ that we are missing in $\hat N_{i+1}$ are contained in $N_{i,i+1}$

\vspace{0.3cm}

\paragraph{Claim:} If $Q\subset N_{i+1}$ is admissible and, up to isotopy, $Q\cap M_i\neq \emp$ then we have that $Q$ is isotopic into $\hat N_{i+1}$.
	
\vspace{0.3cm}

	\bpfc
	Let $Q$ be such a component then by Corollary \ref{infjsj} we have that, up to isotopy, $Q\subset N_i\cup N_{i,i+1}$ and $Q\cap N_i$ is an essential submanifold. Since $Q$ is admissible in $N_{i+1}$ and $Q_i\eqdef Q\cap N_i$ is an essential submanifold we have that $Q_i$ is also admissible in $N_i$. Therefore, $Q_i$ is, up to isotopy, contained in $\hat N_i$. Hence, by the above construction we get that $Q\subset \hat N_{i+1}$.
	\epfc

	Finally, we add to $\hat N_{i+1}$ all admissible solid tori, thickened essential tori and $I$-subbundles contained in $N_{i,i+1}$.

	 By construction we have that $\hat N_{i+1}\cap M_i=\hat N_i$ and $\hat N_i$ is compact thus (i) and (iii) are satisfied. Moreover, for all $i\in\N$ all components of $\hat N_i$ are $I$-bundles over hyperbolic surfaces, solid tori or thickened essential tori. Let $P\subset N_{i+1}$ be admissible then, up to isotopy, we have that $P\subset \hat N_{i+1}$ and so $\hat N_{i+1}$ satisfies (i)-(iii). 
	 
Since by construction $\hat N_i$ does not change as we go through the construction we obtain a collection $\set{\hat N_i}_{i\in\N}$ of nested codimension-zero submanifold satisfying (i)-(iii).

Let, $\hat N\eqdef \cup_{i=1}^\infty \hat N_i\subset  M$. Since every component of $N$ has a natural 3-manifold structure and for all $k\in\N$:
$$\hat N\cap M_k=\cup_{i=1}^\infty \hat N_i\cap M_k\overset{(\text{ii})}=\hat N_k$$

 is compact we get that $N$ is a properly embedded codimension-zero submanifold. Moreover, $\hat N$ contains, up to isotopy, all admissible submanifolds since they appear in some $\hat N_j$.

	Let $P\in\pi_0(\hat N)$ then, by construction, $P$ is either:
	\begin{itemize}
	\item homeomorphic to $F\times \R$ for $F$ a compact surface;
	\item  a nested union of solid tori $T_i\subset N_i$ such that $\overline{T_{i+1}\setminus T_i}$ are essential solid tori;
	\item  a nested union of manifolds $Q_i\subset N_i$ each homeomorphic to $\mathbb T^2\times [0,1]$ and such that $\overline{Q_{i+1}\setminus Q_i}$ are essential solid tori. 
	\end{itemize}

	If $P$ is the limit of solid tori $T_i$ then by Lemma \ref{charlimitsolidtori} we have that $P\cong V\setminus L$ for $V$ a solid torus and $L$ a closed collection of parallel loops in $\partial V$. Similarly, if $P$ is the limit of thickened essential tori: $\mathbb T^2\times [0,1]$ we get by Corollary \ref{charesstori} that $P$ is homeomorphic to $\mathbb T^2\times [0,1]\setminus L$ for $L\subset\mathbb T^2\times \set 0$ a collection of parallel loops.

We will now add boundary to $\hat N$. Let $P\overset\psi\cong F\times\R$ be an $I$-bundle component of $\hat N$ since $P$ is properly embedded by adding $\text{int}(F)\times\set{\pm\infty}$ to $M$ we can compactify $P$ to $\overline P$ in $\hat M\in\cat{Bord}(M)$ so that $\overline P\overset{\bar\psi}\cong F\times I$ is an essential $I$-bundle in $\hat M$. By repeating this for all components of $\hat N$ homeomorphic to $F\times\R$ we obtain a new manifold, which we still denote by $\hat N$, properly embedded in $\hat M$ such that all $I$-bundles components are essential and compact.

Let $P\cong V\setminus L$  or $P\cong\mathbb T^2\times [0,1)\setminus L$ and consider the subset of loops $L^{iso}\subset L$ that are not accumulated by any family of loops $\gamma_i\rar\gamma$. Since each $\gamma\in L^{iso}$ is isolated it means that if we take a closed end neighbourhood $U$ of $\gamma$ it is homeomorphic to a properly embedded annular product $\mathbb A\times [0,\infty)$ which we can compactify in $\hat M$ by adding an open annulus to $\partial \hat M$. Moreover, for all components homeomorphic to $\mathbb T^2\times[0,1)\setminus L$, $L\subset \mathbb T^2\times\set 0$, we add the corresponding boundary torus $\mathbb T^2\times\set 1$ to $\hat M$. We still denote by $\hat N$ the resulting bordified manifold contained in $\hat M\in \cat{Bord}(M)$.

Finally, for each $\gamma$ in the set $L'\eqdef L\setminus L^{iso}$ of an IWSL or an IWET by tracing through the gluing of the solid tori we obtain an embedded product $\P_\gamma:\mathbb A\times [0,\infty)\hookrightarrow \hat M$. We also compactify $\P_\gamma$ by adding an open annulus $A_\gamma$ to $\partial\hat M$ so that $\P_\gamma$ is partially compactified i.e. : 
$$\overline \P_\gamma: (\mathbb A\times [0,\infty)\cup A'\times\infty,A'\times\infty)\hookrightarrow (\hat M\cup A_\gamma, A_\gamma)$$
where $A'\subsetneq\mathbb A$ is an annulus sharing only one boundary component with $\mathbb A$. Moreover, we have that $\P_\gamma$ is properly isotopic into a collar neighbourhood of $A_\gamma$.

Then, we extend the bordification $\hat M\in\cat{Bord}(M)$ to obtain a maximal bordification $\overline M$, see Theorem \ref{bordification}. Finally from $N$ we remove any essential solid torus $T$ that is properly homotopic into the side boundary of an $I$-bundle component.

To show that $N$ is a characteristic submanifold we need to show that any essential annulus $\mathcal A: (\mathbb A,\partial\mathbb A)\rar (\overline M,\partial\overline  M)$ and essential torus $\mathcal T:\mathbb T^2\rar \overline M$ is homotopic into $N$. We first show that annuli can be homotoped into $N$. 

By Proposition \ref{homessentialannuli} for any essential annulus $A\eqdef \mathcal A(\mathbb A)$ in $\overline M$ we have a proper homotopy of $A$ and $i\in\N$ such that all compact sub-annuli of $A_i\eqdef M_i\cap A$ and $A\cap \overline{M\setminus M_i}$ are essential in $M_i$, $\overline{M\setminus M_i}$ respectively. Thus, by Lemma \ref{admcomp} we have that $A_i$ it is homotopic not just into $N_i$ but into $\hat N_i$. Hence, we can assume that $A_i\subset\hat N_i$.

Since all compact components $A_{c}$ of $A\cap \overline{M\setminus M_i}$ are essential there is a proper homotopy of $A_c$ rel $\partial M_i$ such that for all $k>i$ $A_c\cap X_k\subset N_{k,k-1}$. Moreover, since all components $Q$ of $N_{k,k-1}$ containing sub-annuli of $A_c$ match up with admissible components of $N_i$ we get that $Q\subset \hat N_k$. Thus, we have a proper homotopy such that $A_i\cup A_c\subset N$.

We will now do an iterative argument to construct homotopies rel $\partial M_{j-1}$, $j>i$, supported in $M\setminus \text{int}(M_{j-1})$ such that $A\cap X_j\subset \hat N_j\cap X_j$. 

\vspace{0.3cm}

\paragraph{Claim:} If for $n<j$ we have that $A\cap X_n\subset  \hat N_n\cap X_n$ and all annuli of $A\cap X_j $ with boundary on $\partial M_{j-1}$ are essential. Then, there is a proper homotopy rel $\partial M_{j-1}$ supported in $M\setminus \text{int}(M_{j-1})$ such that  $A\cap \overline{M\setminus M_{j-1}}$ are essential and contained in $ \hat N_j\cap X_j$ and all annuli of $A\cap \overline{M\setminus M_j}$ with boundary on $\partial M_{j}$ are essential. 

\vspace{0.3cm}

\bpfc All annuli of $A\cap X_j$ that are $\partial$-parallel are induced by the unbounded components $A_1,A_2$ and have boundary on $\partial M_j$. Thus, by a proper homotopy $\phi_1$ supported in $M\setminus \text{int}(M_{j-1})$ we can remove the $\partial$-parallel annuli and guarantee that all compact components of $(A_1\cup A_2)\cap \overline{M\setminus M_j}$ with boundary on $M_{j}$ are essential so that all annuli of $A\cap \overline{M\setminus M_j} $ with boundary on $\partial M_{j}$ are essential.

Since every annulus of $A\cap X_j$ is essential we have a homotopy $\phi_2$, supported in $X_j$, that by properties of normal family is the identity on $\partial M_{j-1}$ such that $A\cap X_j\subset N_{j,j-1}$. If a sub-annulus $C$ of $A\cap X_j$ satisfies $\partial C\subset \partial M_j$ then by Lemma \ref{admcomp} we have that $A$ is contained in a n admissible component and so $A\subset \hat N_j$. If $\partial A$ has a component contained in $\partial M_{j-1}$ by properties of normal families and the fact that for all $n<j$ $A\cap X_n\subset  \hat N_n\cap X_n$ we get that $A$ matches up with an admissible component and is so admissible.\epfc

The composition $\phi\eqdef \varinjlim_{j\geq i}\phi_j$ gives a proper homotopy of $A$ such that $A\subset \hat N$.

Finally, let $\mathcal T:\mathbb T^2\rar \overline M$ be an essential torus, then by a homotopy we can assume that $\im(\mathcal T)\cap\partial\overline M=\emp$ and by compactness of $\im (\mathcal T)$ we have $M_i$ such that $\im(\mathcal T)\subset M_i$. Since the $M_i$ are atoroidal we have that $\im(\mathcal T)$ is homotopic into a torus component of $\partial_\infty M_i$ which is isotopic into a torus component $T_\infty$ of $\partial\overline M$. Hence, we have a homotopy from $\im(\mathcal T)$ into the component of $N$ corresponding to $T_\infty$. 

\vspace{0.3cm}

\paragraph{Claim:} The manifold $N$ is minimal, that is no component is homotopic into another.

\vspace{0.3cm}

Assume that $P,Q\in\pi_0(N)$ are such that $P$ is properly homotopic into $Q$. By fundamental group reasons we get that $P\cong F\times I$ and $Q\cong\mathbb S^1\times \mathbb D^2$ and that $Q$ is properly homotopic into a side boundary of $P$ and we removed all these redundancies.\epf

\subsubsection{Uniqueness of characteristic submanifolds}
 
 We will now show that any characteristic submanifold of $\overline M$ can be put in a \emph{normal form} so that they are contained in a pre-scribed family of a normal characteristic submanifolds for the gaps of the exhaustion $\set{M_i}_{i\in\N}$. We will use this fact to show that any characteristic submanifold for $\overline M$ is properly isotopic to the one constructed in Proposition \ref{existencechar}.
 
\bdefi
Let $M=\cup_{k\in\N}M_k\in\M$ and $\iota:(N,R)\hookrightarrow (\overline M,\partial\overline M)$ be a characteristic submanifold for the maximal bordification $\overline M$. Let $N'\eqdef \iota(N\setminus R)$ then $N$ is in \emph{pre-normal form} if every component of $\overline{N'\setminus \cup_{k\in\N}\partial M_k}$ is an $I$-bundle, a solid torus or a thickened torus that $\pi_1$-injects in $N'$ and every component of $N'\cap\cup_{k\in\N}\partial M_k$ is a $\pi_1$-injective surface and no component is a $\partial$-parallel annulus or a disk. 

Given a normal family of characteristic submanifolds $\set{N_k}_{k\in\N}$ we say that $N$ is in \emph{normal form} with respect to the $N_k$'s if for all $k\in\N$ we have that each component of $N'\cap \overline{M_k\setminus M_{k-1}}\subset N_k$ is an essential submanifold of $N_k$.
\edefi 
\bese
Note that by construction the characteristic submanifold constructed in Proposition \ref{existencechar} is in normal form.
\eese
\brem\label{pushparallelwing}
The difference between pre-normal form and normal form is that if $N$ is in pre-normal form but not in normal form then there exists a $k\in\N$ and a component $Q$ of $N'\cap M_k$ such that $Q$ is homeomorphic to either a solid torus or a thickened torus and it has at least two parallel wings. Equivalently it means that an annular component of $\partial Q\setminus \partial M_k$ is $\partial$-parallel. 

Moreover, if two wings of a solid torus $Q\subset M_k$ or $Q\subset \overline {M\setminus M_k}$ are parallel by an isotopy supported in $X_k\cup X_{k+1}$ we can slide one over the other and push it in $\overline {M\setminus M_k}$ or $M_k$ respectively.

\erem

We now prove a Lemma needed to show that characteristic submanifolds can be put in pre-normal form.

\blem\label{partialparallelannulichar}
Let $\iota:(N,R)\hookrightarrow (\overline M,\partial\overline M)$ be a characteristic submanifold and let $N'\eqdef N\setminus R$. Let $\mathcal S\eqdef \im(\iota)\cap\cup_{k\in\N}\partial M_k$ and assume that every component of $\mathcal S$ is $\pi_1$-injective and no component of $\mathcal S$ is a disk. Then, there is a proper isotopy $\psi_t$ of $\iota$ supported in $\iota(N')$ such that no component $S$ of $\mathcal S$ is a boundary parallel annulus in $\iota(N')$.

\elem
\bpf 
Let $\mathcal A_k$ be the collection of annuli of $\mathcal S_k\eqdef \mathcal S\cap\partial M_k$ that are $\partial$-parallel in $\iota(N')$. Since $\iota$ is a proper embedding we have that for all $k\in\N$ $\abs{\pi_0(\mathcal A_k)}<\infty$. Moreover, since $N$ is a characteristic submanifold every $\partial$-parallel annulus $A\subset \mathcal A_k$ is contained in a component of $N'$ homeomorphic to either an $\R$-bundle, a missing boundary solid torus $V$ or thickened essential torus $T$.  By an iterative argument it suffices to show the following:

\vspace{0.3cm}

\paragraph{Claim:} If for $1\leq n< k$ $\mathcal A_n=\emp$ then via an isotopy $\phi_k^t$ of $\iota$ supported in $\overline{M\setminus M_{k-1}}\cap\im(\iota)$ we can make $\mathcal A_{k}=\emp$.

\vspace{0.3cm}

\bpfc 
Denote by $A_1,\dotsc, A_n$ the $\partial$-parallel annuli in $\mathcal A_k$ and assume that $A_1,\dotsc , A_{n_1}$ are contained in $\R$-bundle components of $N'$ and $A_{n_1+1},\dotsc, A_n$ are contained in missing boundary solid tori or thickened tori.

Since the annuli contained in $\R$-bundles are finitely many we have a disconnected compact horizontal surface $F_k$ in $N'$ such that $A_1,\dotsc, A_{n_1}$ are contained in $\iota(F_k\times [a_k,b_k])$. By applying Corollary \ref{pushanninter} to each component of $F_k\times[a_k,b_k]$ we have a local isotopy $\phi_k^t$ of $\iota$ that removes all these intersections. The isotopy $\phi_k^t$ is supported in a collection of solid tori $\mathcal V_k'\subset F_k\times[a_k,b_k]$ thus it can be extended to the whole of $\iota$. Moreover, if we consider for $n<k$ a component of intersection of $\partial M_n\cap \P(\mathcal V_k)$ then it is either a boundary parallel annulus or a disk. However, we assumed that for $n<k$ $\mathcal A_n=\emp$ and by hypothesis no component of $\cup_{k\in\N}\partial M_k\cap\im(N')$ is a disk thus, the solid tori $\iota(\mathcal V_k')$ that we push along are contained in $\iota(N')\cap \overline{M\setminus M_{k-1}}$.

Similarly, consider the annuli $A_{n_1+1},\dotsc, A_n$ then by Lemma \ref{outermost} we have a collection of solid tori $\mathcal V_k''\subset N'$ such that by pushing along $\iota(V''_k)$ we obtain an isotopy $\phi_k^t$ of $\iota$ so that $\mathcal A_k=\emp$ and as before the solid tori $\iota(\mathcal V_k'')$ are contained in $\iota(N')\cap \overline{M\setminus M_{k-1}}$ thus $\phi_k^t$ is supported in $\overline{M\setminus M_{k-1}}\cap\im(\iota)$. Moreover, note that $\mathcal V'_k\cap\mathcal V''_k=\emp$ since they are contained in disjoint components of $N$.

Therefore, we get a collection of solid tori $\mathcal V_k\eqdef \mathcal V_k'\coprod \mathcal V_k''$ contained in $\im(\iota)\cap\overline{M\setminus M_{k-1}}$ such that pushing through them gives us an isotopy $\phi_k^t$ of $\P$ that makes $\mathcal A_k=\emp$. \epfc

Since for all $k\in\N$ $\text{supp}(\phi_k^t)=\mathcal V_k$ is contained in $\overline{M\setminus M_{k-1}}$ the limit $\phi^t$ of the $\phi^t_k$ gives us a proper isotopy of $\P$ such that for all $k\in\N$ $\mathcal A_k=\emp$. \epf

By the Lemma we have:

\bprop\label{prenormalform} Given a characteristic submanifold $N$ of the maximal bordification $\overline M$ of $M\in\M$ there is a proper isotopy such that $N$ is in pre-normal form.
\eprop
\bpf Since $N\cap\text{int}(\overline M)$ is a $\pi_1$-injective submanifold of $M\cong \text{int}(\overline M)$ by Lemma \ref{esssubsurface} we have a proper isotopy of $N$ such that for all $k\in\N$ $\partial M_k\cap N$ are $\pi_1$-injective surfaces in $N$ and no component of $\mathcal S\eqdef \im(\iota)\cap\cup_{k\in\N}\partial M_k$ is a disk.  

Then, by Lemma \ref{partialparallelannulichar} we have a proper isotopy of $N'$ such that no component of $\mathcal S$ is a $\partial$-parallel annulus. Therefore, for all components $S$ of $\mathcal S$ the surface $\iota^{-1}(S)$ is an essential surface in $N'$. By the proof of Lemma \ref{wald} we get that up to a proper isotopy of $N'$ supported in the $\R$-bundle components every essential surface $\iota^{-1}(S)$ in an $\R$-bundle component is horizontal. Thus, $\R$-bundles components of $N'$ are decomposed by $\cup_{k\in\N}\partial M_k$ into $I$-bundles contained in $X_k\eqdef\overline{M_k\setminus M_{k-1}}$.

Let $\mathcal S'\subset\mathcal S$ be the collection of components $S$ of $\mathcal S$ such that $\iota^{-1}(S)$ is not contained in an $\R$-bundle. Each component $S$ of $\mathcal S'$ is either an essential annulus or an essential torus. Since all essential tori are contained in products $\mathbb T^2\times [0,\infty)$ by the proof of Lemma \ref{wald} we have a proper isotopy supported inside a, possibly infinite, collection of products over essential tori such that the pre-images under $\iota$ are essential tori of $\mathcal S'$ which co-bound $I$-bundles. 

Let $\mathcal A_k\subset S'$ be the collection of essential annuli of $\mathcal S'\cap X_k$. Then $\iota^{-1}(\mathcal A_k)$ are essential annuli contained in a component $Q$ of $N'$ that is homeomorphic to either a missing boundary solid torus $V$ or a missing boundary thickened essential torus $T$. In either case they co-bound either a solid torus or a thickened torus and the Lemma follows.\epf

We now show that characteristic submanifolds in pre-normal form can be isotoped to be in normal-form. To prove the iterative step we need:

\blem\label{normalformindstep}
Let $M=\cup_{k\in\N}M_k\in\M$ and $\iota:(N,R)\hookrightarrow (\overline M,\partial\overline M)$ be a characteristic submanifold in pre-normal form for the maximal bordification $\overline M$ and let $N'\eqdef \iota(N\setminus R)$. Let $\mathcal A_n$ be the collection of annuli of $\overline{\partial N'\setminus \partial M_n}$ that are $\partial$-parallel in either $M_n$ or $M\setminus \text{int}(M_n)$. If for all $1\leq n <k$ we have that $\mathcal A_n=\emp$ then there is a proper isotopy $\Psi_k^t$ of $\iota$ supported in $\overline{M\setminus M_{k-1}}$ such that $\mathcal A_k=\emp$.\elem

\bpf We have that $\abs{\mathcal A_m}$ is bounded by $b_m\eqdef \abs{\pi_0(\partial M_m\cap\partial N')}$ which is finite by properness of the embedding.

\vspace{0.3cm}

\paragraph{Claim:} Let $A$ be an annulus in $ \mathcal A_n$, for $n\in\N$ such that: 
\be A\cap\cup_{m=1}^n \partial M_m=\partial A\subset\partial M_n\ee
and $A$ is inessential. Then, there exists a solid torus $V\subset \overline{M\setminus M_{n-1}}$ containing $\iota(A)$ such that all components of $\iota(N)\cap V$ are inessential solid tori in either $M_n$ or $\overline{M\setminus M_n}$.

\vspace{0.3cm}

\bpfc The annulus $\iota(A)$ is $\partial$-parallel so it co-bounds with an annulus $C\subset \partial M_n$ a solid torus $V\subset  \overline{M\setminus M_{n-1}}$. Consider a component $Q$ of $\iota(N)\setminus A\cap V$, then $Q\cap \partial M_n\neq\emp$ and since $Q\cap \partial M_n\subset C$ it is a collection of annuli $B$. Since $\iota(N)\cap \partial M_n$ has no $\partial$-parallel annuli we have that all the annuli $B$ are essential in $\iota(N)$ hence since $Q\subset V$ it must be a solid torus contained in a component of $N'$ homeomorphic to either a solid torus or a thickened essential torus, both potentially missing boudary. Moreover, since $Q\subset V$ and $\partial V\cap\partial M_n$ is an annulus we have that $Q$ is inessential in either $M_n$ or $\overline{M\setminus M_n}$.\epfc

Let $A$ be an element of $\mathcal A_k$ and assume that $A\subset M_k$ and that $A\cap \partial M_{k-1}\neq\emp$. Denote by $Q$ the component of $N\cap M_k$ containing $A$. Since $A$ is $\partial$-parallel we have that any essential sub-annulus $A'\subset A$ with boundaries on $\partial M_{k-1}$, such that $A'\subset M_{k-1}$ is also $\partial$-parallel. Since $A\cap\partial M_{k-1}\neq \emp$ and $A\pitchfork \partial M_{k-1}$ the annulus $A$ is decomposed by $\partial M_{k-1}$ into annuli $A_1^{k-1},\dotsc, A_h^{k-1}$ such that $ A_j^{k-1}\cap \partial M_{k-1}\subset A_j^{k-1}\cap \partial M_{k-1}$ and for $1<j<h$ $A_j^{k-1}$ has both boundary components on $\partial M_{k-1}$. Moreover, since $h\geq 3$ and $A$ is $\partial$-parallel in $M_k$ we get that $\mathcal A_{k-1}\neq\emp$ reaching a contradiction. Thus, we can assume that any annulus $A$ satisfies: 
\be A\cap\cup_{m=1}^k \partial M_m=\partial A\subset\partial M_k\ee

and is inessential. By the Claim we have a solid torus $V$ such that $V\cap \iota(N)$ are inessential solid tori $Q_1,\dotsc, Q_{m_k}$. By a proper isotopy of $\iota(N)$ supported in $V\subset \overline{M\setminus M_{k-1}}$ we can push all inessential tori $Q_1,\dotsc, Q_{m_k}$ contained in $V$ in either $M_k$ or $\overline{M\setminus M_k}$ and in either case we reduce $b_k$ by at least $2m_k> 0$. Thus, we can assume that $V\cap \iota(N)=A$ and $A$ is $\partial$-parallel. Let $Q$ be the solid torus or thickened torus of $\iota(N)\cap M_k$ or $\iota(N)\cap \overline{M\setminus M_k}$ containing $A$. Then by Remark \ref{pushparallelwing} we have a proper isotopy supported in $\overline {M\setminus M_{k-1}}$ that reduces $b_k$ by at least one and removes $A$ from $\mathcal A_k$. Finally, since $\mathcal A_k$ has finitely many elements the composition of these isotopies gives us a proper isotopy $\Psi^t_k$ of $\iota$ such that $\mathcal A_m=\emp$ for $m\leq k$. Moreover, since all the isotopies are supported in $\overline {M\setminus M_{k-1}}$ we get that $\Psi^t_k$ is also supported outside $M_{k-1}$ . \epf

\bprop\label{normalform}Let $\iota:(N,R)\hookrightarrow (\overline M,\partial\overline M)$ be a characteristic submanifold for the maximal bordification $\overline M$ of $M=\cup_{i\in\N}M_i\in\M$. Given a normal family of characteristic submanifolds $\set{N_i}_{i\in\N}$ for $X_i\eqdef \overline{M_i\setminus M_{i-1}}$, there is a proper isotopy of $\iota$ such that $N$ is in normal form.\eprop
\bpf By Lemma \ref{prenormalform} we can assume that $N$ is in pre-normal form and let $N'\eqdef\iota( N\setminus R)$. 

\vspace{0.3cm}

\paragraph{Step 1:} Up to a proper isotopy we have that for all $i\in\N$ $\iota(N')\cap X_i$ is a collection of essential, pairwise disjoint $I$-bundles, solid tori and thickened tori.
\vspace{0.3cm}

Since $\iota(N')\cap X_i$ is in pre-normal form we have that $\iota(N')\cap X_1$ is a collection of $I$-bundles, solid tori and thickened tori. 

Let $\A_1$ be the collection of annuli of $\iota(\partial N')\cap X_1$ and $\iota(\partial N')\cap \overline {M\setminus M_1}$ that are $\partial$-parallel. By properness of the embedding we have that $\mathcal A_1$ has finitely many components. Then, by Lemma \ref{normalformindstep} we get a proper isotopy $\Psi_1^t$ such that all annuli of $\iota(\partial N')\cap M_1$ and $\iota(\partial N')\cap  \overline {M\setminus M_1}$ are essential. Therefore, by Remark \ref{pushparallelwing} we have that all components of $\iota(N')\cap X_1$ are essential.

We now proceed iteratively. Assume that we made for $1\leq n< k$ all annuli $Q\in \mathcal A_n$ essential. Then, by applying Lemma \ref{normalformindstep} to $\mathcal A_{k}$ we obtain a proper isotopy $\Psi_{k}^t$ supported in $\overline {M\setminus M_{k-1}}$ that makes all annuli $Q\subset\iota(\partial N')\cap M_k\cup\iota(\partial N')\cap\overline {M\setminus M_k} $. In particular we get that for all $1\leq n\leq k$ all components of $\iota( N')\cap X_n$ are essential.

Since the isotopies $\Psi_k^t$ are supported in $\overline{M\setminus M_{k-1}}$ the limit composition is a proper isotopy $\Psi^t\eqdef\varinjlim_{k\in\N}\Psi^t_k$ of $\iota$ such that for all $k\in\N$ every component of $N'\cap X_k$ is either an essential $I$-bundle, essential solid torus or an essential thickened torus. 

\vspace{0.3cm}

\paragraph{Step 2:} Up to a proper isotopy we have that $N$ is in normal form.

\vspace{0.3cm}

By \textbf{Step 1} we have that for all $i\in\N$ the components of $ N'\cap X_i$ are a collection of essential, pairwise disjoint $I$-bundles, solid tori and thickened tori. Consider $X_1=M_1$ then by JSJ theory we can isotope $N'\cap X_1$ so that $ N'\cap X_1\subset N_1$. Moreover, since $ N'\cap \partial M_1$ is, up to isotopy, contained in both $R_1$ and $R_2$\footnote{We remind the reader that $R_i\eqdef\partial N_i\cap \partial X_i$.} by definition of normal family we can assume that $ N'\cap \partial M_1$ is contained in $R_{1,2}\eqdef R_1\cap R_ 2$. This isotopy is supported in a neighbourhood of $X_1$, hence it can be extended to a proper isotopy $\Psi_1^t$ of $ N'$. We will now work iteratively by doing isotopies relative $R_{k,k+1}\eqdef R_k\cap R_ {k+1}$. 

Assume that we isotoped $ N'$ such that for all $1\leq n\leq k$ we have that $ N'\cap X_n\subset N_n$ and such that $ N'\cap \partial M_n$ is contained $R_{n,n+1}$. Since the components of $ N'\cap X_{k+1}$ are essential, pairwise disjoint $I$-bundles, solid tori and thickened tori of $X_{k+1}$ with some boundary components contained in $R_{k,k+1}$ we can isotope them rel $R_{k,k+1}$ inside $N_{k+1}$ so that their boundaries are contained in $R_{k,k+1}\coprod R_{k+1,k+2}$. This can be extended to an isotopy $\Psi_{k+1}^t$ of $\P$ whose support is contained in $\overline{M\setminus M_k}$, hence the composition of these isotopies gives a proper isotopy of $N'$ such that $\forall i\in\N: N'\cap X_i\subset N_i$, thus completing the proof.  \epf

We now show that characteristic submanifolds are unique up to isotopy.

\bprop If $N$ and $N'$ are two characteristic submanifolds for $\overline M\in\cat{Bord}(M)$, $M\in\M$, then they are properly isotopic.
\eprop
\bpf It suffices to show that any characteristic submanifolds $N'$ is properly isotopic to the one constructed in Theorem \ref{existencechar}, which is in normal form. Say we have another submanifold $(N',R')\subset (\overline M,\partial\overline M)$ satisfying the same properties of $N$ but not properly isotopic to it. By applying Proposition \ref{normalform} to $N'$ we get that up to proper isotopy we can assume that for all $i$ $N'\cap X_i\subset N_{i-1,i}$.

By definition each component of $N'\cap X_i$ is admissible hence it is isotopic into a component of $N\cap X_i$. Then, by an iterative argument and properties of a normal family we can isotope $N'\cap X_i$ into $N\cap X_i$ rel $\partial M_{i-1}$ to get a proper isotopy of $N'$ into $N$. Moreover, $N'\subset N$ has to be equal to $N$, up to another proper isotopy, since otherwise $N'$ does not satisfy the engulfing property.\epf

With this we conclude the proof of Theorem:

\bthm\label{exisjsj} Given $M\in\M$ there exists a maximal bordification $\overline M$ with a characteristic submanifold $(N,R)$. Moreover, any two characteristic submanifolds are properly isotopic.
\ethm

  \subsection{Acylindricity Conditions}\label{acylindricityconditions}
 
From now one we will focus on manifolds in $\M^B$ so that each element of the exhaustion of $M$ has a bound on the genus of its boundary components and the same holds for the components of $\partial\overline M$. Given the existence of characteristic submanifolds we can now study acylindricity properties of manifolds in $\M^B$. In particular we want to construct a system of simple closed curves $P$ in $\partial \overline M$, for $M\in\M^B$, such that $\overline M$ is acylindrical with respect to $P$. 
 	
	\bdefi
	We say that an irreducible 3-manifold $(M,\partial M)$ with incompressible boundary is acylindrical rel $P\subset \partial M$ if $M$ has no essential cylinders $C$ with boundary in $\partial M\setminus P$.
	\edefi

 The example of Section \ref{lochyp} shows that not all manifolds in $\M^B$ admit such a system, but the existence of a doubly peripheral annulus is the only obstruction. To check wether a manifold $M\in\mathcal M^B$ has a doubly peripheral annulus it suffices to check that the characteristic submanifold $(N,R)$ of $\overline M$ does not contain any essential annulus $\mathcal A: (\mathbb A,\partial \mathbb A)\rar (\overline M,\partial \overline M)$ such that $\mathcal A(\partial \mathbb A)$ are both peripheral in $\partial\overline M$.

	\begin{center}\begin{figure}[h!]
		\centering
		\def\svgwidth{300pt}
		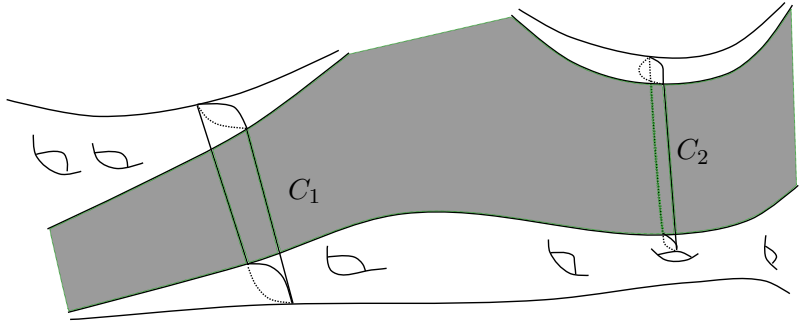\caption{The annulus $C_1$ is an example of a doubly peripheral cylinder while $C_2$ is not since $\partial C_2$ has only one peripheral component in $\partial\overline M$. The interior of $\overline M$ is shaded.}
	\end{figure}\end{center}

  We will show that manifolds $M\in\M^B$ without doubly peripheral annuli have a system of simple closed curves $P$ that make $\overline M$ acylindrical relative to $P$. Once we add all the torus and annular components of $\partial\overline M$ to $P$ we get that $(\overline M,P)$ becomes what is known as a pared manifold (see \cite{Th1986} and definition \ref{paredinfinite}).

In order to find a collection of curves $P$ in $\partial  \overline M$ such that $\overline M$ is acylindrical relative to $P$ it suffices to show that the curves in $P$ `pierce' all the cylinders of the characteristic submanifold $(N,R)\subset (\overline M,\partial\overline M)$.  Before we go on with the construction, we make the following remarks:
 \begin{enumerate}
 	\item[(i)] since we assumed that $M$ has no doubly peripheral cylinder every essential cylinder $C\subset\overline M$ must have at least a non-peripheral boundary component;

 	\item[(ii)] in order to construct $P$ it suffices to find simple closed curves $\set{\gamma_i}_{i\in\mathcal I}\subset \partial M$ such that for any closed curve $\alpha\subset R$ there is $i\in\mathcal I$ such that $i(\gamma_i,\alpha)>0$;
 	\item[(iii)] given a solid torus component $V$ of the characteristic submanifold $N$ we cannot have more than one wing being peripheral in $\partial \overline M$ otherwise we can find a doubly peripheral cylinder, see figure \ref{doubleannu}. Moreover, if $V$ has $n$-wings it suffices to kill all wings but one. Therefore we can always assume that every curve coming from a solid torus to be non-peripheral in $\partial M$.
 	
	\begin{center}\begin{figure}[h!]
		\centering	
		\def\svgwidth{250pt}
		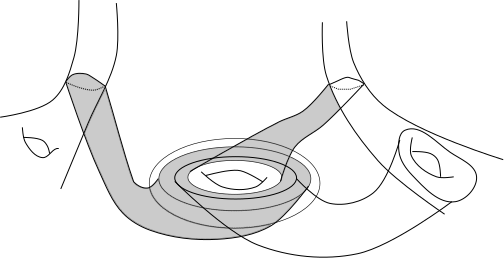
		\caption{The doubly peripheral cylinder is shaded.}\label{doubleannu}
	\end{figure}\end{center}

	 \end{enumerate}

 The construction of the system of curves $P$ will be highly non-canonical since it involves, among other choices, the choice of filling simple closed curves on essential subsurfaces of $\partial\overline M$.

\bprop
Let $(\overline M,\partial\overline M)\in\cat{Bord}(M)$ be maximal bordification for $M\in\M^B$. If $\overline M$ has no doubly peripheral cylinder then we can find a collection $P\subset \partial \overline M$ of pairwise disjoint simple closed curves such that $\overline M$ is acylindrical relative $P$.
\eprop
\bpf

By Theorem \ref{exisjsj} let $(N,R)\subset (\overline M,\partial\overline M)$ be a characteristic submanifold. Consider all components of $N$ that are homeomorphic to $F\times I$ for $F$ a hyperbolic surface, i.e. $\chi(F)<0$, and denote this collection by $N_1$. Denote by $N_2$ the collection of solid and thickened tori in $N$ so that $N=N_1\coprod N_2$. We define $R_1\eqdef R\cap N_1$ and we let  $R_2$ be the intersection of $N_2$ with the non-tori components of $R$ so that all components of $R_2$ are annuli in $\partial\overline M$.

For every component  $F\times I\overset\phi\cong Q\subset N_1$ that is not a pair of pants we have two filling essential simple closed curves $\alpha_0,\alpha_1\hookrightarrow F$, see \cite[3.5]{FM2011}. We can embed $\alpha_i$ in $\phi(F\times\set i)\subset\partial\overline M$, $i=0,1$ and we denote the collection of these simple closed curves in $\partial\overline M$ by $P_1$. Let $\mathcal S$ be the collection of non-annular components of $\partial\overline M$ containing components of $R_1\cup R_2$. Thus, no $S\in\pi_0(\mathcal S)$ is a torus or an annulus. Since $\partial\overline M$ has at most countably many boundary components we can label $\mathcal S$ by the natural numbers so that $\mathcal S=\set{S^n}_{n\in\N}$ and for a component $S^n\subset\mathcal S$ define $\Sigma^n\eqdef\overline{ S^n\setminus N_\epsilon(P_1)}$.

Consider the collection of simple closed curves $\Gamma^n\subset \Sigma^n$ induced by the boundaries of the subsurfaces $R_1\cap S^n$ and the core curves of the annuli in $R_2\cap S^n$. From $\Gamma^n$ throw away all simple closed curves that are peripheral in $\Sigma^n$ and remove redundancies, i.e. if $\beta,\gamma\subset \Gamma^n$ are isotopic we only keep one of them. 

\vspace{0.3cm}

\paragraph{Claim:} Any component $S$ of $\Sigma^n$ containing a component of $\Gamma^n$ has $\chi(S)<0$ and is not a pair of pants.

\vspace{0.3cm}

\bpfc Since no two components of $P_1$ are isotopic and no component of $P_1$ is peripheral in $S^n$ we have that any component $S$ of $\Sigma_n$ has $\chi(S)<0$. Say we have a pair of pants $Q\in\pi_0(\Sigma^n)$ then $\partial Q$ are either peripheral in $\partial\overline M$ or isotopic to elements of $P_1$. Since any component $\gamma$ of $\Gamma^n$ is a simple closed curve we have that $\gamma$ is peripheral in $Q$ hence it is either peripheral in $\partial\overline M$ or isotopic to an element of $P_1$ and neither case can happen.
\epfc

The curves $\Gamma^n\subset \Sigma^n$ are pairwise disjoint, not isotopic and not peripheral. Since $M\in\M^B$ all components of $\partial\overline M$ are of finite type, thus so are the $\Sigma^n$'s. By the previous claim every component of $\Sigma^n$ containing elements of $\Gamma^n$ is not a pair of pants and we denote this components by $X_k^n$, $1\leq k\leq j_n$.

 Since for all $k$ the curve complex $\mathcal C(X_k^n)$ has infinite diameter \cite[2.25]{Schleimer} we can pick an essential simple closed curve $\gamma_k^n\subset X_k^n$ such that for all simple closed curves $\gamma\in\pi_0( \Gamma^n\cap  X_k^n)$ we have that the geometric intersection number $i(\gamma,\gamma_k)>0$. We then add, for all $n,k\in\N$, the curves $\gamma_k^n$ to $P_1$ and denote this new collection by $P$.

By construction we have that every essential cylinder $(C,\partial C)\rar (\overline M,\partial \overline M)$ intersects some component of $P$. Therefore, we have that $\overline M$ is acylindrical with respect to $P$. \epf

								\bdefi
								We say that a manifold $(\overline M,P)$, with $\overline M\in \cat{Bord}(M)$ the maximal bordification for $M\in\M^B$, is an \emph{infinite-type acylindrical pared 3-manifold} if:
								\begin{itemize}
								\item[(i)] the components of $P$ are tori, annuli (closed, open, half open) and all annular and tori components of $\partial\overline M$ are in $P$;
								\item[(ii)] cusp neighbourhoods of $\partial\overline M$ are contained in $P$;
								\item[(iii)] the manifold $\overline M$ is acylindrical rel $P$.
								\end{itemize}
								\edefi

 \bese
 Consider the manifold $M$ of Example \ref{IWSD} and let $\overline M$ be its maximal bordification. Then, $\partial \overline M=A\coprod \amalg_{n=1}^\infty S_n$ where the $S_n$ are genus two surfaces, $A$ is an open annulus and the characteristic submanifold $N$ is given by an infinitely winged solid $T_\infty$ such that $T_\infty\cap \partial\overline M$ are neighbourhoods of the core curve $\beta$ of the annulus $A$ and the separating loops $\alpha_n\subset S_n$ that split the genus two surface into two punctured tori. 
 
 Let $\Gamma$ be a collection of simple closed curves $\gamma_n\subset S_n$ such that $\iota(\gamma_n,\alpha_n)=2$ and $N_\epsilon(\gamma_n,\alpha_n)$ is an essential 4-punctured sphere in $S_n$. Then, by defining $\P\eqdef \set{N_\epsilon(\gamma_n)}_{n\in\N}\amalg A$ we obtain that $(\overline M,P)$ is an infinite type acylindrical pared 3-manifold. In particular, every $M_i$ is acylindrical rel $P_i$ where $P_i$ are the components of $P$ isotopic in $\partial M_i$.
 \eese

								 \subsubsection{Eventual acylindricity of the $M_i$}

In Definition \ref{boundatinf} we gave a decomposition of $\partial M_i$ into two essential sub-surfaces $\partial_\infty M_i$ and $\partial_bM_i$ such that $\partial_\infty M_i$ is the essential subsurface isotopic to infinity in $M\setminus\text{int}(M_i)$. We will now show that $\partial_bM_i$ has ``bounded homotopy class". That is, $\partial_b M_i$ has no essential loops\footnote{Recall that by an essential loop we mean a $\pi_1$-injective loop in a surface $S$ not homotopic into $\partial S$.} homotopic arbitrarily far into $M$. 

\bdefi \label{paredinfinite}
Let $(\overline M,P)$ be an infinite-type pared acylindrical 3-manifold for $M=\cup_{i\in\N} M_i\in\M$. Define, $P_i\subset \partial M_i$ to be the collection of annuli $\hat A\subset \partial M_i$ such that in $\overline M$ we have $I$-bundles $\psi:\mathbb A\times I\hookrightarrow \overline M$ such that $\psi(\mathbb A\times 0)\in\pi_0(\hat A)$, $\psi(\mathbb A\times \set1)$ is a compact annular component of $P$ and for some $\epsilon>0$ $\psi(\mathbb A\times[0,\epsilon))\subset \overline{M\setminus M_i}$. Let $Q_{n_i}$ be the characteristic submanifold of $M_{n_i}$ rel $P_i$ and define $M_{n_i}^{\text{acyl rel }P}$ to be $\overline{M_{n_i}\setminus Q_{n_i}}$.
\edefi
The definition is so that if we consider the cover $K_i$ of $(\overline M,P)$ corresponding to $\pi_1(M_i)$ we have that in the compactification $\overline K_i\cong M_i$ of $K_i$ the lifts of $P$ are isotopic to the $P_i$'s in $\overline K_i\setminus M_i$.

%

\bprop\label{eventualacylindricalcomponents} Given an exhaustion $\set{M_i}_{i\in\N}$ of $M$ by hyperbolizable 3-manifolds with incompressible boundary for all $i$ there exists $n_i$ such that for $n\geq n_i$ no essential loop of $\partial_b M_i$ is homotopic into $\partial M_{n}$ in $M\setminus\text{int}(M_i)$. Moreover, if $(\overline M,P)$ is an infinite-type pared acylindrical 3-manifold there exists $n_i$ such that $M_i\subset M_{n_i}^{\text{acyl rel }P}$.
\eprop
\bpf For $k\in\N$ we define $X_k\eqdef \overline{M_{k+1}\setminus M_k}$. 
%
							\vspace{0.3cm}	
								
						\paragraph{Step 1:}  For all $i$ there exists $n_i$ such that no essential curve in $\partial_ {b} M_i$ is homotopic  in $M\setminus\text{int}(M_i)$ into $\partial M_{m}$ for $m\geq n_i$.
						
						\vspace{0.3cm}

Assume that we have a collection $\set{\gamma_k}_{k\in\N}$ of homotopically distinct essential curves in $\partial_ {b} M_i$ such that for each $n>i$ we can find $k_i\in\N$ such that for $k\geq k_i$ the curve $\gamma_k$ is homotopic into $\partial  M_n$. We want to show that a tail of  $\set{\gamma_k}_{k\in\N}$ always fills a fixed essential subsurface $F$ of $S$. That is, for all $k\in\N$ we have that $\text{Fill}(\cup_{n\geq k}\gamma_n)=F$. Let $n_k$ be such that $\gamma_k$ is homotopic into $\partial M_{n_k}$ via a homotopy $H_k$. Without loss of generality we can assume that $\im(H_k)$ is contained in $\overline{M _{n_k}\setminus M_i}$.

								 Moreover, up to reordering, we can assume that for $i<j$ we have $n_i\leq n_j$ and $n_i\rar\infty$. Let $F_0$ be the essential compact subsurface of $\partial_bM_i$ filled by the $\set{\gamma_k}_{k\in\N}$. Since the $\set{\gamma_k}_{k\in\N}$ fill $F_0$ and are homotopic in $\overline{ M_{n_1}\setminus  M_i}$ into $\partial M_{n_1}$ we have that $F_0$ is homotopic into $\partial M_{n_1}$ in $\overline{ M_{n_1}\setminus  M_i}$. Thus, by JSJ theory we have an essential $I$-bundle $F_0\times I \rar \overline{ M_{n_1}\setminus  M_i} $ and let $F_0'$ be the induced surface in $\partial M_{n_1}$.

								  Let $k_1\in\N$ be such that $\gamma_k$ is homotopic to $\partial M_{n_1+1}$ for $k\geq k_1$. Denote by $\set{\gamma^{(2)}_{k}}_{k\geq k_1}$ the curves in $\partial M_{n_1}$ homotopic to the $\gamma_k$ via the $I$-bundle. Then the curves $\set{\gamma^{(2)}_{k}}_{k\geq k_1}$ fill an essential subsurface $F_1'\subset F_0'$ in $\partial M_{n_1}$ which is isotopic in $\overline{M\setminus M_i}$ to an essential subsurface $F_1\subset F_0$ in $\partial M_i$. Thus we have $ \abs{\chi(F_1)}\leq \abs{\chi(F_0)}$. By iterating this process we obtain a nested sequence of connected essential subsurfaces $\set{F_j}_{j\in\N}$ of $\partial_bM_i$ such that $\forall j\in\N$: $0\leq\abs{\chi(F_{j+1})}\leq\abs{\chi(F_j)}$ and $F_j$ is isotopic in $\overline{M\setminus M_i}$ into $M_{m_j}$ for $m_j\rar\infty$. Since every-time the surface shrinks the absolute value of the Euler characteristic goes down the sequence $\set{F_j}_{j\in\N}$ must stabilise to an essential subsurface $F$. Moreover, since no $\gamma_k$ is peripheral $F$ is not a peripheral annulus in $\partial_bM_i$. Therefore, the tail of the $\set{\gamma_k}_{k\in\N}$ always fills a fixed subsurface $F$ of $\partial_bM_i$. 
								  
								  By making the homotopies $H_k$ immersions and transverse to $\partial M_n$'s we can homotope them to be essential in each $\overline{M_{n+1}\setminus M_n}$ for $i\leq n\leq n_k-1$. Then, by picking a collection of the $\gamma_k$'s that fill $F$ we get that $F$ is homotopic in $M\setminus\text{int}(M_i)$ into $\partial M_{i+n}$ for all $n\in\N$.
								 
Then by JSJ theory we have a collection of essential $I$-bundles: 
$$\phi_n:F\times [0,a_n]\hookrightarrow M_{i+n}\setminus\text{int}(M_i)$$
 $a_n\in\N$ and $a_n\geq n$, such that $\phi_n(F\times\set0)=\text{Fill}(\gamma)\subset\partial_b M_i$ and $\phi_n(F\times\set {a_n})\subset M_{i+n}$. Moreover, up to an isotopy of the $\phi_n$'s we can assume that for all $n\in\N$:
$$\phi^{-1}_n(\im (\phi_n)\cap\cup_{m=i}^{i+n}\partial M_m)=\cup_{0\leq a\leq a_n}F\times\set{a}$$ 
and such that for all $0\leq a\leq a_n$ $\phi_n(F\times [a,a+1])$ is an essential $I$-bundle in some $X_{k}$, $i+1\leq k\leq i+n$, see the proof of Theorem \ref{prodstandardform}.

					By applying Lemma \ref{isoIbundle} to $\set{\phi_n\vert_{F\times[0,1]}}_{n\in\N}$, i.e. to the restriction of the $\phi_n$ so that $\phi_n(F\times[0,1])\subset X_{i+1}$, we obtain, up to isotopy, a sub-sequence such that for all $n\in\N: \phi_n(F\times[0,1])=P_1\subset X_{i+1}$. We then repeat the argument to $\set{\phi_n\vert_{F\times[1,2]}}$ in $X_{i+2}$ to obtain a subsequence such that $\phi_n(F\times[1,2])=P_2$. Moreover, we have that $P_1\cup P_2$ is naturally homeomorphic to an $I$-bundle over $F$. Then, one works in either $X_{i+1}$ or $X_{i+2}$ depending on wether the lid $\phi_n(F\times \set 2)$ of $P_2$ is contained in $\partial M_{i+2}$ or $\partial M_{i+1}$. Note that the new $I$-bundle $P_3$ obtained in $X_{i+1}$ is disjoint from $P_1$ since otherwise the $\phi_n$ were not isotopic to embeddings.

						By iterating this argument we obtain a collection of pairwise disjoint $I$-bundles $P_n$ such that $P_n$ and $P_{n+1}$ have a matching lid. Then $P\eqdef \cup_{n\in\N} P_n\subset M\setminus \text{int}(M_i)$ gives a product $\P:F\times[0,\infty)\hookrightarrow M$, contradicting the fact that $F\subset\partial_b M_i$ was not peripheral.			
								 
\vspace{0.3cm}

We only need to prove the last claim. By \textbf{Step 1} we get that the only cylinders from $\partial M_i$ to $\partial M_n$ with $n\geq n_i$ have boundaries that are peripheral in $\partial_b M_i$.

If $M_i$ is not in $M_{n}^{\text{acyl rel }P}$, $n\geq n_i$ by Lemma \ref{char-submanifold} it means that we can find an annulus of the following type:
\begin{itemize}
\item an embedded essential cylinder $C\subset M_{n}$ such that $C_i^n\eqdef C\cap M_i$ is essential in $M_i$ with boundary homotopic to $\partial_bM_i$;
\item an immersed annulus $C$ formed by an embedded annulus $C_1\subset M_n\setminus\text{int}(M_i)$ with one boundary $\gamma$ homotopic to $\partial_b M_i$ and the non-trivial homotopy contained in a solid torus $C^i_n(r)$ obtained by collapsing $\gamma\simeq r^n$ to the root $r$.
\end{itemize}

\vspace{0.3cm}

\paragraph{Step 2:} If $(\overline M,P)$ is an infinite-type pared acylindrical 3-manifold there exists $n_i$ such that $M_i\subset M_{n_i}^{\text{acyl rel }P}$.

\vspace{0.3cm}

Assume that the statement is not true and assume that no peripheral element of $\partial_b M_i$ has a root in $M_i$ or $X_j\eqdef \overline{M_j\setminus M_i}$, $j>i$. Then, by Lemma \ref{char-submanifold} we have essential cylinders: 
$$C_n:(\mathbb S^1\times I,\mathbb S^1\times\partial I)\hookrightarrow(M_n, \partial M_n\setminus P_n), \qquad n\geq n_i$$
 such that $\im( C_n)\cap M_i$ is a collection of pairwise disjoint essential cylinders whose boundaries are peripheral in $\partial_b M_i$. By doing isotopies, supported in $M_i$, of the $\set{C_n}_{n\in\N}$ and picking a sub-sequence $\set{C^{(1)}_n}_{n\in\N}$ we can assume that the $\set{C^{(1)}_n}_{n\in\N}$ satisfy $\im(C^{(1)}_{n})\cap M_i=\im(C^{(1)}_m)\cap M_i$ for all $n,m\geq n_i$.

By repeating this argument in $M_k$, $k\geq i$, and picking a diagonal subsequence $\set{C_n}_{n\in\N}$ we obtain a bi-infinite cylinder $\hat C=\varinjlim_n C_n$ with $\hat C:\mathbb S^1\times\R\hookrightarrow M$ such that for all $k\in\N$ $\hat C_k\eqdef \hat C\cap M_k$ are cylinders and $\hat C_i$ contains an essential cylinder. Let $\hat M$ be a bordification where $\hat C$ compactifies and denote by $\alpha_1,\alpha_2$ the boundaries of $\hat C$. Since $\hat C_i$ is essential we have that $\hat C$ is also essential in the bordification $\hat M$. By uniqueness of the maximal bordification we can assume that $\overline M=\hat M$. Since $(\overline M,P)$ is an infinite-type pared acylindrical 3-manifold and $\hat M\cong \overline M$ we have at least one component $\gamma$ of $P$ such that $\iota(\gamma,\partial\hat C)>0$, say that $\iota(\gamma,\alpha_1)>0$. Moreover, since $\gamma$ is not peripheral in $\partial\overline M$ so is $\alpha_1$ and let $S\in\pi_0(\partial\overline M)$ be the component containing them.

Pick $k>n_i$ such that we have $\gamma'\iso\gamma$ contained in $P_k\subset \partial M_k$\footnote{See Definition \label{paredinfinite}.} and such that the unbounded component $\hat C^1$ of $\hat C\cap\overline{ M\setminus M_k}$ compactifying to $\alpha_1$ is contained in a product in standard form: 
$$\mathcal Q:(F\times[0,\infty),F\times\set 0)\hookrightarrow( \overline {M\setminus M_k},\partial_\infty M_k)$$
where $F$ is a surface isotopic into $S\subset \partial\overline M$ in $ \overline {\overline M\setminus M_k}$ containing $\gamma'$. Moreover, up to a proper isotopy of $\hat C$, we can assume that $\hat C^1=\mathcal Q(\alpha_1\times[0,\infty))$.

Since $\varinjlim C_n=\hat C$ we can pick a cylinder $C_n$ such that $\im(\hat C)\cap M_k\subset \im(C_n)\cap M_k$ and let $C^1_n$ be the component of $\im(C_n)\cap \overline{ M\setminus M_k}$ containing $\mathcal Q(\alpha_1\times\set0)$. Because $\im(\mathcal Q)\subset \overline{ M\setminus M_k}$ and $\mathcal Q$ is in standard form we have a minimal $t=t_n>0$ such that $\mathcal Q(F\times\set t)\subset \partial M_n$, $n\geq k$, hence:
$$\mathcal Q:(F\times [0,t], F\times\set{0,t})\hookrightarrow (\overline{M_n\setminus M_k},\partial(\overline{M_n\setminus M_k}))$$
\noindent is an essential $I$-bundle. 

Since $\alpha_1$ is not peripheral in $F$ and $C_n^1\cap\partial M_k=\mathcal Q(\alpha_1\times\set 0)$ we have that the annulus $C_n^1$, is up to isotopy, vertical in $\mathcal Q (F\times[0,t])$. The simple closed curve $\mathcal Q(\gamma'\times \set{t})$ is isotopic through $\mathcal Q$ to $\gamma\subset P$ and we have some $\epsilon>0$ such that $\mathcal Q(\gamma'\times[t,t+\epsilon))$ is contained in $M\setminus \text{int}(M_n)$. Hence, we have that $\mathcal Q(\gamma'\times \set t)\subset \P_n$. Therefore, $\partial C_n$ is not in $\partial M_n\setminus P_n$ reaching a contradiction.

Say that peripheral elements of $\partial_b M_i$ have roots in $M_i$ and that we do not have any collection of embedded cylinders so that we have:
$$C_n:(\mathbb S^1\times I,\mathbb S^1\times\partial I)\rar(M_n, \partial M_n\setminus P_n), \qquad n\geq n_i$$
whose images in $M_i$ contain a root of $\partial _b M_i$. By Lemma \ref{char-submanifold}, up to a homotopy, we have that all $C_n$ contain a sub-cylinder $C_n^i(r)\subset M_i$ which is the non-trivial homotopy obtained by collapsing a wing of a solid torus to a power of the core and re-expanding it and are embedded otherwise. Since $\partial_b M_i$ has finitely many peripheral elements and they have unique roots up to a subsequence we can assume that all $C_n^i(r)$ are the same and then the previous argument applies to the embedded part of the $C_n$'s. The argument for the case in which the root is contained in $X_j$ is similar. \epf

By combining Proposition \ref{eventualacylindricalcomponents} and Theorem \ref{prodstandardform} we obtain the following corollary:

\bcor\label{minimalexhaustion}
Let $M=\cup_{i=1}^\infty M_i\in\M^B$ then for all $i$ $ M_i$ is contained in an open 3-manifold $\hat M_i\cong M_i\cup\partial_\infty M_i\times [0,\infty)$ such that $\partial \hat M_i$ is isotopic to $ \partial_b M_i$ and $\partial_\infty M_i\times [0,\infty)$ is in standard form. Moreover, for all $n>i$ we have that $ \partial M_n\cap M\setminus\text{int}(\hat M_i)$ are properly embedded surfaces and there is $j=j(i)\in\N$ such that for all $n\geq j$ no essential loop of $ \partial\hat M_i$ is homotopic into $\partial M_n$.
\ecor 
\bpf By Lemma \ref{eventualacylindricalcomponents} we have a maximal essential subsurface $S\cong\partial_\infty M_i\subset \partial M_i$ that generates a properly embedded product 
$$\P:(S\times [0,\infty),S\times\set 0)\hookrightarrow (M,\partial _bM_i)$$
in $M\setminus \text{int}(M_i)$ and we have $j=j(i)\in \N$ such that for all $n\geq j$ no essential loop of $ \partial_b M_i\iso\partial\hat M_i$ is homotopic into $\partial  M_n$. By Theorem \ref{prodstandardform} we can assume $\P$ to be in standard form. Thus, $\hat M_i\eqdef M_i\cup_{\partial_\infty M_i}\P\cong M_i\cup\partial_\infty M_i\times [0,\infty)$ is the required submanifold. \epf

							\subsubsection{Homotopy equivalences of 3-manifolds}	 
								 We now proceed with two topological Lemmata that we need in the final proof. The first is a generalisation of Lemma 2.2 in \cite{SS2013} and the latter is a relative version of Johansson homeomorphism Theorem \cite{Jo1979}.
								 
								 \blem\label{2.2}
								 Let $M$ be a complete, open hyperbolic 3-manifold and $K$ a compact, atoroidal, aspherical, 3-manifold with incompressible boundary such that $\iota: K\rar M$ is a homotopy equivalence and $\iota\vert_{\partial K}:\partial K\hookrightarrow M$ is an embedding. Then $\iota$ is homotopic relative to $\partial K$ to an embedding $\iota':K\rar M$.
								 \elem
								 \bpf
								 Since $M$ is homotopy equivalent to $K$ we have that $\pi_1(M)$ is finitely generated. By the Tameness Theorem \cite{AG2004,CG2006} it follows that $M\cong \text{int}(\overline M)$ for some compact manifold $\overline M$. Therefore, we have $\iota: K\rar \overline M$ with the same properties as in the statement.
								 
								 We need to show that $\iota(\partial K)$ is peripheral in $\overline M$ since then by Waldhausen's Theorem \cite{Wa1968} follows that $\iota$ is homotopic to an embedding. Since $\iota$ is a homotopy equivalence the map on homology $\iota_*$ induces an isomorphism: $\iota_*:H_*(K)\rar H_*(\overline M)$. Let $\partial K=\coprod_{i=1}^n S_i$ then we can write: $H_3(K,\partial K)\cong\Z\langle [K,\partial K]\rangle$ with $\partial[K,\partial K]=\sum_{i=1}^n[S_i]$ and $H_2(\partial K)\cong \Z\langle [S_i]\rangle$. By the long exact sequence of the pair $(K,\partial K)$: 
								 \be 0\rar H_3(K,\partial K)\rar H_2(\partial K)\rar H_2(K)\ee
								 we obtain the following injection:
								 \be \quotient{\Z\langle [S_i]\rangle}{\sum_{i=1}^n[S_i]}\hookrightarrow H_2(K)\overset{\iota_*}\cong H_2(\overline M)\cong H_2(M)\ee
								  This means that no linear combination of the $[S_i]$ except $\sum_{i=1}^n [S_i]$ is null-homologous in $\overline M$. Moreover, the $\set{S_i}_{i=1}^n$ are separating in $K$ since they are not dual to any 1-cycle in $H_1(K)$. Since $\iota$ preserves homological conditions the same holds for the $\set{\iota(S_i)}_{i=1}^n$. Therefore, all the $\iota(S_i)_{i=1}^n$ are separating in $M$ and no linear combination except $\sum_{i=1}^n[\iota(S_i)]$ is null-homologous in $M$. Hence, if we start splitting along the $\set{[\iota(S_i)]}_{i=1}^n$ we get a connected submanifold $N\subset M$ whose boundary is $\sum_{i=1}^n [\iota (S_i)]$. If we show that $\overline {M\setminus N}$ is a product manifold over $\partial N\cong\partial K$ we are done.
								 
								 Consider a homotopy inverse $f:\overline M\rar K$ then $f\vert_N:N\rar K$ is a $\pi_1$-injective map and up to homotopy it sends $\partial N\rar \partial K$ homeomorphically. The map $f\vert_N$ is also degree one since $f\vert_N=f\circ\iota_N$ and $f$ has degree one. Thus $f\vert_N:N\rar K$ is a $\pi_1$-injective degree one map. We now claim that $f\vert_N$ is a surjection on $\pi_1$ and so a homotopy equivalence. If $f_*$ is not surjective let $H\eqdef f_*(\pi_1(N))\subgroup \pi_1(K)$ and consider the cover $\pi:K_H\twoheadrightarrow K$ corresponding to $H$. Since by construction the map $f\vert_N$ lifts to $\tilde f\vert_N$ we have: $f\vert_N=\pi\circ\tilde f\vert_N$, but $f\vert_N$ has degree one hence $\deg\pi=\pm 1$. Therefore, $K_H\cong K$ which implies that $(f\vert_N)_*$ is a surjection on $\pi_1$, hence an isomorphism.
								 
								 By Whitehead's Theorem we get that $f\vert_N$ is a homotopy equivalence which implies that $\iota: N\hookrightarrow M$ is also a homotopy equivalence and so $N$ is a Scott core. Moreover, since $\iota(\partial K)=\partial N$ is incompressible by Corollary \cite[5.5]{Wa1968} we have that each component of $\overline {M\setminus N}$ is a product over $\iota(S_i)$. Hence the homotopy equivalence $\iota: K\rar \overline M$ is homotopic rel boundary to a homotopy equivalence $\iota':K\rar N\subset M$ that is a homeomorphism on $\partial K\rar \partial N$. Thus by Waldhausen's Theorem \cite{Wa1968} we get that it is homotopic rel boundary to a homeomorphism from $K$ to $N$. Hence, we get an embedding $\iota':K\hookrightarrow M$ with $\iota'\vert_{\partial K}=\iota\vert_{\partial K}$. 								 \epf

								 The following is an application of the classification theorem of Johansson \cite[2.11.1]{CM2006} for homotopy equivalences between 3-manifolds.

								 \blem\label{reljoh}
								 Let $\phi: M \rar N$ be a homotopy equivalence between compact, irreducible, orientable 3-manifolds and let $X\subset M$ be a codimension-zero submanifold. If  $S$ is a collection of essential subsurfaces of $\partial M$ such that $X$ is contained in the acylindrical part of $M$ relative to $S$ and $\phi\vert_{S}:S\rar\partial N$ is an embedding, then we can homotope $\phi$ to $\psi$ so that $\psi\vert_X:X\rar N$ is an embedding and the homotopy is constant on $S$.
								 \elem
								 \bpf
								 Complete $S$ and $\phi(S)$ to useful boundary patterns for $\partial M,\partial N$, which we denote by $\overline S$ and $\overline{ \phi(S)}$ respectively. Let $V,Z$ be the characteristic submanifolds corresponding to $\overline S$ and $\overline{ \phi(S)}$  then we have that $X\subset \overline{M\setminus V}$. By the Johansson Classification theorem \cite{Jo1979} we have that $\phi$ is admissibly homotopic to a homeomorphism $\psi: (\overline{M\setminus V}, S)\rar (\overline { N\setminus Z},\phi(S))$. An admissible homotopy is a homotopy by pair maps hence since $\phi\vert_S$ is already a homeomorphism we can choose it to be constant on $S$. Since we assumed that $X\subset\overline{ M\setminus V}$ we get that $X$ is embedded by $\psi$ and $\psi\vert_S=\phi\vert_S$. 								 \epf

\section{Proof of the Main Theorem}
In this section we prove our main Theorem:
\begin{customthm}{1} \label{maintheorem} Let $M\in\M^B$. Then, $M$ is homeomorphic to a complete hyperbolic 3-manifold if and only if  the associated maximal bordified manifold $\overline M$ does not admit any doubly peripheral annulus. 	\end{customthm}

In the next subsection we show that not having doubly peripheral annuli is a necessary condition. Specifically, we prove that if $M\in\M^B$ does not have a doubly peripheral cylinders then it is homotopy equivalent to a hyperbolic 3-manifolds $N$. Finally we show that particula homotopy equivalences between $M$ and $N$ are homotopic to homeomorphisms.

\subsection{Necessary condition in the main Theorem}\label{necessarycond}

Using the same techniques of \cite{C20171} we prove the necessary condition on the annulus in Theorem \ref{maintheorem}. We start with a remark on characteristic submanifolds of manifolds in $\M^B$.

 \brem\label{doublyperipheralannuli}
By Theorem \ref{exisjsj} for a manifold $M\in\M^B$ we have a characteristic submanifold $N$ for the maximal bordification $\overline M$. Then, any doubly peripheral cylinder $C$ is homotopic into one of the following components of $N$:
\begin{itemize}
\item[(i)] a solid torus with at least one peripheral wing in $\partial\overline M$ that wraps around the soul $n>2$ times; 
\item[(ii)] a solid torus with at least two peripheral wings in $\partial\overline M$ each of which wraps around the soul once;
\item[(iii)] a thickened essential torus with at least one wing that is peripheral in $\partial\overline M$;
\item[(iv)] an $I$-bundle $P\cong F\times I$ such that at least one component of $\partial F\times I$ is doubly peripheral.
\end{itemize} 

For (i),(ii) and (iii) the cases with infinitely many wings are also allowed. However, for all cases except (i) if $C$ is a doubly peripheral annulus then there exists a properly embedded annulus $C'$ that is also doubly peripheral and such that $C$ is homotopic into $C'$, i.e. $C$ is a \emph{power} of $C'$.
\erem

We first show that if $M$ is hyperbolizable, i.e. $M\cong\hyp\Gamma$, the elements of $\pi_1(M)$ that are peripheral in $\partial \overline M$ are represented by parabolic elements in the Kleinian group $\Gamma$.

\blem\label{parabolicelement}
For $M\in\mathcal M^B$ let $\overline M\in \cat{Bor}(M)$ be the maximal bordification. If $M\cong\hyp\Gamma$ admits a complete hyperbolic metric and $\gamma\in \pi_1(M)$ is homotopic to $\overline \gamma\subset \partial \overline M$ such that $\overline\gamma$ is peripheral in $\partial \overline M$ then $\gamma$ is represented by a parabolic element in $\Gamma$.
\elem
		
\bpf  Let $\set{ M_i}_{i\in\N}$ be the exhaustion of $M$ and let $G$ be the bound on the Euler characteristic of the boundary components of the $ M_i$.  Without loss of generality it suffices to consider the case where $\overline\gamma$ is a simple closed curve. If $\gamma$ is peripheral in $\partial\overline M$ the components of $\cup_{i\in\N} \partial M_i $ that have a simple closed curve $\gamma_n$ isotopic to $\gamma$ in $\overline M\setminus\text{int}(M_i)$ form a properly embedded sequence of hyperbolic incompressible surfaces $\set{\Sigma_n}_{n\in\N}$ with $\Sigma_n\in\pi_0(\partial  M_{i_n})$. Moreover, up to picking a subsequence we can assume that $i_n< i_{n+1}$.

If $\gamma$ is not represented by a parabolic element it has a geodesic representative $\hat\gamma$ in $M$ which is contained in some $ M_i$. Let $\tau_n$ be a 1-vertex triangulation of $\Sigma_n$ realising $\gamma_n$. Since the $\Sigma_n$ are incompressible closed surfaces we can realise them, in their homotopy class, via simplicial hyperbolic surfaces $(S_n,f_n)$ in which $\gamma_n$ is mapped to $\hat\gamma$ (see \cite{Ca1996,Bo1986}). Since the $\Sigma_n$ are hyperbolic surfaces and $\gamma_n$ is simple each $S_n$ contains at least one pair of pants with one boundary component $\gamma_n$.

By Gauss-Bonnet for all $n$ we have: 
$$A(S_n)\leq 2\pi\abs{\chi(S_n)}\leq 2\pi G$$
Since the $S_n$ have uniformly bounded area,  we have that any maximally embedded \emph{one-sided} collar neighbourhood of $\gamma_n$ in $S_n$ has radius uniformly bounded by $K\eqdef K(\ell_M(\hat\gamma),G)$. 
Therefore, for any $\epsilon>0$ an essential pair of pants $P_n\subset S_n$ with $\gamma_n$ in its boundary is contained into a $K +\epsilon$ neighbourhood of $\gamma_n\subset S_n$. Since the maps $f_n$ are 1-Lipschitz for all $n$ the $f_n(P_n)$ are contained in a $K+\epsilon$ neighbourhood of $\hat\gamma=f_n(\gamma_n)$ in $M$.

Since for any pair of pants $P_n\subset S_n$ with $\gamma_n\subset\partial P_n$ the $f_n(P_n)$'s are at a uniformly bounded distance from $\hat\gamma$ we can assume that all such pair of pants are contained in $\hat M_k$ for some $k$, where $\hat M_k$ is as in Corollary \ref{minimalexhaustion} and we also let $j\eqdef j(k)$ be as in the Corollary. Moreover, we can assume that there is a cusp component $Q$ of $\partial\hat M_k$ that compactifies to a simple closed loop $\alpha\subset\partial\overline M$ isotopic to $\overline\gamma$ in $\partial\overline M$. We now want to show that we can find $n\in\N$ such that $\Sigma_n$ has a pair of pants $P_n$ with boundary $\gamma_n$ contained outside $\hat M_k$.

\vspace{0.3cm}

\paragraph{Claim:} There is a component of $\cup_{n\geq j} \partial M_n\setminus\text{int}( \hat M_k)$ that has $\gamma_n\subset Q$ as a boundary component and is not an annulus.

\vspace{0.3cm}

\bpfc Since $\partial_\infty M_k\times [0,\infty)$ is in standard form we have that all components $S$ of $\cup_{i\in\N} M_i$ intersect $\partial_\infty M_k\times [0,\infty)$ in level surfaces or are disjoint from it. If the claim is not true we have $j\in\N$ such that every component of $\cup_{i\geq j} M_i\setminus\text{int}( \hat M_k)$ having $\gamma_n$ as a boundary component is an annulus. Each such component $A$ has boundary on $\Lambda\eqdef \partial (\partial_\infty M_k)\times [0,\infty)$.

Since $\Lambda$ has finitely many components we can assume that we have a collection of annuli $\set{A_\ell}_{\ell\in\N}$ with $\partial A_\ell=\alpha\times\set{t_\ell}\cup\beta\times \set{t_\ell}\subset \Lambda$. Then, the $\set{A_\ell}_{\ell\in\N}$ have to be in at most two homotopy classes. If not we have two essential tori that have homotopic boundaries and this cannot happen in hyperbolic 3-manifolds, see Remark \ref{annuliintersection}. Thus, we get that eventually we can enlarge the product $\partial_\infty M_i\times [0,\infty)$ so that $Q$ is not peripheral anymore. This gives a new bordification $M'$ such that $\overline M\subsetneq M'$ contradicting the maximality of $\overline M$ and the fact that $\bar\gamma$ was peripheral. \epfc

Thus, we have a surface $F_n\subset \pi_0(\partial M_n)$ with $n\geq j$ such that $F_n\cap\hat M_k^C$ contains a pair of pants $P_n$ with a boundary component homotopic in $M$ to $\hat \gamma$. Therefore, the corresponding simplicial hyperbolic surface $f_n:S_n\rar M$ has a pair of pants $P_n$ with $\gamma_n\subset P_n$ such that $f_n(P_n)\subset\hat M_k$. However, this means that $P_n$ is homotopic into $\partial\hat M_k\iso \partial _b M_k$ and since $\pi_1(P_n)$ cannot inject into $\Z$ we reach a contradiction with Corollary \ref{minimalexhaustion}. \epf

We will now show that manifolds in $\M^B$ with a doubly peripheral annulus are not hyperbolic. First we need the following topological Lemma saying that if $\alpha,\beta\subset\partial\overline M$ are peripheral simple closed curves and isotopic in $\overline M$ then we can separate their homotopy class by a compact subset.

\blem\label{separation} Let $\overline M\in\cat{Bord}(M)$, $M\in\M^B$ be the maximal bordification and $\mathcal A:(\mathbb A,\partial\mathbb A)\hookrightarrow (\overline M,\partial\overline M)$ be an essential doubly peripheral annulus. Then, there exists $M_i$ such that in $\overline M\setminus\text{int}( M_i)$ the peripheral loops $\mathcal A(\partial\mathbb A)\eqdef \gamma_0\coprod\gamma_1$ have no essential homotopy and are isotopic in $\overline M\setminus\text{int}( M_i)$ to peripheral loops in $\partial_\infty M_i$. 
\elem
\bpf Since $A$ is embedded it is isotopic in a component $P$ of the characteristic submanifold $N$ of $\overline M$. Thus, by Remark \ref{doublyperipheralannuli} we have three possibilities for $P$:
 \begin{itemize}
 \item[(i)] a solid torus with at least two peripheral wings in $\partial\overline M$; 
\item[(ii)] an $I$-bundle $P\cong F\times I$ such that at least one component of $\partial F\times I$ is doubly peripheral;
\item[(iii)] an essential torus with at least one wing that is peripheral in $\partial\overline M$.
\end{itemize} 

\paragraph{Case (iii).} Let $P$ be the component of the characteristic submanifold $N$ corresponding to the essential torus $T\subset \partial\overline M$ that $A\eqdef \mathcal A(\mathbb A)$ is homotopic into. Then, by Lemma \ref{essentialannuli} there exists a minimal $i$ such that the essential torus $T$ is isotopic in $\overline M\setminus \text{int}( M_i)$ into an essential torus $\overline T$ of $\partial_\infty M_i$ and such that the compact components of $A\cap M_i$ and $A\cap M\setminus\text{int}(M_i)$ are essential and $\gamma_0,\gamma_1$ are isotopic in $\overline M\setminus\text{int}( M_i)$ into peripheral loops $\gamma_0^i,\gamma_1^i$ of $\partial_\infty M_i$. For now assume that $\gamma_0$ is not isotopic to $\gamma_1$ in $\partial \overline M$.

For $j\geq i$ let $X_j^\infty\eqdef \overline M\setminus \text{int} (M_j)$ and assume we have an essential homotopy $C_j:(\mathbb A, \partial\mathbb A)\hookrightarrow (X_j^\infty,\partial X^\infty_j)$ from $\gamma_0$ to $\gamma_1$. Then, $C_j\cup_\partial A$ forms a torus $\hat T$ which is essential since otherwise $A\cap M_i$ was inessential. Moreover, up to a homotopy of $\hat T$ pushing it off $\partial\overline M$ we can assume that $\hat T$ and $\overline T$ are contained in $M_k$ for some $k> i$. Thus, since $\hat T$ and $\overline T$ have homotopic simple closed curves by hyperbolicity of $M_k$ we must have that $\hat T$ is homotopic into $\overline T$ in $M_k$. Hence, we have that $A\cap X_i^\infty$ does not contain any compact annuli since they would be homotopic into $M_i$ contradicting Lemma \ref{essentialannuli}.

Let $C'_j\subset X_i^\infty$ be the subannulus of $\hat T$ obtained by taking $C_j$ and going to $\gamma_0^i,\gamma_1^i\subset\partial M_i$ along $A$. Since $\hat T$ is homotopic into $M_i$ we get that $C'_j$ is inessential in $X_i^\infty$, thus $\partial C'=\gamma_0^i\coprod \gamma_1^i$ are parallel in $\partial M_i$ and co-bound an annulus $A'\subset \partial M_i$. Moreover, we have that $C_j$ is parallel to $A'$ in $X_i^\infty$. Thus we proved:

\vspace{0.3cm}

\paragraph{Claim:} If $C_j:(\mathbb A, \partial\mathbb A)\hookrightarrow (X_j^\infty,\partial X^\infty_j)$ is an essential annulus connecting $\gamma_0,\gamma_1$ then it is isotopic to the annulus $A'\subset\partial M_i$ connecting $\gamma_0^i,\gamma_1^i$. Moreover, for $j\neq\ell$ we also have that $C_j\iso C_\ell$.

\vspace{0.3cm}

If we cannot separate the homotopy class of $\gamma_0,\gamma_1$ we have a collection of annuli $C_j:(\mathbb A, \partial\mathbb A)\hookrightarrow (X_j^\infty,\partial X^\infty_j)$, $j\geq i$, such that for all $k> i$ there exists $j$ such that $C_j\cap M_k=\emp$. Then, we get that in $\partial\overline M$ there are cusps neighbourhoods $P_1,P_2$ of the components $S_1,S_2$ of $\partial\overline M$ co-bounding  with $A$ a submanifold of the form $\mathbb A\times[0,\infty)$ contradicting the properties of a maximal bordification. 

Now assume that $\gamma_0,\gamma_1$ are isotopic in $\partial\overline M$ so that they co-bound an annulus $C\subset \partial\overline M$ and assume that we have an essential annulus $C':(\mathbb A, \partial\mathbb A)\hookrightarrow (X_i^\infty,\partial X^\infty_i)$ with boundary $\partial C$. Since $C'$ is essential we have that $\hat T\eqdef C\cup_\partial C'(\mathbb A)$ is an essential torus which is then homotopic to the torus $T\subset P$ in $\overline M$. Moreover, since $T,\hat T$ are embedded, incompressible and homotopic we can assume by \cite[5.1]{Wa1968} that they are isotopic in $M$ so they co-bound an $I$-bundle $J$. Since $\hat T\subset X_i^\infty$ up to an isotopy of $J$ we can assume that $J\cap\partial M_k$ are level surfaces hence all components of $J\cap\partial M_k$ are essential tori.

Then, either $\hat T$ is contained in the $I$-bundle $Q\cong \mathbb T^2\times I$ generated by the boundary torus $\overline T$ of $\partial M_i$ or it is contained in some other component of $X_i^\infty$. If it is contained in $Q$ we get a contradiction since then $\gamma_0,\gamma_1$ are contained in a torus component of $\partial \overline M$. Thus, since $\hat T\subset X_i^\infty\setminus Q$ and $J\cap \partial M_i\setminus Q$ are essential tori we get that $M_i \cong \mathbb T^2\times I$. In turn, this gives us that $M\cong\mathbb T^2\times \R$ and $\overline M\cong \mathbb T^2
\times I$ which does not contain any doubly peripheral annulus.

\vspace{0.3cm}

We will now deal with annuli of type (i) and (ii) and we can assume that we have no doubly peripheral annulus of type (iii).

\vspace{0.3cm}

\paragraph{Case (i) and (ii).}  Let $A$ be as before an essential annulus connecting $\gamma_0$ to $\gamma_1$ in $\overline M$. By Lemma \ref{essentialannuli} we have an isotopy of $A$ and a minimal $M_i$ such that compact components of $A\cap M_i$ and $A\cap X_i^\infty$ are essential annuli. 

Assume we have an essential annulus $C$ connecting $\gamma_0$ to $\gamma_1$ in $X_i^\infty$. The annuli $C$ and $A$ cannot be parallel since otherwise $A\cap M_i$ would have no essential components. Therefore, by taking a push-off $C'$ of $C$ and connecting it to $A$ along $\gamma_1$ we obtain an essential annulus $A'$ that has both boundaries isotopic to $\gamma_0$ in $\partial\overline M$. Therefore, we contradict the fact that we had no type (iii) annuli and so $M_i$ disconnects the homotopy class of $\gamma_0$ and $\gamma_1$ in $M$. \epf

\bthm \label{nothyp}
If $M\in\mathcal M^B$ is hyperbolic, $M\cong\hyp\Gamma$, then $\overline M$ cannot have an essential doubly peripheral annulus.
\ethm
\bpf Since $M\in \M^B$ we have $G\in\N$ such that for all $i\in\N$ all components $\Sigma$ of $\partial M_n$ have $\abs{\chi(\Sigma)}\leq G$. Let $\mathcal A:(\mathbb A,\partial\mathbb  A)\rar (\overline M,\partial\overline M)$ be an essential annulus such that $\mathcal A(\partial \mathbb A)\eqdef\gamma_1\cup\gamma_2$ are peripheral in $\partial\overline M$. Let $\gamma\in\pi_1(\overline M)$ be the element that generates $\pi_1(\mathbb A)\hookrightarrow\pi_1(\overline M)$. By Lemma \ref{parabolicelement} $\gamma$ has to be represented by a parabolic element and by Remark \ref{doublyperipheralannuli} we only have to consider the following four cases:
 
 \begin{itemize}
\item[(i)] a solid torus with at least one peripheral wing in $\partial\overline M$ that wraps around the soul $n>2$ times; 
\item[(ii)] a solid torus with at least two peripheral wings in $\partial\overline M$ each of which wraps around the soul once;
\item[(iii)] an $I$-bundle $P\cong F\times I$ such that at least one component of $\partial F\times I$ is doubly peripheral;
\item[(iv)] an essential torus with at least one wing that is peripheral in $\partial\overline M$.
\end{itemize} 

Except for (i) we can assume that $\mathcal A$ is an embedding.

 \vspace{0.3cm}

\paragraph{Step 1}$M$ cannot have a doubly peripheral cylinder $\mathcal A:(\mathbb A,\partial \mathbb A)\rar (\overline M,\partial\overline M)$ with $\mathcal A(\partial \mathbb A)\eqdef\gamma_1\cup\gamma_2$ of type (i).

\vspace{0.3cm}

In this case we have a doubly peripheral cylinder $C\eqdef \mathcal A(\mathbb A)$ whose boundaries are isotopic in $\partial\overline M$. Let  $S\subset\partial \overline M$ be the component containing $\partial C$. By construction of the characteristic submanifold and of the maximal bordification we have $M_i$ such that $\partial_{\infty}M_i $ contains a component isotopic to $S$ in $\overline M\setminus \text{int}(M_i)$ and an essential solid torus $V\subset M_i\cup\partial_\infty M_i\times[0,\infty)$ with a wing $w$ whose boundary $\partial C$ is isotopic to a collar neighbourhood of $\gamma_i$ and such that $w$ wraps around the soul $\gamma$ of $V$ $n>1$ times. Also note that in this case $\gamma$ is primitive in $\pi_1(M_i)$.

 Since the cover $\tilde M$ of $M$ corresponding to $\pi_1(M_i)$ is homeomorphic to $\text{int}(M_i)$, see Lemma \ref{lochyp}, we have that in the pared hyperbolic 3-manifold $(N,P)$, $N\cong M_i$, such that $\text{int}(N)=\tilde M$ there are disjoint embedded annuli $A,B\in \pi_0(P)$ such that $\gamma$ is homotopic to the soul $b$ of $B$ and the soul $a$ of $A$ is isotopic to $\gamma_i$. 

Thus, the $n$-th power of the soul $b$ of $B$ is homotopic to the soul $a$ of $A$. Therefore, we have a component of the characteristic submanifold of $N$ that realises this homotopy. However, since $N$ is hyperbolizable this cannot happen because no component $Q$ of the characteristic submanifold has elements in the boundary such that one is a root of the other in $Q$.

\vspace{0.3cm}

Since we dealt with the non-embedded case from now on we can assume that the doubly peripheral cylinder is embedded. The idea is to use a collection of simplicial hyperbolic surfaces $S_n$'s, as in Lemma \ref{parabolicelement}, all intersecting the doubly peripheral annulus in simple loops $\gamma_n$ isotopic to a peripheral loop $\gamma$ in the boundary $\partial\overline M$. By Lemma \ref{separation} the simplicial hyperbolic surfaces $S_n$'s will be forced to go through some $M_i$. By using the hyperbolicity of $M$ this will force loops $\alpha_n$ transverse to the $\gamma_n$ to have uniformly bounded length and this will allow us to construct a product $\P$ whose compactification makes $\gamma$ not peripheral.

\vspace{0.3cm}

\paragraph{Step 2} A hyperbolic $M$ cannot have a doubly peripheral cylinder $(A,\partial A)\rar (\overline M,\partial\overline M)$ with $\partial A\eqdef\gamma_1\cup\gamma_2$ of type (ii)-(iv).
 
\vspace{0.3cm}

Let $M_i$ be as in Lemma \ref{separation} so that $\partial_\infty M_i$ contains peripheral loops $\gamma_1^i,\gamma_2^i$ isotopic to $\gamma_1,\gamma_2$ respectively in $\overline M\setminus 
\text{int}(M_i)$. Moreover, if $\gamma_1\iso\gamma_2$ in $\partial\overline M$ we can assume, by picking a larger $i$, that the essential torus $T$ induced by the doubly peripheral annulus $C$ is contained in $M_i$. Let $Q$ be the cusp neighbourhood, in $M$, of the parabolic element corresponding to the homotopy class of $\gamma$. 

In the case that we have an essential torus $T$ the cusp $Q$ is contained in a component of $\overline{M\setminus M_i}$ homeomorphic to $\mathbb T^2\times [0,\infty)$. Otherwise, it corresponds to a neighbourhood of $\gamma_1,\gamma_2$ in $\partial \overline M$ in which case we will assume it is $\gamma_2$.

Let $\hat M_i$ be the manifold from Lemma \ref{minimalexhaustion} and let $\set{\Sigma_n}_{n> i\in\N}$ be the sequence of surfaces in $\partial M_n\cap \hat M_i^C$ coming from \textbf{Claim 1} of Lemma \ref{parabolicelement} and let $F_n\subset \partial M_n$ the component containing $\Sigma_n$. 
	
Let $\tau_n$ be ideal triangulations of the $\hat\Sigma_n\eqdef\Sigma_n\setminus\gamma_1^n$ where the cusps corresponding to $\gamma_1^n$ have exactly one vertex each and homotope them so to obtain proper maps of the punctured surfaces. Then by \cite{Ca1996,Bo1986} we can realise the embeddings $\hat\Sigma_n\hookrightarrow M$ by simplicial hyperbolic surfaces $(S_n,f_n)$ in which $\gamma_1^n$ is sent to the cusp $Q$. Therefore, since $M_i$ separates in $ M$ the homotopy class of $\gamma_1$ and $\gamma_2$ we have that all the $S_n$'s must intersect $\partial M_i$. Moreover, we still have that $\abs{\chi(\hat\Sigma_n)}\leq G$.
	
	Let $\mu\eqdef \min\set{\mu_3,\text{ inj}_M(\partial M_i)}$ then for all $n$ the $\mu$-thick part of $S_n$ has a component intersecting $\partial M_i$. Moreover, each such component contains the image $f_n(\hat P_n)$ of a pair of pants that has $\gamma^n_1$ in its boundary. Then, by the Bounded diameter Lemma \cite{Th1978} we have that $\hat P_n$ has diameter bounded by $D_1\eqdef D_1(G,\mu)$ and let $D_2$ be maximal diameter of a component of $\partial M_i$.

	Consider the loops $\set{\alpha_n}_{n\in\N}$ that are contained in the surfaces $F_n$ with $i(\alpha_n,\gamma_1^n)=2$ and such that $N_r(\alpha_n\cup\gamma_1^n)\subset F_n$, $r>0$, is an essential four punctured sphere. Then, the $\set{\alpha_n}_{n\in\N}$'s have representatives in $M$ whose length is bounded by:
	$$\ell_M(\alpha_n)\leq D\eqdef2( D_1+D_2)$$
this is because we can push $\alpha_n$ to be obtained as two arcs $\beta_n,\beta_n'\subset S_n$ connected by two arcs in $\partial_bM_i$.

	Therefore, by discreetness of $\Gamma$ we have that they are in finitely many homotopy classes. Thus, we have a subsequence $\set{\alpha_{n_k}}_{k\in\N}$ such that for all $k, h\in\N$ $\alpha_{n_h}\simeq \alpha_{n_k}$ in $M\setminus \text{int}(M_i)$. Thus, by taking a sub-sequence of embedded annuli connecting $\alpha_{n_1}$ to $\alpha_{n_k}$ we obtain an embedded product with base a neighbourhood of $\alpha_{n_1}$ in $\Sigma_{n_1}$. This means that in the compactification $\overline M$ we have that $\gamma_1$ was not peripheral since it is a separating curve of a four punctured sphere embedded in $\partial\overline M$.\epf

In order to prove the characterisation Theorem \ref{maintheorem} we only need to show that if $\overline M\in\cat{Bord}(M)$, for  $M\in\M^B$, does not have any double peripheral cylinder $C$ then $M$ admits a complete hyperbolic metric.

\subsection{Relatively Acylindrical are Hyperbolic}\label{section3}

Now that we have completed the necessary topological construction we can show that manifolds $M\in\M^B$ with $(\overline M,P)$ an infinite-type pared acylindrical 3-manifold admit a complete hyperbolic metric. We first show that such manifolds are homotopic equivalent to a complete hyperbolic manifold $N$ with $\pi_1(
P)$ represented by parabolic elements. We achieve this by using the relative compactness of algebraic sequences developed by Thurston. Afterwards, with techniques similar to \cite{SS2013}, we show that the homotopy equivalence is homotopy to a homeomorphism $\psi:M\cong N$. 

\subsubsection{Relatively acylindrical are homotopy equivalent to hyperbolic}
To prove the homotopy equivalence we need a couple of technical result about sequences of non-elementary representations.
	\blem\label{fact1}
	Let $\rho_n:G\rar \text{PSL}_2(\C)$ be non-elementary discrete and faithful representations such that $\rho_n\rar \rho$ and let $\set{g_n}_{n\in\N}\subset \text{PSL}_2(\C)$. Then,  if we have a converging sub-sequence $g_{n_k}\rho_{n_k} g_{n_k}^{-1}\rar \rho'$ we have that up to an ulterior sub-sequence: $g_{n_k'}\rar g$ and $\rho'=g\rho g^{-1}$. The converse also holds.
	\elem 
	\bpf
	If the $\set{g_n}_{n\in\N}$ have a converging subsequence we are done. So assume that $g_n\rho_n g_n$ has a converging subsequence, which we denote by $g_n\rho_n g_n\rar \rho'$. Since the $g_n\rho_n g_n$ are non-elementary their algebraic limit $\rho'$ is non-elementary as well (see \cite{JM1990}). Therefore, we can find $\alpha,\beta\in G$ loxodromic elements that generate a discrete free subgroup $\langle \rho'(\alpha),\rho'(\beta)\rangle$. By algebraic convergence we have $g_n\rho_n(\alpha)g_n^{-1}\rar \rho'(\alpha)$, $g_n\rho_n(\beta)g^{-1}_n\rar\rho'(\beta)$ with $\rho_n(\alpha)\rar \rho(\alpha)$ and $\rho_n(\beta)\rar\rho(\beta)$.
	
	Since traces are preserved under conjugation and we assumed that $\rho'(\alpha),\rho'(\beta)$ were loxodromic so are $\rho(\alpha)$, $\rho(\beta)$. Denote by $x^\pm_\infty$, $y^\pm_\infty$ the attracting/repelling fixed points of $\rho(\alpha),\rho(\beta)$ respectively and similarly define $a^\pm,b^\pm$ for $\rho'(\alpha),\rho'(\beta)$. Moreover, we have that eventually $g_n\rho_n(\alpha)g_n^{-1}$ are all loxodromic and similarly for $g_n\rho_n(\beta)g_n^{-1}$. Therefore, the attracting (repelling) fixed points of $g_n\rho_n(\alpha)g_n^{-1}$ converge to the attracting (repelling) fixed point of $\rho'(\alpha)$. The fixed points of $g_n\rho_n(\alpha)g_n^{-1}$ are $g_n(x^\pm_n)$ for $x^\pm_n$ the fixed point of $\rho_n(\alpha)$, hence we have that $x^\pm_n\rar x^\pm_\infty$ and $g_n(x^\pm_n)\rar a^\pm$. By triangle inequality:
		\begin{align*} d_{\bH}(g_n x^\pm_\infty, a^\pm)&\leq d_{\bH}(g_nx^\pm_\infty,g_nx^\pm_n)+d_{\bH}(g_nx^\pm_n,a^\pm)\\&=d_{\bH}(x^\pm_\infty,x^\pm_n)+d_{\bH}(g_nx^\pm_n,a^\pm)\end{align*}
			Thus it follows that $g_n(x^\pm_\infty)\rar a^\pm$ and this also holds for the $y^\pm_\infty$ and $b^\pm$. Since $\langle \rho(\alpha),\rho(\beta) \rangle$ is discrete at least three of the $x^\pm_\infty$ and $y^\pm_\infty$ are distinct. Then, by Theorem \cite[3.6.5]{Be1983} we have that the $\set{g_n}_{n\in\N}$ form a normal family.	\epf

	\bprop\label{compactness conv}
	Let $\rho_i:G=\cup_{j=1}^\infty G_j\rar \text{PSL}_2(\C)$ be non-elementary discrete and faithful representations and let $\rho^j_i$ be the restriction of $ \rho_i$ to ${G_j}$. Then, given $g^j_i\in \text{PSL}_2(\C)$ such that $\forall i,j: g^j_i \rho^j_i (g^j_i)^{-1}:G_j\rar \text{PSL}_2(\C)$ converge up to subsequence we have that $\forall j:g^1_i \rho^j_i(g_i^1)^{-1}$ converge up to subsequence.
	\eprop	\bpf
		
	We first show that $\forall j:\set{g^1_i(g_i^j)^{-1}}_{i\geq j}$ converge up to subsequence. For all $j$ consider:
	\be g^1_i \rho^1_i(g_i^1)^{-1}  =\left( g^1_i (g^j_i)^{-1}\right) g^j_i \rho^1_i (g^j_i)^{-1} \left(g^j_i(g^1_i)^{-1}\right)\qquad i\geq j\ee
	By assumption we have the the left-hand side has a converging subsequence which we call $i_n$. Since the $\set{g^j_i \rho^j_i (g^j_i)^{-1}}_{i\geq j}$ have a converging subsequence so do their restrictions on $G_1$: $\set{g^j_{i} \rho^1_{i} (g^j_{i})^{-1}}_{i\geq j}$. Therefore, we can extract another subsequence $\set{i_n'}_{n\in\N}$ such that both $ g^1_{i_n'} \rho^1_{i_n'}(g_{i_n'}^1)^{-1} $ and $g^j_{i_n'}\rho^1_{i_n'} (g^j_{i_n'})^{-1}$ are converging. Then we are in the setting of Lemma \ref{fact1} thus, we have that $g^1_{i_n'} (g^j_{i_n'})^{-1}$ are converging, up to an ulterior subsequence, as well and we call the limit $g^j_1$. Since we are only concerned about subsequences we assume that $g^1_i (g^j_i)^{-1}$ is the converging subsequence for which also $ g^j_i \rho^j_i (g^j_i)^{-1} $ is converging. Then we have:
	\be \forall j: g^1_i \rho^j_i(g_i^1)^{-1}=g^1_i (g^j_i)^{-1}\left( g^j_i \rho^j_i (g^j_i)^{-1}\right) g^j_i (g^1_i )^{-1}\ee
	but everything on the right hand side is converging, hence the left hand side does. Since the $j$ was arbitrary this concludes the proof.\epf
	
									We recall the following Theorem by Thurston \cite[0.1]{Th1998b} for hyperbolizable compact pared 3-manifolds:
									
									\begin{Theorem*}\label{algcomp} Let $(M,P)$ be a pared compact hyperbolizable 3-manifold. Then the set of representations induced by $AH(M,P)$ on the fundamental group of any component of $M^{\text{acyl rel } P}$ $AH(M,P)$ is compact up to conjugation.
									\end{Theorem*}
									
						Where he shows that given a sequence of discrete and faithful representations $\set{\rho_n}_{n\in\N}$ of an acylindrical 3-manifold we can find elements $\set{g_n}_{n\in\N}$ of $\text{PSL}_2(\C)$ so that the sequence $\set{g_n\rho_ng^{-1}_n}_{n\in\N}$ has a converging subsequence. 
						
						Then for $(\overline M,P)$ we have:

								 \bthm\label{hypmetric}
									Let $(\overline M,P)$ be an infinite-type pared acylindrical hyperbolic 3-manifold for $\overline M$ the bordification of a manifold in $\M^B$. Then $ M$ is homotopic to a complete hyperbolic 3-manifold $N$ such that $P$ is represented by parabolic elements.
									\ethm
									\bpf Let $\set{ M_i}_{i\in\N}$ be the exhaustion of $M$, then by Proposition \ref{eventualacylindricalcomponents} for each $i$ we can find $n_i$ such that $M_i\subset  M_{n_i}^{\text{acyl rel }P}$. Since each $ M_i$ is hyperbolizable we have a discrete and faithful representation $\rho_i\in AH( M_i,P)$.									
									Let $X_i\eqdef  M_{n_i}^{\text{acyl rel }P}$, then by Theorem \cite[0.1]{Th1998b} applied to the sequence $\set{\rho_k\vert_{\pi_1(X_i)}}_{k\geq n_i}$ we can find $\set{g_k^j}\subset PSL_2(\C)$ such that for all $j$ the restriction of $\set{g_k ^j\rho_k\left( g_k^j\right)^{-1}}$ to $\pi_1( M_j)\subgroup \pi_1( X_j)$ have a converging subsequence. By Proposition \ref{compactness conv} we can assume that the $g^j_k$ do not depend on $j$ so that we have representations $\set{g_k\rho_kg_k^{-1}}_{k\in\N}$ that subconverge on each $\pi_1(M_j)$. By picking a diagonal subsequence we can define:
									\be \forall\gamma\in\pi_1( M),\gamma\in\pi_1( M_i):\rho_\infty(\gamma)\eqdef\lim_{n\geq n_i}g_n \rho_n(\gamma)g_n^{-1}\ee
									Since $\pi_1(  M)=\cup_{i\in\N}\pi_1( M_i)$ we get a representation $\rho_\infty:\pi_1( M)\rar PSL_2(\C)$ which is discrete and faithful by \cite{JM1990}. 
									
									Thus if we define $N\eqdef \quotient{\bH}{\rho_\infty(\pi_1( M))}$ we have $\pi\eqdef\pi_1(N)\cong \pi_1( M)$ and since they are both $K(\pi,1)$ there is a homotopy equivalence between them. By construction all elements of $P$ are parabolic in $\rho_\infty$. \epf

									We define:
									
									\bdefi
									Let $(\overline M,P)$ with $\overline M\in \cat{Bord}(M)$ an infinite type pared acylindrical 3-manifold and $N\cong\hyp\Gamma$ be a hyperbolic 3-manifold. Then, a homotopy equivalence $\phi:M\rar N$ is said to \emph{preserve parabolics} if $\forall \gamma\in\pi_1(M)$ homotopic in a component of $P$ we have that $\phi_*(\gamma)$ is represented by a parabolic element in $\Gamma$.								\edefi

									\blem\label{finiteapprox}
									Let $(\overline M,P)$ be a pared infinite type acylindrical 3-manifold for $M\in\M^B$ and let $\phi:M\rar N$ be a homotopy equivalence preserving $P$. Then, for all $n\in\N$ we have that the cover $N_n\rar N$ corresponding to $\phi_*(\pi_1(M_n))$ the lift: $\phi: M_n\rar N_n$ is homotopic to an embedding and $N_n\cong \text{int}(M_n)$.									\elem
									\bpf

									By Proposition \ref{eventualacylindricalcomponents} for all $n$ we have $k_n$ such that $M_n\subset M_{k_n}^{\text{acyl rel }P}$ and let $P_{k_n}$ be the annuli in $\partial M_{k_n}$ induced by $P$. Consider the cover $N_{k_n}\rar N$ corresponding to $\pi_1(M_{k_n})$ and let $N_{k_n}'$ be its manifold compactification, which exists by Tameness \cite{AG2004,CG2006}, and let $Q\subset\partial N'_{k_n}$ the parabolic locus. Since $\phi$ preserves parabolics we can homotope $\tilde\phi\vert_{P_{k_n}}:P_{k_n}\rar Q$ to be a homeomorphism onto its image. Then, by Lemma \ref{reljoh} the homotopy equivalence:
									$$ \tilde\phi: M_{k_n}\rar N_{k_n}'$$
									is homotopic to a map $\psi$ that is an embedding on $M_n$. Then, $\psi\vert_{M_n}\simeq\phi\vert_{M_n}$ lifts to the cover $N_n\rar N_{k_n}\rar N$ and its image forms a Scott core for $N_n$. Since the homotopy equivalence $\tilde\psi\vert_{M_n}: M_n\rar N_n$ is an embedding and $\tilde\psi(\partial M_n)$ is incompressible by Lemma \ref{scottcore} we get that $N_n\cong \text{int}(M_n)$. \epf

			 \subsubsection{Relatively hyperbolic are homeomorphic to hyperbolic}

	A key step in the proof of Theorem \ref{maintheorem} is that the homotopy equivalence $\phi: M\rar N$ mapping the elements corresponding to $P\subset\partial\overline M$ to parabolics in $N$ is homotopic to a proper homotopy equivalence and it embeds the boundary components of a subsequence of a minimal exhaustion\footnote{See Definition \ref{minimalexhdef}.} $\set{M_n}_{n\in\N}$.
	
	Our first objective is to show the following Theorem:
		\bthm
Given the maximal bordification $\overline M\in\cat{Bord}(M)$ of $M\in\M^B$ and a minimal exhaustion $\set{ M_i}_{i\in\N}$ , if $\overline M$ has no doubly peripheral cylinders let $P\hookrightarrow \partial\overline M$ be such that $(\overline M,P)$ is an infinite-type pared acylindrical manifold. Then, for any hyperbolic $N$ and $\phi:  M\rar N$ a homotopy equivalence preserving $P$ we have a proper homotopy equivalence $\widehat\phi: M\rar N$ preserving $P$ such that $\widehat\phi$ is a proper embedding on tame ends of $M$ and on $\mathcal S\eqdef\cup_{i\in\N} \partial M_{a_i}$ for $\set{a_i}_{i\in\N}$ an increasing subsequence.
\ethm

Before doing a full proof we deal with a couple of preliminary Lemmata. The first Lemma says that if $\overline M$ is a maximal bordification induced by a maximal product $\P_{max}:S\times[0,\infty)\hookrightarrow M$ in standard form then we cannot have products $\P:\mathbb A\times [0,\infty)\hookrightarrow M$, also in standard form, such that  for $k$ sufficiently large the components of $\im(\P)\cap\cup_{k\in\N}\partial M_k$ are not peripheral in $\partial M_k\setminus \im(\P_{max})$.

\blem\label{noperipheralproducts}
Let $\P_{max}$ be a maximal product in $M\in\M$ and let $\P:\mathbb A\times[0,\infty)\hookrightarrow M$ be a product. If they are both in standard form then there exists some $i\in\N$ such that for $k>i$ the component of intersections of $\partial M_k\cap \im(\P)$ are peripheral in $\partial M_k\setminus \im(\P_{max})$.
\elem
\bpf By maximality of $\P_{max}$ we have a proper isotopy of $\P$ and a connected sub-product $\mathcal Q$ of $\P_{max}$ such that $\P$ is contained in an $r$-neighbourhood $\mathcal Q'$ of $\mathcal Q$. Moreover, without loss of generality we can also assume that this new $\P$ is in standard form. Let $i$ be the minimal $i$ such that $\im(\P)\cap\partial M_i\neq\emp$. Then, for $k\geq i$ all components of $\partial M_k\cap \im(\P)$ are isotopic in $\mathcal Q'$ into $\mathcal Q\cap S$ and so they are peripheral in $\partial M_k\setminus \im(\P_{max})$ reaching a contradiction. \epf

The next Lemma says that tame ends of $M$ relative to $P$ embed, up to homotopy, into $N$ for any homotopy equivalence $\phi:M\rar N$ preserving $P$.

\blem\label{tameends} Let $(\overline M, P)$, with $\overline M\in \cat{Bord}(M)$ and $M\in\M^B$, be an infinite-type pared acylindrical 3-manifold,  $N$ be a hyperbolic 3-manifold, $\phi:M\rar N$ be a homotopy equivalence preserving $P$ and let $\P_{max}:S\times[0,\infty)\hookrightarrow M$ be a maximal product in $M$ inducing $\partial\overline M$. Given $\epsilon<\mu_3$\footnote{For $\mu_3$ the 3-dimensional Margulis constant.} we have a homotopy equivalence $\psi:M\rar N$, homotopic to $\phi$, that is a proper embedding on $\im(\P_{max})$ and with the property that $\psi\circ\P_{max}(\partial S\times[0,\infty))\subset \partial  Q_\epsilon$ for $ Q_\epsilon$ the $\epsilon$-boundary of the parabolic locus of $N$.
\elem
\bpf
Let $\P'\subset \im(\P_{max})$ be a regular neighbourhood of $P\subset\partial \overline M$ in $\im(\P_{max})$ and let $N^0\eqdef N\setminus \text{int}(Q_\epsilon)$. We will now show that $\phi$ can be homotoped to be an embedding on $\P'$. 

\vspace{0.3cm} 

\paragraph{Step 1:} We can homotope $\phi$ to be a homeomorphism on $\P'$ mapping $\partial \P'$ onto $\partial Q_\epsilon$.

\vspace{0.3cm}

Each component of $\P'$ is homeomorphic to either $\mathbb T^2\times [0,\infty)$ or $\mathbb A\times[0,\infty)$. Since $(\overline M, P)$ is relatively acylindrical no two components of $\P'$ are mapped by $\phi$ to the same component of $Q_\epsilon$. If a component of $\P'$ is homeomorphic to $\mathbb T^2\times [0,\infty)$ the fact that $\phi$ can be homotoped to be an embedding follows directly by the hyperbolicity of $N$. So from now on we only consider components $P_i$ of $\P'$ homeomorphic to $\mathbb A\times [0,\infty)$. Since $\phi$ preserves parabolics, the image of the fundamental group of any component $P_i$ of $\P'$ is contained in a parabolic subgroup $\langle\gamma\rangle$ in $\pi_1(N)$ for $\gamma$ a primitive element. Moreover, since $(\overline M, P)$ is an infinite-type pared acylindrical 3-manifold any generator of $\pi_1(P_i)$ is primitive in $\pi_1(M)$ and hence so is its image in $\pi_1(N)$ thus $\phi_*:\pi_1(P)\rar \langle\gamma\rangle$ is an isomorphism.

Therefore, we see that the core of $P_i$ is homotopic through $\phi$ to a simple closed curve $\gamma$ in $\partial Q_\epsilon$. We claim that $\gamma$ is contained in an annular component of $\partial Q_\epsilon$. If $\gamma$ is in a torus component then we would have that $\pi_1(P)$ is contained in a $\Z^2$ factor in $\pi_1( M)$ which again is not possible by the fact that $(\overline M, P)$ is an infinite-type pared acylindrical 3-manifold. Since $\gamma$ is in an annular component the claim follows and we can assume that, up to a homotopy, $\phi$ is a homeomorphism on $\P'$. Moreover, we can assume that $\phi(\partial \P')\subset\partial  Q_\epsilon$.

\vspace{0.3cm}

Let $S\subset M$ be a push-in of a component of $\partial\overline M\setminus P$. Since $S$ is incompressible and $\phi$ is a homotopy equivalencem up to a homotopy of $\phi$ that is constant on $\P'$, we have that $\phi\vert_S:(S,\partial S)\rar (N^0,\partial Q_\epsilon)$ and let $U\eqdef N_r(\phi(S))\subset N^0$ be a regular neighbourhood of $\phi(S)$. By the existence results for PL-least area surfaces in \cite{JaRu1988} and \cite[1.26]{Kap2001} for a triangulation $\tau_S$ of $N^0$ there are PL-least area representatives $S'$ of $\phi(S)$ such that $S'\subset U$. Consider the cover $\pi_j:N_j\rar N^0$ corresponding to $\phi_*(\pi_1( M_j))$ such that the set $U$ lifts homeomorphically to $\tilde U$ in the cover and denote by $\tilde S'$ the lift of $S'$. By Lemma \ref{finiteapprox}, up to homotopy we have that the homotopy equivalence $\tilde\phi: M_j\rar N_j$ is an embedding, thus we see that all the $\tilde S'$ are homotopic to embedded surfaces. Hence, by the results of \cite{JaRu1988} we have that the $S'$'s are embedded as well. Since the covering projection $\pi_j$ is a homeomorphism on $\tilde U$ with image $ U$ we get that $S'=\pi_j(\tilde S')$ is embedded as well.

Therefore, up to a homotopy of $\phi$ we can assume that $\phi$ embeds in $N^0$ a collection of surfaces $\mathcal S\eqdef\set{S_n}_{n\in\N}$ that are pushed-in of components of $\partial\overline M\setminus\P$. Since $[S]\in H_2(M,\partial \P')$ is a separating $\pi_1$-injective surface so is $[\phi(S)]\in H_2(N,\partial Q_\epsilon)$.

Let $\set{E_n}_{n\in\N}$ be the collection of tame ends of $M\setminus\P'$ facing the surfaces $S_n$ and let $\Sigma_n\eqdef \phi(S_n)$.

\vspace{0.3cm}

\paragraph{Step 2:} Up to a homotopy of $\phi$ we have that for all $n$ $\phi:E_n\rar N$ is an embedding, $\phi(E_n)$ are pairwise disjoint and $\phi(\partial\P')\subset \partial Q_\epsilon$.

\vspace{0.3cm}

 We first show that each $E_n$ embeds. Let $X_1,X_2$ be the connected components of $N^0\vert\Sigma_n$ and assume that $\pi_1(\Sigma_n)\hookrightarrow \pi_1(X_h)$ is not a surjection for $h=1,2$. Thus, neither one of $X_1$, $X_2$ is homeomorphic to $\Sigma_n\times[0,\infty)$. Then, we can find some compact submanifold $K_n\subset N$ containing $\Sigma_n$ such that $K_n$ is not homeomorphic to a product and it contains topology on both sides of $\Sigma_n$. Let $m>n$ such that $\pi_1(K_n)\subset\phi_*(\pi_1(M_m))$ and consider the cover $N_m\rar N$ corresponding to $\phi_*(\pi_1(M_m))$. By Lemma \ref{finiteapprox} we have that $\tilde\phi:M_m\rar N_m$ is homotopic to a homeomorphism $\Phi:\text{int}(M_m)\rar N_m$. In $M_m$ the incompressible surfaces $S_n,\Phi^{-1}(\Sigma_n)$ are homotopic and since they are incompressible they are isotopic by \cite{Wa1968}. Thus, at least one of $X_1,X_2$ is homeomorphic to a product and so we can homotope $\phi$ rel $\P'$ so that $\phi:E_n\rar N$ is an embedding and we denote by $\mathcal E_n$ its image. Therefore, up to a homotopy of $\phi$ we get that all ends $E_n$ embed in $N$.

We now need to show that up to an ulterior homotopy we have that the $\phi(E_n)$'s are pairwise disjoint. Up to relabelling the $E_n$ we can assume that for all $n$ there exists $k_n$ such that the ends $E_1,\dotsc, E_n$ are ends of $\hat M_{k_n}\setminus\P'$ for $\hat M_{k_n}\cong M_{k_n}\cup\partial_\infty 
M_{k_n}\times[0,\infty)$. Consider the covers $N_{k_n}$ of $N$ corresponding to $\phi_*(\pi_1(M_{k_n}))$ and denote the homeomorphic lifts of $\mathcal E_n$ in $N_{k_n}$ by $\tilde{\mathcal E}_n$. By Lemma \ref{finiteapprox} $\phi:\text{int}(M_{k_n})\rar N_{k_n}$ is homotopic to a homeomorphism. Hence, all ends $E_1,\dotsc, 
E_n$ are mapped to distinct ends in $N_{k_n}$ and so the $\tilde{\mathcal E}_n$ correspond to distinct ends of $N^0_{k_n}$. Therefore, up to pushing the $\phi(E_n)$ inside the $\tilde{\mathcal E}_n$ we can assume that $\phi: M_{k_n}\rar N_{k_n}$ is an embedding on $E_1\cup \dotsc\cup E_n$. Since the projection is a homeomorphism on $\tilde {\mathcal E}_n$ and maps $\phi(\partial \mathcal P')$ into $\partial Q_\epsilon$ by iterating this construction we conclude the proof.\epf

Before showing that if we have a homotopy equivalence $\phi:M\rar N$ respecting parabolics we can homotope it so that it is a proper homotopy equivalence we need to understand how loops in components $S$ of $\cup_{k\in\N} \partial M_k$ are homotopic into $P$. 

\bdefi
Let $\P_m:S\hookrightarrow M$, $M\in\M^B$, be a maximal product in standard form with respect to a normal family $\set{(N_k,R_k)}_{k\in\N}$ and such that $\P_m(S\times\set 0)\subset\cup_{k\in\N}\partial M_k$ and let $\P_s\eqdef\P_m(\partial S\times[0,\infty))$. We say that an essential torus $V\subset N_k\setminus \text{int}(\im(\P_m))$ is a \emph{parabolic solid torus} (PST) if:
\begin{itemize}
\item[(i)] no annulus of $\partial V\setminus R_k$ is parallel to $\P_s$;
\item[(ii)] $\partial V\cap R_k$ has a component that is isotopic into $N_\epsilon (\P_s)\cap \partial X_k$ in $\partial X_k$;
\item[(iii)] $V$ is maximal with respect to (ii), i.e. if $V,Q\subset N\in\pi_0(N_k)$ are both PST with an isotopic wing then $V=Q$.
\end{itemize}
\edefi		
		\brem
		If $M$ does not have any doubly peripheral cylinder, then if $V$ is a PST we have that by (ii) the component of $V\cap R_k$ isotopic into $N_\epsilon (\P_s)\cap \partial X_k$ is unique.
		\erem

				\bdefi
				Let $\set{N_k}_{k\in\N}$ be a normal family of JSJ for $M=\cup_{k\in\N}M_k\in\M^B$ and let $V\subset N_k$ be a parabolic solid torus. We define a \emph{maximal parabolic solid torus} (MPST) $\hat V$ as the direct limit $\varinjlim V_i$ where:
				\begin{itemize}
				\item[(i)] $V_1=V$;
				\item[(ii)] $V_i\setminus V_{i-1}$ are essential solid tori contained in $N_j$, $j\in\N$, whose wings wrap once around the soul;
				\item[(iii)] $V_i$ is obtained from $V_{i-1}$ by adding all essential solid tori $Q\subset N_j,j\in\N$, that have a wing matching up with one of $V_{i-1}$ and such that $\partial Q\setminus R_j$ has no annuli parallel to $\P_s$ in $X_j$.
				\end{itemize}
				\edefi
				\brem
				Let $\hat V$ be a maximal parabolic solid torus, since for all $k$ we have that by (ii) $V\cap X_k$ are essential tori by Lemma \ref{charlimitsolidtori} we get that $\hat V\cong \mathbb S^1\times\mathbb D^2\setminus L$ for $L\subset\partial( \mathbb S^1\times\mathbb D^2)$ a collection of parallel simple closed curves.
				\erem
				We now show that if we do not have doubly peripheral cylinders then maximal parabolic solid tori are compact.
				
				\blem\label{mpstcompact}
				Let $V$ be a parabolic solid torus in $M\in\M^B$. If $M$ does not have any doubly peripheral annulus then $\hat V$ is compact.
				\elem
				\bpf
				If $\hat V$ is not compact by Lemma \ref{charlimitsolidtori} we have that $\hat V\cong \mathbb S^1\times\mathbb D^2\setminus L$ for $L\neq\partial ( \mathbb S^1\times\mathbb D^2)$ and $L\neq \emp$. Thus, we obtain a product in standard form $\P:\mathbb A\times[0,\infty)\hookrightarrow M$ whose image is contained in $\hat V$. By construction, in particular property (iii), we have that no annular component of $\partial\hat V\cap X_k$ is parallel in $X_k$ into $\P_s$, thus no component of $\im(\P)\cap X_k$ is isotopic in $X_k$ into $\P_m$. Therefore, by Lemma \ref{noperipheralproducts} we reach a contradiction with the fact that $\P_m$ was maximal. 
				\epf
						 		
			We now show that if $M\in\M^B$ has no double peripheral annuli maximal parabolic solid tori corresponding to distinct parabolic solid tori are disjoint.

			\blem\label{distincmpst} Let $V\neq Q$ be disjoint parabolic solid tori contained in $X_k,X_j$ respectively and assume that $M$ has no doubly peripheral annuli. Then, $\hat V\cap\hat Q=\emp$.
			\elem				
				\bpf If $\hat V\cap \hat Q\neq\emp$ then by construction we get that $\hat V=\hat Q$. Let $A_1\subset \partial V\cap R_k$ and $A_2\subset \partial Q\cap R_j$ be the annuli isotopic into $N_\epsilon(\P_s)\cap X_k$ and $N_\epsilon(\P_s)\cap X_j$ respectively. Since $V\neq Q$ by (iii) of the definition of PST we get that if $j=k$ then $A_1$ and $A_2$ are non-isotopic annuli in $R_k$. 
				
				Since $\hat V=\hat Q$ by (ii) of the definition we have an annulus $C\subset \hat V$ connecting $A_1$ to $A_2$. By extending the annulus $C$ to an annulus $\hat C$ by going to infinity along the components of $\P_s$ that $A_1,A_2$ are homotopic to we get a properly embedded annulus $\hat C\subset M$ which compactifies to an annulus $\overline C$ in $\overline M$.
				
				\vspace{0.3cm}
				
				\textbf{Claim:} The annulus $\overline C$ is essential.
				
				\vspace{0.3cm}
				
				\bpfc
				If $A_1,A_2$ are isotopic into distinct components of $\P_s$ then $\overline C$ is essential in $\overline M$ and we are done. Thus, we can assume that $\partial\overline C$ are isotopic in $\partial \overline M$. Then, if $\overline C$ is $\partial$-parallel we have $k\in\N$ such that $\overline C\cap X_k$ is isotopic into $\P_s$ contradicting the construction of $\hat V$ and $\hat Q$ or the fact that $A_1$ was not isotopic to $A_2$ in $X_k$ and so that $V$ and $ Q$ were distinct parabolic solid tori (we contradict property (iii)).
				\epfc
				
				Thus, since $\overline C$ is essential and has both boundaries peripheral in $\partial\overline M$ we get that $\overline M$ has a doubly peripheral annulus reaching a contradiction.\epf

				Our final preparatory Lemma is:
				\blem\label{noperipheral}
Let $\P_m\subset M\in\M^B$ be the image of a maximal product in standard form and let $\mathcal V$ be a maximal collection of MPST. Then, for $S\in\pi_0(\cup_{k\in\N} \partial M_k\setminus (\P_m\cup\mathcal V))$ we have that any essential non-peripheral simple loop $\gamma\subset S$ is not homotopic into $\P_s$, for $\P_s$ the side boundary of $\P_m$.
				\elem								
				\bpf Since $\gamma\subset S\in\pi_0(\cup_{k\in\N} \partial M_k\setminus (\P_m\cup\mathcal V))$ is non-peripheral it is not isotopic into any torus of $\mathcal V$. Let $H$ be the cylinder connecting $\gamma$ to $\P_s$, up to a homotopy of $H$ rel $\partial H$ we can assume that for all $k\in\N$ $H\cap X_k$ is essential. By an iterative argument and the Annulus Theorem we have that for all $k$ $H\cap X_k$ is homotopic to an embedded annulus. Then a thickening $P$ of $H$ is a PST and since $\gamma\cap \mathcal V=\emp$ by Lemma \ref{distincmpst} we have that $\hat P\cap \mathcal V=\emp$ contradicting the maximality of $\mathcal V$.
				\epf								
					Then, by working gap by gap and the above Lemma we have that:
					
					\bcor\label{MPST} Let $M\in\M^B$ and $\P:\Sigma\times[0,\infty)\hookrightarrow M$ be a maximal product. Then, there exists a maximal collection $\hat V$ of pairwise disjoint MPST.
					\ecor

				We can now show that given a homotopy equivalence $\phi$ between the interior of an infinite-type acylindrical pared 3-manifold and a hyperbolic 3-manifold respecting parabolics we have that $\phi$ is homotopic to a proper homotopy equivalence.

\bthm\label{homotopyequivalence} Let $(\overline M, P)$ with $\overline M\in \cat{Bord}(M)$ and $M\in\M^B$ an infinite-type pared acylindrical 3-manifold. Then, there exists a complete hyperbolic 3-manifold $N$ and a proper homotopy equivalence $\phi:M\rar N$ respecting $P$ such that $\phi$ is an embedding on $\mathcal S\eqdef\cup_{i\in\N} \partial M_{a_i}$ for $\set{a_i}_{i\in\N}$ an increasing subsequence. Moreover, we can also assume $\phi$ to be a proper embedding on any tame end of $M$.
\ethm

\bpf 
By Theorem \ref{hypmetric} we have a homotopy equivalence $\phi:M\rar N$ respecting $P$. Let $\P_{max}$ be a maximal product inducing $\overline M$, $\P'\subset \im(\P_{max})$ be a neighbourhood of $P$ and for $\epsilon<\mu_3$ let $Q_\epsilon$ be the $\epsilon$-thin part of $N$. By Lemma \ref{tameends} we have:

\vspace{0.3cm}

\textbf{Step 1:} Up to a homotopy of $\phi$ we can assume that $\phi\vert_{\im(\P_{max})}$ is an embedding and that it maps $\partial \P'$  into $\partial Q_\epsilon$.

\vspace{0.3cm}

Let $\mathcal S\eqdef \cup_{k\in\N}\partial M_k\setminus \im(\P_{max})$, since $\im(\P_{max})\cap\cup_{k\in\N}\partial M_k$ are images of level surfaces and are essential in $\partial M_k$ we have that every component of $\mathcal S$ is an essential subsurface of $\partial M_k$ whose boundary is contained in $\partial\P'$. Then, for any component $S$ of $\mathcal S$ we have that $\phi\vert_S: S\rar N$ maps $\partial S$ homeomorphically into $\partial Q_\epsilon$ and without loss of generality we can assume each component of $\phi(\partial S)$ to be a horocycle in $\partial Q_\epsilon$. 

Let $\hat {\mathcal V}\subset M$ be a maximal collection of pairwise disjoint MPST, see Corollary \ref{MPST}, and define:
$$\hat {\P}\eqdef N_r(\im(\P_{max}))\cup\hat {\mathcal V}$$
Since each $V\in\pi_0(\hat{\mathcal V})$ is homeomorphic to a solid torus $\mathbb S^1\times\mathbb D^2 $ then $\hat{\P}\cong\im(\P_{max})$ and for all $k\in\N$: $\partial M_k\setminus \hat{\P}$ is a collection of essential subsurfaces of $\partial M_k$, thus since $M\in\M^B$ they have a uniform bound on their complexity:
$$\forall k,\forall S\in\pi_0(\partial M_k\setminus \hat{\P}): \abs{\chi(S)}\leq G$$
\noindent and we define $\hat{\mathcal S}\eqdef \cup_{k\in\N}\partial M_k\setminus \hat{\P}$.

Since $(\overline M, P)$ is an infinite-type acylindrical 3-manifold we have that no component $S$ of $\hat{\mathcal S}$ is an annulus. If it where then, by going through the boundary of $\hat \P\iso \im(\P_{max})$ we would have a doubly peripheral essential annulus.


\vspace{0.3cm}

\textbf{Step 2:} There exists a proper homotopy equivalence $\psi: M\rar  N$ with $\psi\simeq \phi$.

\vspace{0.3cm}

The aim will be to show that up to homotopy we have that $\phi$ is proper when restricted on $\cup_{k\in\N}\partial M_k$. We first show that we can make $\phi\vert_{\hat{\mathcal S}}$ proper and then by doing homotopies of the annuli of $\mathcal S\setminus\hat{\mathcal S}$ we will get $\phi$ is proper when restricted on $\mathcal S$ and then by doing homotopies in the tame ends of $N$ we will obtain the required result.

\vspace{0.3cm}

\textbf{Claim:} Up to homotopy $\psi\vert_{\hat{\mathcal S}}$ is proper map.

\vspace{0.3cm}

For any essential subsurface $S\subset \partial M_k$ in $\pi_0(\hat{\mathcal S}) $ we can pick a triangulation $\tau$ such that each component of $\partial S$ is realised as a single edge in $\tau$ and all vertices are contained in $\partial S$. Each component of $\phi(\partial S)$ is homotopic into a unique component of $ \partial Q_\epsilon$. Let $S'$ be an open regular neighbourhood of $S$ in $\partial M_k$, then we can homotope $\phi\vert_{S'}$ so that $\phi$ maps the cusps of $S'$ into cusps region contained in $\phi(\P_m)\cap Q_\epsilon$, i.e. into cusps of $N$. Since $\phi\vert_{S'}$ is type preserving proper map we can realise the proper homotopy class of $\phi(S')$ by a simplicial hyperbolic surface sending cusps to cusps, see \cite{Ca1996,Bo1986}. Moreover, we can do this for all $S$ in $\hat{\mathcal S}$ via a homotopy of $\phi$. With an abuse of notation we still denote by $\phi$ the resulting map. We now claim that $\phi$ is a proper map when restricted to $\hat{\mathcal S}=\set{\Sigma_k}_{k\in\N}$ whose image is contained in the simplicial hyperbolic surfaces $\set{S_k}_{k\in\N}$ we constructed.

If $\phi$ is not proper we can find a sequence $\set{p_k\in\Sigma_k}_{k\in\N}\subset\hat{\mathcal S}$ of points and surfaces such that for $i\neq k$ we have $\Sigma_i\neq\Sigma_k$ and $\phi(p_k)$ has a limit point $p\in N$. Each $\Sigma_k$ is contained in a simplicial hyperbolic surface $S_k$ of bounded topological type. Since the $S_k$ have uniformly bounded complexity by Gauss-Bonnet we get that their area is uniformly bounded by some $A\eqdef A(G)$.

				\vspace{0.3cm}					
						
\paragraph{Case 1:} There is a sub-sequence of the $\Sigma_k$ such that $\Sigma_k$ is not a pair of pants.

				\vspace{0.3cm}

Since $Area(S_k)\leq A$ and $S_k$ is not a pair of pants we can find a constant $D$ such that for all $k\in\N$ there is a non-peripheral essential simple closed loop $\gamma_k\subset S_k$ based at $p_k$ such that $\ell_N(\phi(\gamma_k))\leq D$. Since we assumed that the points $\phi(p_k)\rar p$ in $N^0$ we have that the  $\set{\phi(\gamma_k)}_{k\in\N}$ have to be in finitely many distinct homotopy classes. Since $\phi$ is a homotopy equivalence the same must happen to the $\set{\gamma_k}_{k\in\N}$.
									
									Then by picking a subsequence of the $\Sigma_k$'s we can assume that they all have a homotopic curve $\gamma$. This curve was essential and non-peripheral in each $\Sigma_k$ and so is not homotopic in $\partial \Sigma_k$, thus we get a product $\P$ that either is not contained in $\P_{max}$ and is not peripheral $\partial M_m\setminus\im(\P_{max})$ for $m\geq n$, for $n$ the smallest $m$ such that a component of $\set{\Sigma_k}_{k\in\N}$ is in $\partial M_m$, thus contradicting Lemma \ref{noperipheralproducts} or $\gamma$ is isotopic in $\P_s$ the side boundary of $\P_m$ contradicting Lemma \ref{noperipheral}.
				
				\vspace{0.3cm}					
									
\paragraph{Case 2:} All but finitely many $\Sigma_k$ are pair of pants.

				\vspace{0.3cm}					
						
			Let $\phi_k\eqdef \phi\vert_{S_k}$ be the simplicial hyperbolic surface corresponding to the thrice punctured sphere $\Sigma$, so that we have simplicial hyperbolic surfaces: 
			$$\psi_k:(\Sigma,p_k)\rar N$$
			 such that $\psi_k(p_k)\rar p\in N$. Since $inj_N(p)>0$, the $\phi_k$ are $1$-Lipschitz maps and $\liminf inj_{\Sigma}(p_k)>0$ the $\set{p_k}_{k\in\N}\geq inj_N(p)$ are contained in a compact core $K\subset \Sigma$, homeomorphic to a pair of pants. This, means that we can find a compact set $K'\subset N$ containing $p$ and with the property that for all $k\in\N$ $\psi_k(K)\subset K'$. Pick $i$ such that $\pi_1(K')\subset \pi_1(\phi_*(\pi_1(M_i)))$, then all the $\Sigma_k$ lift to the cover $M_i$ and we get that they are eventually parallel by the Kneser-Haken finiteness theroem, giving us a product over a pair of pants which cannot be properly isotopic into $\P_m$.
			 
			 				\vspace{0.3cm}

									Thus $\phi$ is a proper map when restricted on $\hat{\mathcal S}$ and $\im(\P_m)\cap\cup_{k\in\N} \partial M_k$ and every component $A$ of $\cup_{k\in\N}\partial M_k\setminus (\hat{\mathcal S}\cup\im(\P_m))$ is an annulus that is mapped into a cusp region of $N$. Then, by mapping the annuli further and further in the cusp end we obtain that $\phi$ is a proper map when restricted on $\cup_{k\in\N}\partial M_k$.


									Since the restriction of $\psi$ to $\mathcal S'\eqdef \cup_{k\in\N}M_k$ is a proper map for every compact set $K\subset N$ the preimage $\phi\vert_{\mathcal S'}^{-1}(K)=\phi^{-1}( K)\cap\cup_{k\in\N}M_k$ is compact and so is contained in $\cup_{i\leq k}\partial M_i$ for some $k\in\N$. Thus we have that $\phi^{-1}(K)\subset  M_k$ and so $\phi^{-1}(K)$ is a compact since it is closed.

\vspace{0.3cm}

\textbf{Step 3:} Up to picking a sub-sequence of the $\set{ M_i}_{i\in\N}$ and a proper homotopy of $\psi$ we can assume that all surfaces $\mathcal S\eqdef \cup_{k\in\N}M_k$ are properly embedded in $N$.

\vspace{0.3cm}

Since $\psi$ is a proper map when restricted to $\cup_{i\in\N} \partial M_i$ we have that for all $i$ we can find neighbourhoods $U_i\subset N$ of $\psi(\partial M_i)$ with compact closure such that the open sets $\set{U_i}_{i\in\N}\subset N$ are properly embedded. This means that up to picking a subsequence, which we still denote by $i\in\N$, we can assume that the $\set{U_i}_{i\in\N}\subset N$ are pairwise disjoint.

Then we have a $\pi_1$-injective map: $\psi: \partial M_i\rar U_i$. By the existence results for PL-least area surfaces in \cite{JaRu1988} and \cite[1.26]{Kap2001} for a triangulation $\tau_i$ of $N$ there are PL-least area representatives $S'$ of the $\psi(S)$, $S\in\pi_0(\partial M_i)$, such that $S'\subset U_i$. Now consider the cover $N_j$ of $N$ corresponding to $\pi_1( M_j)$ such that the set $U_i$ lifts homeomorphically to $\tilde U_i$ in the cover and denote by $\tilde S'$ the lift of $S'$. By Lemma \ref{finiteapprox}, up to homotopy we have that $\tilde\psi: M_j\rar N_j$ is an embedding, thus we see that all the $\tilde S'$ are homotopic to pairwise disjoint embedded surfaces $\Sigma$. Moreover, since by properties of a minimal exhaustion we have that $[\Sigma]\neq 0$ in $H_2(M)$ by the results of \cite[Thm 6]{JaRu1988} the $S'$'s are embedded as well. By properties of a minimal exhaustion we have that no two $\tilde S'$ are covering of an embedded surface thus by \cite[Thm 7]{JaRu1988} the  PL-least area surfaces are disjoint. Since the covering projection is a homeomorphism on $\pi_j:\tilde U_i\rar U_i$ we get that $S'=\pi_i(\tilde S')$ are embedded as well. By repeating this for all $i\in\N$ we obtain a proper homotopy of $\psi$ such that for $\set{ M_i}_{i\in\N}$ the restriction of $\psi$ to $\partial  M_i$ is an embedding. Moreover, since the $U_i$'s were pairwise disjoint we see that $\cup_{i\in\N}\partial  M_i$ actually embeds in $N$.

 The last claim in the statement, that $\psi$ is a proper embedding on tame ends of $M$ follows by the same argument of Lemma \ref{tameends}. \epf
 									
									We now want to promote $\phi$ to a homeomorphism from $M$ to $N$. This will complete the proof of the main theorem:
									
									\bthm \label{maintheorem2} Let $M\in\M^B$, then $M$ is homeomorphic to a complete hyperbolic 3-manifold if and only if the maximal bordification $\overline M$ does not admit any doubly peripheral cylinder.									\ethm

									We now prove the final part of Theorem \ref{maintheorem} which is that if the maximal bordification $\overline M$ of $M$ does not admit any doubly peripheral cylinders then $M$ is hyperbolizable. 
											
									\bthm\label{homeqhomeo}
									Let $M\in\M^B$ and $\phi:  M\rar N$ be a homotopy equivalence with $N$ a complete hyperbolic manifold. If $\overline M$ does not have any doubly peripheral annulus, then we have a homeomorphism $\psi:M\rar N$.
									\ethm
									
									\bpf
									By Lemma \ref{minimalexhaustion} let $\set{M_i}_{i\in\N}$ be a minimal exhaustion of $M$. By Theorem \ref{homotopyequivalence} we have a proper homotopy equivalence preserving $P$: 
																		$$\phi:M\rar N$$ 
that is an embedding on $\cup_{k\in\N}\partial  M_{i_k}$ and tame ends of $M$. The submanifold $P\subset \partial\overline M$ is a collection of annuli and tori that make $(\overline M,P)$ an infinite-type pared acylindrical 3-manifold. Thus, without loss of generality we can assume that the exhaustion of $M$ is the one given by $\set{\partial  M_{i_k}}_{k\in\N}$. Therefore, we have a proper homotopy equivalence $\phi:M\rar N$ respecting $P$ that is an embedding on the boundary components of a minimal exhaustion $\set{ M_i}_{i\in\N}$ and any tame end of $M$.
Since $\phi $ is a homotopy equivalence we get that $\forall i :\phi(\partial  M_i)$ bounds a 3-dimensional compact submanifold $K_i$ of $N$. We now want to show that the $K_i$ are nested. Since $\phi\vert_{\cup_{i\in\N}\partial M_i}$ is an embedding we have that for all $i\neq j$: $\phi(\partial M_i)\cap\phi(\partial M_j)=\emp$ we only need to show that $\phi(\partial M_{i+1})\not\subset K_i$.

\vspace{0.3cm}

\paragraph{Claim:} For all $ i: \phi(\partial  M_{i+1})\not\subset K_i$ and up to a proper homotopy $M_i$ embeds in $N$ with image $K_i$.

\vspace{0.3cm}
\bpfc
Assume we have some $i$ such that the above does not happen so that there is $S\in\pi_0(\partial  M_{i+1})$ such that $\phi(S)\subset K_i$. Pick $L>i$ so that $K_i$ lifts homeomorphically $\tilde\iota (K_i)\hookrightarrow N_L$ for $N_L$ the cover corresponding to $\phi_*(\pi_1(M_L))$. Then in the cover we see $\tilde\iota(K_i)$ and $\tilde\phi(S)$ inside it. By Lemma \ref{finiteapprox} we have that the map $\tilde\phi: M_L\rar \overline N_L$, for $\overline N_L$ the manifold compactification \cite{AG2004,CG2006} of $N_L$, is homotopic to a homeomorphism $\psi$. Particularly, we have $\psi(M_{i})\subset \overline N_L$ and, up to isotopy, we can assume $\psi(\partial M_i)=\tilde\iota(\partial K_i)$. Since $\overline N_L$ is not a closed 3-manifold we must have $\tilde\iota(K_i)= \psi( M_i)$. Moreover, since $\tilde\iota(K_i)$ projects down homeomorphically we get that up to a proper homotopy $\phi $ embeds $ M_i$ in $N$.

In particular from the homeomorphism $\psi: M_L\rar\overline N_L$ we see that $\tilde\phi(S)$ is homotopic outside $\tilde\iota(K_i)=\psi( M_i)$ and therefore it must be homotopic into $\tilde\iota(\partial K_i)=\psi(\partial M_i)$. Since we had a minimal exhaustion this can only happen if $S$ co-bounds with $S'\subset \partial M_i$ an $I$-bundle contained in a tame end of $M$. Therefore, since $\phi$ was a proper embedding on tame ends of $M$ we reach a contradiction.\epfc

Thus we can assume that we have an exhaustion $\set {K_i}_{i\in\N}$ of $N$ with $\phi(\partial  M_i)=\partial K_i$ and we define $K_{j,i}\eqdef \overline{K_j\setminus K_i}$. Moreover, by the claim we also have that $\pi_1(K_{j,i})\overset{\phi_*}\cong \pi_1(U_{j,i})$ for $U_{j,i}\eqdef \overline{M_j\setminus M_i}$ and $U_{1,0}=M_1$.

By the claim we can also assume that up to a proper homotopy of $\phi $ the restriction $\phi\vert_{M_1}$ is an embedding with image $ K_1$.

							To conclude the proof we need to show that the map $\phi$ is properly homotopic to an embedding.

We will now show with an inductive argument that up to proper homotopy $\phi$ is an embedding. Our base case is that $ M_1$ embeds. By an iterative argument we need to show that we can embed  $U_{i+1,i}$ relative to the previous embedding, hence rel $\partial  M_i$.

Consider the following diagram:
	\be \xymatrix{ & \quotient{\bH}{\rho_\infty(\pi_1(U_{i+1,i}))}\ar[dd]^\pi & \\ U_{i+1,i}\ar[ur]^{\tilde\phi} \ar[dr]_\phi & & K_{i+1,i}\ar[dl]^\iota\ar[ul]_{\tilde\iota}\\ & N & }
											\ee 
By Lemma \ref{2.2} we have that $\tilde\phi$ is homotopic to an embedding $\psi$ rel boundary and we have that $\psi(\partial M_i)=\psi(\partial K_i)$. Then we can isotope $\psi$ so that $\psi(\partial  M_{i+1})=\tilde\iota(\partial K_{i+1})$. Hence we have that $\psi(U_{i+1,i})=\tilde\iota(K_{i+1,i})$, they are compact submanifolds with the same boundary in an open manifold. Therefore we get that $\pi\circ\psi$ is properly homotopic to $\phi$, the homotopy is constant outside a compact set, and embeds $U_{i+1,i}$ rel the previous embedding. We can then glue all this proper homotopies together to get a proper embedding $\psi: M\hookrightarrow N$. Since the embedding is proper and $N$ is connected we get that $\psi$ is a homeomorphism from $M$ to $N$ completing the proof. \epf

								By combining Theorem \ref{nothyp} and Theorem \ref{homeqhomeo} we complete the proof of Theorem \ref{maintheorem2}.	
								
								\brem\label{catzero}
								Using our main result we can show that manifolds in $\M$ have $\cat{CAT}(0)$ metric in which $\pi_1(M)$ acts by semisimple isometries. Then by \cite[p.86]{BGS1985} for $\gamma\in\pi_1(M)$ we have that the centraliser $C(\gamma)$ is isomorphic to $\Z$. Since all roots of $\gamma$ are in $C(\gamma)$ we would get that $\Z$ has a divisible element which is impossible.	
																
								To construct the $\cat{CAT}(0)$ structure let $\A$ be the collection of doubly peripheral annuli in $\overline M$. Let $X\eqdef\coprod_{i\in\N}X_i$ be the manifold obtained by splitting $\overline M$ along the annuli $\A$. Each manifold $X_i$ has a collection $\A_i$ of annuli in $\partial X_i$ corresponding to annuli in $\A$. By Theorem \ref{maintheorem} we can construct a complete hyperbolic metric on $X_i$. Moreover, we can rig the hyperbolic metric so that all $\A_i$ correspond to rank one cusps\footnote{This is achieved by choosing the curves $\P_i$ that make $X_i$ acylindrical in $\partial X_i\setminus \A_i$ and then adding $\A_i$ to the collection $\P_i$. }. Then by flattening all the cusps we obtain complete $\cat{CAT}(0)$ metrics in every $X_i$.  Then by gluing back by euclidean isometries along the $\A_i$ one can obtain a singular $\cat{CAT}(0)$ metric on $M$ in which every element is represented by an hyperbolic isometry, since they all have an axis.  												\erem

\newpage

\appendix
\section{The manifold $N$ and $X$}\label{appendix A}

		A \emph{knot} is an embedding $K:\mathbb S^1\hookrightarrow M$. Given a knot we denote the complement of a regular neighbourhood $N_r(K)$ of $\im(K)$ in $M$ by $M_K\eqdef\overline{M\setminus N_r(K)}$. We recall the following theorem by Myers \cite[6.1]{My1982}:

			\begin{Theorem*} Let $M$ be a compact, orientable, 3-manifold whose boundary contains no 2-spheres. Then $M$ has a knot $K:\mathbb S^1\hookrightarrow M$ such that $M_K$ is irreducible, with incompressible boundary and without any non-boundary parallel annuli or tori.\end{Theorem*}
	 		
			\noindent We call such a knot a \textit{simple knot}.
			  
Consider the hyperbolizable 3-manifold $M\eqdef \Sigma_{1,1}\times I$. By Myer's Theorem we can pick a simple knot $K$ in $M$ such that $M_K$ is acylindrical, atoroidal, irreducible and with incompressible boundary $\partial M_k=\Sigma_2\coprod\mathbb T^2$. We now want to fill in the torus boundary $\mathbb T^2$ while keeping the resulting manifold acylindrical and hyperbolizable.

Consider the manifold $M'$ obtained by gluing two copies of $M_K$ along the genus two boundaries. The manifold $M'$ has two torus boundaries corresponding to the copies of the knot $K$ and since each $M_K$ is acylindrical and atoroidal we have that $M'$ is atoroidal as well. Thus, by the Hyperbolization Theorem we have that $M'$ is hyperbolizable. By Hyperbolic Dehn Filling Theory \cite{BP1992,Th1978} we can find a high enough Dehn Filling of type $\f pq$ so that the manifold $M'$ filled by $\f pq $ surgery on the tori is hyperbolic. We denote by $M'\left(\f p q, \f pq \right)$ the resulting manifold.

The manifold $M'\left(\f pq,\f pq\right)$ is homeomorphic to the double of $M_{K}$ along the genus two boundary with $\f p q$ filling in the tori boundaries. The double being atoroidal implies that $M_{K}\left(\f p q\right)$ is acylindrical. Therefore, by doing $\f p q$ Dehn filling on $M_K$ we obtain an acylindrical and hyperbolizable 3-manifold $N$ that has for boundary an incompressible genus two surface.

By gluing two copies of $N$ along a separating annulus $A$ in their boundary we get a manifold $X$ such that $\partial X$ are two genus two surfaces. We also denote by $A$ the essential annulus in $X$ obtained by the gluing. We now claim that $X$ is hyperbolizable, since $\pi_1(X)$ is infinite by the Hyperbolization Theorem \cite{Kap2001} it suffices to show that $X$ is atoroidal. Let $T\subset X$ be an essential torus. Then, up to isotopy, $T\vert A\eqdef \overline{T\setminus A}$ is a collection of embedded cylinders in $N$ that are $\pi_1$-injective. Therefore, since $N$ is acylindrical we have that all the components of $T\vert A$ are boundary parallel. Moreover, since their boundaries are contained in $A$ all components of $T\vert A$ are isotopic into $A$. Therefore, the torus $T$ is not essential since $\iota(\pi_1(T))\cong\Z$.

Hence, $X$ is an hyperbolizable 3-manifold with a unique essential cylinder $A$ and two genus two incompressible boundaries.

\section{Example with Divisible Element}\label{appendix B}
Recall:
\bdefi\label{divisiblelementdefi}
An element $\gamma\in G$ is said to be \emph{divisible} if for all $n\in\N$ there is $a\in G$ such that $\gamma=a^n$.
\edefi

			Let $N$ be the acylindrical hyperbolizable 3-manifold with genus two incompressible boundary constructed in Appendix \ref{appendix A}. Consider an infinite one sided \textit{thick} cylinder $C\cong(\mathbb S^1\times [0,1])\times  [0,\infty)$ and let $\gamma$ be the generator of $\pi_1(C)$.
			
			Let $\set{T_n}_{n=2}^\infty$ be a collection of solid tori such that $T_n$ has one wing winding around the soul $n$ times. Glue the boundary of the wing of $T_n$ to $C$ along a small neighbourhood of $\mathbb S^1\times \set 1\times \set {n-\f 32}$. The resulting 3-manifold $\widehat C$ has a divisible element given by $\gamma$ and is not atoroidal since if we consider the portion containing $T_n,T_{n+1}$ it contains an incompressible non-boundary parallel torus.  
			
			By construction we still have that: $\partial \widehat C\cong\mathbb  S^1\times (-\infty,\infty)$ and we can think of the boundary of the solid torus as being a neighbourhood of $\mathbb S^1\times \set 1\times \set {n-\f 32}$ so that the $\mathbb S^1\times \set 1\times[n,n+1]$ pieces now contain part of the boundary of the solid tori. 
			
			On each $N$ we mark a closed neighbourhood $A$ of the simple closed curve in $\partial N$ splitting the genus two surface into two punctured tori $\Sigma^\pm$. We then glue countably many copies $\set{N_n}_{n\in\N}$ of $N$ such that the marked annulus $A_n$ in $\partial N_n$ is glued to $\mathbb S^1\times \set 1 \times [n,n+1]$ and then glue $\Sigma_n^+$ to $\Sigma_{n+1}^-$ via the identity. On $\mathbb S^1\times\set 0\times[0,\infty)$ we glue countably many copies $\set{N_k}_{k\in\N}$ of $N$ by gluing $A_k$ in $\partial N_k$ to $\mathbb S^1\times \set 1 \times [k,k+1]$ and $\Sigma_k^+$ to $\Sigma_{k+1}^-$ via the identity.

\begin{center}\begin{figure}[h!]
	 										\centering
	 										\def\svgwidth{300pt}
										
	 										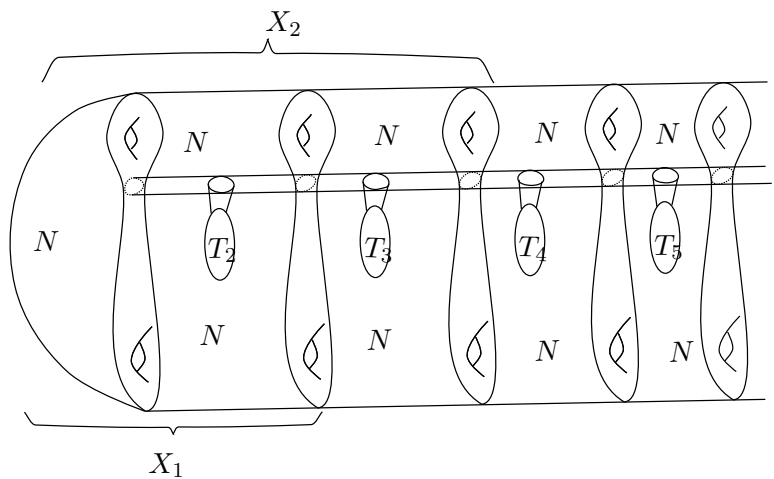	
											\caption{The manifold $X$ with the first two elements of the exhaustion.}\label{fig19}
\end{figure}\end{center}

			Finally we glue another $N$ to the remaining genus $2$ boundary component. The result is a 3-manifold $X$ that has an exhaustion $X_i$, $i\in\N$, given by taking the manifolds $\set{N_k}_{k\leq i},\set{N_n}_{n\leq i}$, $\mathbb S^1\times I\times [0,n]$ and the bottom copy of $N$.  A schematic of the manifold is given in Figure \ref{fig19}.

The gaps $\overline{X_n\setminus X_{n-1}}$ are hyperbolizable since they are homeomorphic to two copies of $N$ glued along the solid torus $T_{n+1}$ (the proof of the manifold being atoroidal is similar to the proof of the atoroidality of $X$). The element of the exhaustion are not hyperbolizable, for example if we look at $X_2$ we see that it has a Torus subgroup $\langle \alpha,\beta\vert \alpha^2=\beta^3\rangle$. 

\nocite{BP1992,CEM2006,CM2006,Ha2002,He1976,Sh1975,Sc1972,SY2013,Th1978,Ja1980,MT1998,Ca1993,MT1998,So2006,Th1982}
			
\thispagestyle{empty}
{\small
\markboth{References}{References}
\bibliographystyle{alpha}
\bibliography{mybib}{}
}

\bigskip

\noindent Department of Mathematics, Boston College.

\noindent 140 Commonwealth Avenue Chestnut Hill, MA 02467.

\noindent Maloney Hall
\newline \noindent
email: \texttt{cremasch@bc.edu}

\end{document}